\documentclass[10pt]{amsart}
\usepackage{latexsym,amsmath,amssymb,amscd}
\def\today{\ifcase \month \or
   January \or February \or March \or April \or
   May \or June \or July \or August \or
   September \or October \or November \or December \fi
   \space\number\day , \number\year}

\begin{document}

\makeatletter
\@addtoreset{figure}{section}
\def\thefigure{\thesection.\@arabic\c@figure}
\def\fps@figure{h,t}
\@addtoreset{table}{bsection}

\def\thetable{\thesection.\@arabic\c@table}
\def\fps@table{h, t}
\@addtoreset{equation}{section}
\def\theequation{
\arabic{equation}}
\makeatother

\newcommand{\bfi}{\bfseries\itshape}

\newtheorem{theorem}{Theorem}
\newtheorem{acknowledgment}[theorem]{Acknowledgment}
\newtheorem{algorithm}[theorem]{Algorithm}
\newtheorem{axiom}[theorem]{Axiom}
\newtheorem{case}[theorem]{Case}
\newtheorem{claim}[theorem]{Claim}
\newtheorem{conclusion}[theorem]{Conclusion}
\newtheorem{condition}[theorem]{Condition}
\newtheorem{conjecture}[theorem]{Conjecture}
\newtheorem{construction}[theorem]{Construction}
\newtheorem{corollary}[theorem]{Corollary}
\newtheorem{criterion}[theorem]{Criterion}
\newtheorem{data}[theorem]{Data}
\newtheorem{definition}[theorem]{Definition}
\newtheorem{example}[theorem]{Example}
\newtheorem{lemma}[theorem]{Lemma}
\newtheorem{notation}[theorem]{Notation}
\newtheorem{problem}[theorem]{Problem}
\newtheorem{proposition}[theorem]{Proposition}
\newtheorem{question}[theorem]{Question}
\newtheorem{remark}[theorem]{Remark}
\newtheorem{setting}[theorem]{Setting}
\numberwithin{theorem}{section}
\numberwithin{equation}{section}

\newcommand{\todo}[1]{\vspace{5 mm}\par \noindent
\framebox{\begin{minipage}[c]{0.85 \textwidth}
\tt #1 \end{minipage}}\vspace{5 mm}\par}

\newcommand{\1}{{\bf 1}}

\newcommand{\hotimes}{\widehat\otimes}

\newcommand{\Ad}{{\rm Ad}}
\newcommand{\ad}{{\rm ad}}
\newcommand{\Alt}{{\rm Alt}\,}
\newcommand{\Aut}{{\rm Aut}}
\newcommand{\Ci}{{\mathcal C}^\infty}
\newcommand{\Co}{{\mathcal C}^\omega}
\newcommand{\comp}{\circ}
\newcommand{\C}{\mathbb C} 
\newcommand{\CB}{{\mathcal C}{\mathcal B}}
\newcommand{\D}{\text{\bf D}}
\newcommand{\de}{{\rm d}}
\newcommand{\ee}{{\rm e}}
\newcommand{\ev}{{\rm ev}}
\newcommand{\fimes}{\mathop{\times}\limits}
\newcommand{\id}{{\rm id}}
\newcommand{\ie}{{\rm i}}
\newcommand{\End}{{\rm End}\,}
\newcommand{\Gr}{{\rm Gr}}
\newcommand{\GL}{{\rm GL}}
\newcommand{\Hilb}{{\bf Hilb}\,}
\newcommand{\Hom}{{\rm Hom}\,}
\newcommand{\Iso}{{\rm Iso}}
\newcommand{\Ker}{{\rm Ker}\,}
\newcommand{\Kern}{\textbf{Kern}}
\newcommand{\Lie}{\textbf{L}}
\newcommand{\lf}{{\rm l}}
\newcommand{\N}{\mathbb N}
\newcommand{\pr}{{\rm pr}}
\newcommand{\R}{\mathbb R} 
\newcommand{\Ran}{{\rm Ran}\,}
\newcommand{\Rep}{{\rm Rep}}
\newcommand{\RK}{{\mathcal P}{\mathcal K}^{-*}}
\newcommand{\spann}{{\rm span}}
\newcommand{\T}{\text{\bf T}}
\newcommand{\Tr}{{\rm Tr}\,}

\newcommand{\G}{{\rm G}}
\newcommand{\U}{{\rm U}}
\newcommand{\VB}{{\rm VB}}

\newcommand{\Bc}{{\mathcal B}}
\newcommand{\Cc}{{\mathcal C}}
\newcommand{\Dc}{{\mathcal D}}
\newcommand{\Ec}{{\mathcal E}}
\newcommand{\Fc}{{\mathcal F}}
\newcommand{\Gc}{{\mathcal G}}
\newcommand{\Hc}{{\mathcal H}}
\newcommand{\Kc}{{\mathcal K}}
\newcommand{\Lc}{{\mathcal L}}
\newcommand{\Nc}{{\mathcal N}}
\newcommand{\Oc}{{\mathcal O}}
\newcommand{\Pc}{{\mathcal P}}
\newcommand{\Qc}{{\mathcal Q}}
\newcommand{\Sc}{{\mathcal S}}
\newcommand{\Tc}{{\mathcal T}}
\newcommand{\Uc}{{\mathcal U}}
\newcommand{\Vc}{{\mathcal V}}
\newcommand{\Xc}{{\mathcal X}}
\newcommand{\Yc}{{\mathcal Y}}
\newcommand{\Zc}{{\mathcal Z}}

\newcommand{\Ag}{{\mathfrak A}}
\newcommand{\Bg}{{\mathfrak B}}
\renewcommand{\gg}{{\mathfrak g}}
\newcommand{\hg}{{\mathfrak h}}
\newcommand{\Hg}{{\mathfrak H}}
\newcommand{\mg}{{\mathfrak m}}
\newcommand{\Mg}{{\mathfrak M}}
\newcommand{\nng}{{\mathfrak n}}
\newcommand{\pg}{{\mathfrak p}}
\newcommand{\Sg}{{\mathfrak S}}
\newcommand{\ug}{{\mathfrak u}}
\newcommand{\Ug}{{\mathfrak U}}
\newcommand{\Xg}{{\mathfrak X}}

\pagestyle{myheadings}
\markboth{}{}


\makeatletter
\title[Holomorphic geometric models for representations of $C^*$-algebras]{Holomorphic geometric models for representations of $C^*$-algebras}
\author{Daniel Belti\c t\u a} 
\author{Jos\'e E. Gal\'e}
\address{Institute of Mathematics ``Simion
Stoilow'' of the Romanian Academy, 
RO-014700 Bucharest, Romania}
\email{Daniel.Beltita@imar.ro}
\address{Departamento de matem\'aticas, Facultad de Ciencias, 
Universidad de Zaragoza, 50009 Zaragoza, Spain}
\email{gale@unizar.es}

\keywords{representation; Banach-Lie group; $C^*$-algebra; conditional expectation; 
homogeneous vector bundle; 
reproducing kernel; Stinespring dilation; amenable Banach algebra}

\subjclass[2000]{Primary 46L05; 
      Secondary 46E22; 47B38; 46L07; 46L55; 58B12; 43A85; 22E65}

\begin{abstract}
Representations of $C^*$-algebras are realized on section spaces 
of holomorphic homogeneous vector bundles. 
The corresponding section spaces are investigated by means of a new 
notion of reproducing kernel, suitable for dealing with involutive diffeomorphisms 
defined on the base spaces of the bundles.
Applications of this technique to dilation theory of completely positive maps 
are explored and the critical role of complexified homogeneous spaces 
in connection with the Stinespring dilations is pointed out. 
The general results are further illustrated by a discussion 
of several specific topics, including similarity orbits of 
representations of amenable Banach algebras, 
similarity orbits of conditional expectations, 
geometric models of representations of Cuntz algebras, 
the relationship to endomorphisms of $\Bc(\Hc)$, 
and non-commutative stochastic analysis. 
\end{abstract}

\date{\tt July 5, 2007}
\makeatother
\maketitle

\section{Introduction}\label{intro}

Originally, the interest in the study of representations of algebras and groups of operators
on infinite-dimensional Hilbert or Banach spaces is to be found, as one of the main  
motivations, in problems arising from Quantum Physics. In this setting, 
unitary groups of operators can be interpreted as symmetry groups while  
the self-adjoint operators are thought of as observable objects,  
hence the direct approach to such questions leads naturally to 
representations both involving algebras generated by commutative or non-commutative
canonical relations, and groups of unitaries on Hilbert spaces; 
see for instance \cite{GW54}, \cite{Sh62}, or \cite{Se57}. 
Over the years, there have been important developments of this initial approach, 
in papers devoted to analyze or classify a wide variety of representations, 
and yet many questions remain open in the subject. 
It is certainly desirable to transfer to this field methods, 
or at least ideas, of the rich representation theory of finite-dimensional Lie groups.

In this respect, recall that geometric representation theory is a classical topic 
in finite dimensions. Its purpose is to shed light on certain classes of representations 
by means of their geometric realizations  
(see for instance \cite{Ne00}). 
Thus the construction of geometric models of representations 
lies at the heart of that topic, and one of the classical 
results obtained in this direction is the Bott-Borel-Weil theorem 
concerning realizations of irreducible representations of compact Lie groups 
in spaces of sections (or higher cohomology groups) 
of holomorphic vector bundles over flag manifolds; 
see \cite{Bo57}. 
Section spaces of vector bundles also appear in methods of induction of representations, of Lie groups, 
from representations initially defined on appropriate subgroups. 
Induced representations are required
for instance by the so-called orbit method, consisting of establishing a neat link 
between general representations of a Lie group and the symplectic geometry of its coadjoint orbits; 
see \cite{Ki04} or \cite{Fo95}. 
Such sections are quite often obtained out of 
suitable square-integrable functions on the base space of the bundle. 

Sometimes these ideas work well in the setting of infinite-dimensional Lie groups, in special situations 
or for particular aims; see for example \cite{Bo80}, \cite{Ki04}, and \cite{Ne04}.
However, several difficult points are encountered when one tries 
to extend these ideas in general, and perhaps the most important one is related to 
the lack of an algebraic structure theory for representations of these groups. 
Also, it is not a minor question the fact that, in infinite
dimensions, there is no sufficiently well-suited theory of integration. 
The most reasonable way to deal with these problems seems to be to restrict both the class of groups 
and the class of representations one is working with. 
Moreover, one is led quite frequently to employ methods of operator algebras. 
See for example \cite{SV02}, where the study of factor representations of 
the group $\U(\infty)$ and AF-algebras is undertaken. 
There, a key role is played by the 
Gelfand-Naimark-Segal (or GNS, for short) representations constructed out of 
states of suitable maximal abelian self-adjoint subalgebras.

The importance of GNS representations as well as that of the geometric properties of state spaces 
in operator theory are well known. 
In \cite{BR07}, geometric realizations of restrictions of GNS representations to 
groups of unitaries in $C^*$-algebras are investigated. 
Suitable versions of reproducing kernels 
on vector bundles are considered, in order to build representation spaces formed by sections. 
This technique has already a well-established place in 
representation theory of finite-dimensional Lie groups; 
see for instance the monograph \cite{Ne00}. 
In some more detail, let $\1\in B\subseteq A$ be unital $C^*$-algebras such that there exists 
a conditional expectation $E\colon A\to B$. 
Let  $\U_A$ and $\U_B$ be the unitary groups of $A$ and $B$ respectively,  
and $\varphi$ a state of $A$ such that $\varphi\circ E=\varphi$. 
A reproducing kernel Hilbert 
space $\Hc_{\varphi,E}$ can be constructed out of $\varphi$ and $E$, 
consisting of $\Ci$ sections 
of a certain Hermitian vector bundle with base the homogeneous space $\U_A/\U_B$, 
and the restriction to $\U_A$ of the GNS representation associated 
with $\varphi$ can be realized (by means of a certain intertwining operator) as the natural multiplication of 
$\U_A$ on $\Hc_{\varphi,E}$; see Theorem~5.4 in \cite{BR07}.   
This theorem relates the GNS representations to the geometric representation theory, 
in the spirit of the Bott-Borel-Weil theorem.  
In view of this result and of the powerful method of induction developed in \cite{Bo80}, 
it is most natural to ask about similar results
for more general representations of infinite-dimensional Lie groups. 

Another circle of ideas is connected with holomorphy. 
Recall that this is the classical setting of the Bott-Borel-Weil theorem of \cite{Bo57} 
involving the flag manifolds,  
and it reinforces the strength of applications. 
On the other hand, the idea of \emph{complexification} plays a central role in this area, 
inasmuch as one of the ways to describe the complex structure of the flag manifolds is 
to view the latter as homogeneous spaces of complexifications of compact Lie groups. 
(See \cite{LS91}, \cite{Bi03}, \cite{Bi04}, and \cite{Sz04} 
for recent advances in understanding the differential geometric flavor of 
the process of complexification.) 
In some cases involving finiteness properties of spectra and traces of elements in a $C^*$-algebra, 
it is possible to prove that the aforementioned infinite-dimensional homogeneous space $\U_A/\U_B$ 
is  
a \emph{complex} manifold as well and the Hilbert space $\Hc_{\varphi,E}$ is formed by 
\emph{holomorphic} sections 
(see Theorem 5.8 in \cite{BR07}). 

One can find in Section~\ref{sect1} of the present paper some related results of holomorphy   
in the important special case of tautological bundles over Grassmann manifolds associated 
with involutive algebras. 
For the reader's convenience, these results are exposed in some detail  
since they illustrate the main ideas underlying the present investigation. 
(A complementary perspective on these manifolds 
can be found in \cite{BN05}). 
Apart from the above two examples, the holomorphic character of the manifolds
$\U_A/\U_B$ (and associated bundles) is far from being clear in general.
On the other hand, the aforementioned conditional expectation 
$E\colon A\to B$ has a strong geometric
meaning as a connection form defining a reductive structure in the homogeneous space 
$\G_A/\G_B$; see \cite{ACS95} and \cite{CG99}. 
Since $X$ is the Lie algebra of the complex Banach-Lie group $\G_X$ for $X=A$ and $X=B$, 
it is important to incorporate full groups of invertibles to 
the framework established in \cite{BR07}. 
Note also that $\G_X$ is the universal complexification of $\U_X$, 
according to the discussion of \cite{Ne02}.  

\textbf{Brief description of the present paper.}
One of our aims in the present paper is to extend the geometric representation theory of 
unitary groups of operator algebras to the complex setting of full groups of invertible elements. 
For this purpose we need a method to realize the representation spaces as 
Hilbert spaces of sections in holomorphic vector bundles. 
If one tries to mimic the arguments of \cite{BR07} then 
one runs into troubles very soon  
(regarding the construction of appropriate reproducing kernels), 
due to the fact that general invertible elements of a $C^*$-algebra lack, when
considered in an inner product, helpful cancellative properties (that unitaries have). 
This can be overcome by 
using certain involutions $z\mapsto z^{-*}$ (that come from the involutions of $C^*$-algebras) 
on the bases of the
bundles, but then the problem is that our bundles lose their Hermitian character.  

So we are naturally led toward developing a special theory of reproducing kernels on vector bundles. 
Section~\ref{sect2} includes a discussion of a version of Hermitian vector bundles 
suitable for our purposes. 
We call them \emph{like-Hermitian}. 
The bases of such vector bundles 
are equipped with involutive diffeomorphisms $z\mapsto z^{-*}$, so that 
we need to find out a class of reproducing kernels, compatible in 
a suitable sense with the corresponding diffeomorphisms, 
which we call here reproducing $(-*)$-kernels. 
The very basic elements for the theory of reproducing $(-*)$-kernels 
are presented in Section~\ref{sect3} (it is our intention to develop such a theory more
sistematically in forthcoming papers). 
In Section~\ref{sect4} we discuss examples of the above notions 
which arise in relation to homogeneous manifolds obtained by (smooth) actions of complex Banach-Lie groups (see
Definition~\ref{homog}). 
These examples play a critical role for our main constructions of geometric models of representations; 
see Theorem~\ref{induction} and Theorem~\ref{realization_hom}. 
In particular, Theorem~\ref{realization_hom} 
provides the holomorphic versions of such realizations. In order to include the homogeneous spaces of unitary
groups $\U_A/\U_B$ in the theory and to avoid the fact that they are not necessarily complex manifolds, we had to
view them as embedded into their natural complexifications $\G_A/\G_B$. 

It is remarkable that, using a significant polar
decomposition of $\G_A$ found by Porta and Rech, relative to a prescribed conditional expectation 
(see \cite{PR94}), it is
possible to interpret the manifold $\G_A/\G_B$ as (diffeomorphic to) the tangent bundle of $\U_A/\U_B$, see Theorem~\ref{suppl1} and Theorem~\ref{supplcor2} below. 
These properties resemble very much similar properties
enjoyed by complexifications of manifolds of compact type in finite dimensions. 
This may well mean that the homogeneous spaces 
$\U_A/\U_B$ and $\G_A/\G_B$ are good substitutes for compact homogeneous spaces 
in the infinite-dimensional setting. 

The set of ideas previously exposed can be used to investigate geometric models for representations
which arise as Stinespring dilations of completely positive maps on $C^*$-algebras $A$. 
In this way we shall actually end up with a geometric dilation theory 
of completely positive maps. 
This in particular enables us to get more examples of 
representations of Banach-Lie groups (namely, $\U_A$, $\G_A$) which admit geometric realizations 
in the sense of \cite{BR07}. 
Also, just by differentiating it is possible to recover the whole dilation on $A$
and not only its restriction to $\U_A$ or $\G_A$, see Theorem~\ref{stinerealiz}, 
and this provides a geometric interpretation of the classical methods of extension and induction of representations of $C^*$-algebras 
(see \cite{Di64} and \cite{Ri74}). 
We should point out here that there exist earlier approaches 
to questions in dilation theory with a geometric flavor 
---see for instance \cite{ALRS97}, \cite{Ar00}, \cite{Po01}, or \cite{MS03}--- 
however they are different from the present line of investigation. 

The last section of the paper, Section~\ref{sect6}, is devoted to showing, by means of 
several specific examples, that the theory established here has interesting links with 
quite a number of different subjects in operator theory and related areas.

For the sake of better explanation, we conclude this introduction by 
a summary of the main points considered 
in the paper. 
These are: 
\begin{enumerate}
\item[-] a theory of reproducing kernels on vector bundles 
that takes into account prescribed involutions of the bundle bases 
(Section~\ref{sect3}); 
\item[-] in the case of homogeneous vector bundles we investigate 
a circle of ideas centered on the relationship between reproducing kernels 
and complexifications of homogeneous spaces 
(Theorems \ref{realization_hom}~and~\ref{suppl1}); 
\item[-] by using the previous items we model the representation spaces of 
Stinespring dilations as spaces of holomorphic sections 
in certain homogeneous vector bundles; 
thereby we set forth a rich panel of differential geometric structures 
accompanying the dilations of completely positive maps 
(Section~\ref{sect5}); 
for one thing, we provide a geometric perspective 
on induced representations of $C^*$-algebras (cf. \cite{Ri74}); 
\item[-] as an illustration of our results we describe in Section~\ref{sect6} 
a number of geometric properties of orbits of representations of nuclear $C^*$-algebras 
and injective von Neumann algebras (Corollary~\ref{simiunita}), similarity orbits of conditional
expectations, and some relationships with representations of Cuntz algebras and endomorphisms of $\Bc(\Hc)$, 
as well as with non-commutative stochastic analysis.  
\end{enumerate}

\section{Grassmannians and homogeneous Hermitian vector bundles}\label{sect1}

We begin with several elementary considerations about 
idempotents in complex associative algebras.

\begin{notation}\label{set}
\normalfont
We are going to use the following notation: 
$A$ is a unital associative algebra over ${\mathbb C}$ with 
unit $\1$ and set of idempotents 
$\Pc(A)=\{p\in A\mid p^2=p\}$; 
for $p_1,p_2\in\Pc$ the notation $p_1\sim p_2$ means that 
we have both $p_1p_2=p_2$ and
$p_2p_1=p_1$.  For each $p\in\Pc(A)$ we denote its equivalence class by 
$[p]:=\{q\in\Pc(A)\mid q\sim p\}$. 
The quotient set is denoted 
by $\Gr(A)=\Pc(A)/ \sim$ (the \textit{Grassmannian} of $A$) 
and the quotient map by $\pi\colon\,p\mapsto[p]$, $\Pc\to\Gr(A)$.  

The group of invertible elements of $A$ is denoted by $\G_A$,  
and it has a natural action on $\Pc(A)$ by 
$$
\alpha\colon\,(u,q)\mapsto uqu^{-1},\quad \G_A\times\Pc(A)\to\Pc(A).
$$
The corresponding isotropy group at $p\in\Pc(A)$ is 
$\{u\in \G_A\mid\alpha(u,p)=p\}=\G_A\cap\{p\}'=\G_{\{p\}'}=:\G(p)$ 
where we denote by $\{p\}'$ 
the commutant
subalgebra of $p$ in $A$ (see page~484 in \cite{DG02}).
\qed
\end{notation}

\begin{lemma}\label{beta}
There exists a well-defined action 
of the group $\G_A$ upon $\Gr(A)$ like this: 
$$
\beta\colon\,(u,[p])\mapsto[upu^{-1}],\quad \G_A\times\Gr(A)\to\Gr(A),
$$
and the diagram 
$$
\begin{CD} 
\G_A\times\Pc(A) @>{\alpha}>> \Pc(A) \\
@V{\id_\G\times\pi}VV @VV{\pi}V \\
\G_A\times\Gr(A) @>{\beta}>> \Gr(A)
\end{CD}
$$
is commutative. 
\end{lemma}

\begin{proof} 
See for instance the end of Section~3 in \cite{DG01}. 
\end{proof}

\begin{definition}\label{iso}
\normalfont
For every idempotent $p\in\Pc(A)$ we denote by $\G_A([p])$ 
the isotropy group of the action 
$\beta\colon \G_A\times\Gr(A)\to\Gr(A)$ at the point 
$[p]\in\Gr(A)$, 
that is, 
$\G_A([p])=\{u\in \G_A\mid[upu^{-1}]=[p]\}$. 
\qed
\end{definition}

The following statement concerns the relationship between the isotropy groups 
of the actions $\alpha$ and $\beta$ of $\G_A$ upon 
$\Pc(A)$ and $\Gr(A)$, respectively. 

\begin{proposition}\label{isotropy} 
The following assertions hold. 
\begin{itemize}
\item[{\rm(i)}] For every $p\in\Pc(A)$ we have 
$\G_A([p])\cap \G_A([\1-p])=\G(p)$.
\item[{\rm(ii)}] If $\U$ is a subgroup of $G_A$ and $p\in\Pc(A)$ is 
such that $\U\cap \G_A([p])=\U\cap \G_A([\1-p])$, then 
$\U\cap \G_A([p])=\U\cap\{p\}'=:\U(p)$. 
\end{itemize}
\end{proposition}

\begin{proof}
(i) We have 
$$
\G_A([p])=\{u\in \G_A\mid[upu^{-1}]=[p]\}\text{ and }
\G_A([\1-p])=\{u\in \G_A\mid[u(\1-p)u^{-1}]=[\1-p]\},
$$
so that clearly $\G_A([p])\cap \G_A([\1-p])\supseteq \G_A\cap\{p\}'$. 
For the converse inclusion let $u\in \G_A([p])\cap \G_A([\1-p])$ arbitrary. 
In particular $u\in \G_A([p])$, whence $upu^{-1}\sim p$, 
which is equivalent to 
the fact that $(upu^{-1})p=p$ and $p(upu^{-1})=upu^{-1}$. 
Consequently we have both 
\begin{equation}\label{1}
pu^{-1}p=u^{-1}p
\end{equation}
and 
\begin{equation}\label{2b}
pup=up. 
\end{equation}
On the other hand, since $u\in \G_A([\1-p])$ as well, 
it follows that 
$(\1-p)u^{-1}(\1-p)=u^{-1}(\1-p)$ 
and $(\1-p)u(\1-p)=u(\1-p)$. 
The later equation is equivalent to 
$u-up-pu+pup=u-up$, 
that is, $pup=pu$. 
Then \eqref{2b} implies that $up=pu$, that is, $u\in \G(p)$.

(ii) This follows at once from part (i). 
\end{proof}

\begin{remark}\label{fixed}
\normalfont
For instance, Proposition~\ref{isotropy}(ii) can be applied if the algebra $A$ 
is equipped with an involution $a\mapsto a^*$ such that $p=p^*$, 
and $\U=\U_A:=\{u\in \G_A\mid u^{-1}=u^*\}$ is the corresponding unitary group. 
In this case, it follows by \eqref{1} and \eqref{2b} 
that $up=pu$ whenever $u\in \U_A\cap \G_A([p])$, 
hence $\U_A\cap \G_A([p])=\U_A\cap \G_A([\1-p])=\U_A\cap\{p\}'=:\U_A(p)$.
\end{remark} 

For $q\in\Pc(A)$, put $\hat q:=\1-q$ and $A^q:=\{a\in A\mid \hat qaq=0\}$. 
The following result is partly a counterpart, for algebras, 
of Proposition~\ref{isotropy}.

\begin{proposition}\label{almost} 
Assume that $A$ is equipped with an involution and 
let $p\in\Pc(A)$ such that $p=p^*$. 
Then the following assertions hold: 
\begin{itemize}
\item[{\rm(i)}] $uA^pu^{-1}= A^p$, for every $u\in U_A(p)$ ;
\item[{\rm(ii)}] $A^p\cap A^{\hat p}=\{p\}'$ ;
\item[{\rm(iii)}] $A^p+A^{\hat p}=A$;
\item[{\rm(iv)}] $(A^p)^*=A^{\hat p}$.
\end{itemize}
\end{proposition}

\begin{proof} (i) This is readily seen.

(ii) Firstly, note that, for $a\in A$, we have $\hat pap=pa\hat p$ 
if and only if $ap=pa$.
Moreover, if $ap=pa$ then $\hat pap=ap-pap=ap-ap=0$ and 
analogously $pa\hat p=0$. 
From this, the equality of the statement follows.

(iii) For every $a\in A$ and $q\in\Pc(A)$ we have $qa\in A^{q}$. 
Hence $a=pa+\hat pa\in A^p+A^{\hat p}$, as we wanted to show.

(iv) Take $a\in A^p$. Then $pap=ap$, that is, $pa^*p=pa^*$. Hence, 
$\hat p a^*\hat p=(a^*-pa^*)(\1-p)=a^*-pa^*-a^*p+a^*p=a^*-a^*p=a^*\hat p$. 
This means that $a^*\in A^{\hat p}$. Conversely, if $a\in A^{\hat p}$ then, as above, 
$pa^*p=a^*p$; that is, $a=(a^*)^*$ with $a^*\in A^p$. 
\end{proof}

Assume from now on that $A$ is a unital $C^*$-algebra. 
Then $\G_A$ is a Banach-Lie group whose Lie algebra coincides with $A$. 
The $\G_A$-orbits in $\Gr(A)$, obtained by
the action $\beta$ and equipped  with the topology inherited from 
$\Gr(A)$, are holomorphic Banach manifolds diffeomorphic to $\G_A/G_A([p])$ 
(endowed with its quotient topology), see Theorem~2.2 in \cite{DG02}. 
Also, the Grassmannian
$\Gr(A)$ can be described as  the discrete union of these $\G_A$-orbits, 
see \cite{DG01} and Theorem~2.3 in
\cite{DG02}. 
Moreover, $\U_A$ is a Banach-Lie subgroup of $\G_A$ with 
the Lie algebra $\ug_A:=\{a\in A\mid a^*=-a\}$. 
As it is well known, 
the complexification of $\ug_A$ is $A$, via the decomposition
$a=\{(a-a^*)/2\}+\ie\{(a+a^*)/2\ie\}$, ($a\in A$). 
Thus the conjugation of $A$ is given by 
$a\mapsto\overline a:=\{(a-a^*)/2\}-\ie\{(a+a^*)/2\ie\}=-a^*$.
We seek for possible topological and/or differentiable relationships
between the $\G_A$-orbits and the $\U_A$-orbits $\U_A/\U_A(p)$ in $\Gr(A)$. 

Let $p=p^*\in\Pc(A)$ and $\ug_A(p):=\ug_A\cap\{p\}'$. 
It is clear that $\ug_A(p)+\ie\ug_A(p)=\{p\}'$. 
Also, there is a natural identification between 
$\ug_A/\ug_A(p)$ and the tangent space $T_{[p]}(\U_A/\U_A(p))$.  

The above observations and Proposition~\ref{almost} yield
immediately the following result. 
Let $\Ad_\U$ denote the adjoint representation of $\U_A$.
 
\begin{proposition}\label{complexes} 
With the above notations,
$$
\Ad_\U(u)A^p\subset A^p,\ (u\in \U_A(p));
\ A^p\cap\overline{A^{\hat p}}=\ug_A(p)+\ie\ug_A(p);
\ A^p+\overline{A^{\hat p}}=A.
$$
In particular, $\ug_A/\ug_A(p)\simeq A/A^p$ whence 
we obtain that $\U_A/\U_A(p)$ and
$\G_A/\G_A([p])$ are locally diffeomorphic, 
and so $\U_A/\U_A(p)$ inherits the complex 
structure induced by $\G(A)/\G_A([p])$.
\end{proposition}

\begin{proof} 
The first part of the statement is just a rewriting of
Proposition~\ref{almost}. 
Then the result follows from Theorem~6.1 in \cite{Be06}.
\end{proof}

\begin{remark}\label{picture}
\normalfont
Since $\G_A(p)\subset \G_A([p])$, there exists the canonical projection 
$\G_A/\G_A(p)\to \G_A/\G_A([p])$. 
It is clear that its restriction to  
$\U_A/\U_A(p)$ becomes the identity map $\U_A/\U_A(p)\to \U_A/\U_A([p])$. 
We have seen that 
$\U_A/\U_A(p)$ enjoys a holomorphic structure 
inherited from that one of $\G_A/\G_A([p])$. 
Moreover, $\G_A/\G_A(p)$ is a {\sl complexification} of $\U_A/\U_A(p)$, 
in the sense that there
exists an anti-holomorphic diffeomorphism in $\G_A/\G_A(p)$ 
whose set of fixed points coincides
with $\U_A/\U_A(p)$: 

The mapping 
$a\G_A(p)\mapsto (a^*)^{-1}\G_A(p), \ \G_A/\G_A(p)\to \G_A/\G_A(p)$ 
is an anti-holomorphic
diffeomorphism (which corresponds to the mapping 
$apa^{-1}\mapsto (a^*)^{-1}pa^*$ in terms of
orbits). Then $a\G_A(p)=(a^*)^{-1}\G_A(p)$ if and only if 
$(a^*a)\G_A(p)=\G_A(p)$, that is, 
$(a^*a)p=p(a^*a)$. Using the functional calculus for $C^*$-algebras, 
we can pick 
$b:=+\sqrt{a^*a}$ in $A$ and obtain $bp=pb$. Since $a^*a=b^2=b^*b$ we have 
$(ab^{-1})^*=(b^{-1})^*a^*=(b^*)^{-1}a^*=ba^{-1}=(ab^{-1})^{-1}$ and 
therefore 
$u:=ab^{-1}\in \U_A$. 
Finally, $a\G_A(p)=ub\G_A(p)=u\G_A(p)\equiv u\U_A(p)\in \U_A/\U_A(p)$.

According to Proposition~\ref{isotropy} (i), 
idempotents like $apa^{-1}\equiv a\G_A(p)$, for
$a\in \G_A$, can be alternatively represented as pairs 
$(a[p]a^{-1},(a^*)^{-1}[p]a^*)$ so that
the 
\lq\lq orbit" $\G_A/\G_A(p)$ becomes a subset of 
the Cartesian product 
$\G_A([p])\times \G_A([p])$. 
In this perspective, the preceding projection and diffeomorphism are
given, respectively, by
$$
(a[p]a^{-1},(a^*)^{-1}[p]a^*)\mapsto a[p]a^{-1}\equiv 
(a[p]a^{-1},a[p]a^{-1}),\ 
\G_A/\G_A(p)\to \G_A/\G_A([p])
$$
(so $(u[p]u^{-1},u[p]u^{-1})\mapsto upu^{-1}\equiv (u[p]u^{-1},u[p]u^{-1})$, 
when $u\in \U_A$)
and
$$
(a[p]a^{-1},(a^*)^{-1}[p]a^*)\mapsto ((a^*)^{-1}[p]a^*, a[p]a^{-1}),
$$
for every $a\in \G_A$.
\qed
\end{remark}

\begin{remark}\label{quasitauto}
\normalfont
Proposition~\ref{complexes} relates to the setting of \cite{BR07}. 
Namely, assume that $B$ is a $C^*$-subalgebra of $A$, 
with $\1\in B\subseteq A$, 
for which there exist a conditional expectation $E\colon A\to B$ and a state 
$\varphi\colon A\to{\mathbb C}$ such that $\varphi\circ E=\varphi$. 
For $X\in\{A,B\}$, 
we denote by $\varphi_X$ the state $\varphi$ restricted to $X$. 
Let $\Hc_X$ be the Hilbert space, and let $\pi_X\colon X\to\Bc(\Hc_X)$ be the corresponding 
cyclic representation obtained by the Gelfand-Naimark-Segal (GNS, for short) 
construction applied to the state $\varphi_X\colon X\to{\mathbb C}$. 
Thus, $\Hc_X$ is 
the completion of $X/N_X$ with respect to the norm 
$\Vert y+N_X\Vert_\varphi:=\varphi(y^*y)$, where $N_X:=\{y\in X\mid\varphi(y^*y)=0\}$. 
The representation 
$\pi_X$ is then defined as the extension to $\Hc_X$ of 
the left multiplication of $X$ on $X/N_X$. 
Let $P$ denote the orthogonal projection $P\colon\Hc_A\to\Hc_B$. 

An equivalence relation can be defined in $\G_A\times\Hc_B$ by setting that 
$(g_1,h_1)\sim(g_2,h_2)$ (with $g_1,g_2\in \G_A$, $h_1,h_2\in\Hc_B$) 
if and only if there exists
$w\in \G_B$ such that $g_2=g_1w$ and $h_2=\pi_B(w^{-1})h_1$. 
The corresponding quotient space will be denoted by $\G_A\times_{\G_B}\Hc_B$, 
and the equivalnce class 
in $\G_A\times_{\G_B}\Hc_B$ of the element $(g,h)\in \G_A\times\Hc_B$ will be denoted
by $[(g,h)]$. 
Define $\U_A\times_{\U_B}\Hc_B$ in an analogous fashion. 
Then the mappings 
$\Pi_G\colon[(g,h)]\mapsto g\G_B,\  \G_A\times_{\G_B}\Hc_B\to \G_A/\G_B$ and 
$\Pi_U\colon[(u,h)]\mapsto u\U_B,\ \U_A\times_{\U_B}\Hc_B\to \U_A/\U_B$ 
are vector bundles, 
$\Pi_U$ being Hermitian, in fact. 
Moreover, $\Pi_U$ admits a reproducing kernel $K$ 
with the associated 
Hilbert space $\Hc^K$, formed by continuous sections of $\Pi_U$, 
such that the restriction of the GNS representation $\pi_A$ to $\U_A$ can be realized on $\Hc^K$, 
see \cite{BR07}.

Let us apply the above theory to the case when, for a given unital $C^*$-algebra $A$, 
we take $B:=\{p\}'$ in $A$, where $p=p^*\in\Pc(A)$. 
Then $E_p\colon a\mapsto pap+\hat ap\hat a,\ A\to B$ is 
a conditional expectation from $A$ onto $B$. 
Let $\Hc$ be a Hilbert space such that 
$A\hookrightarrow\Bc(\Hc)$. 
Pick $x_0\in p\Hc$ such that $\Vert x_0\Vert=1$. 
Then 
$\varphi_0\colon A\to{\mathbb C}$, 
given by $\varphi_0(a):=(ax_0\mid x_0)_{\Hc}$ for all $a\in A$, 
is a state of $A$ such that $\varphi_0\circ E_p=\varphi_0$. 
The GNS representation of $A$ associated with $\varphi_0$ is as follows. 
Set 
$(a_1\mid a_2)_0:=\varphi_0(a_2^*a_1)=(a_2^*a_1x_0\mid x_0)_{\Hc}
=(a_1x_0\mid a_2x_0)_{\Hc}$ 
for every $a_1,a_2\in A$. So $\varphi_0(a^*a)=\Vert a(x_0)\Vert^2$ for all
$a\in A$, whence the null space of $(\cdot\mid\cdot)_0$ is 
$N_0:=\{a\in A:(a\mid a)_0=0\}=\{a\in A:a(x_0)=0\}$. 
The norm $\Vert\cdot\Vert_0$ induced by 
$(\cdot\mid\cdot)_0$ on $A/N_0$ is given by 
$\Vert h\Vert_0\equiv\Vert a+N_0\Vert_0:=\varphi_0(a^*a)^{1/2}
=\Vert a(x_0)\Vert_{\Hc}=\Vert h\Vert_{\Hc}$
for every $h\in A(x_0)\subset\Hc$, where $a(x_0)=h\leftrightarrow a+N_0$. 
Hence $\Hc_A$ 
is a closed 
subspace of $\Hc$ such that $a\Hc_A\subset \Hc_A$ for every $a\in A$. 
Note that $\Hc_A$ 
coincides with 
$\Hc$ provided that we can choose $x_0$ in $\Hc$ such that $A(x_0)$ is
 dense in $\Hc$. 
This will be of interest 
in Remark~\ref{motiv} below.

Analogously, we can consider the restriction of 
$(\cdot\mid\cdot)_0$ to $B$ and 
proceed in the same way as above. 
Thus we obtain that the corresponding null space is $B\cap N_0$, 
that the norm in $B/(B\cap N_0)$ is that one of $p\Hc$ (so that one of $\Hc$), 
and that $\Hc_B$ is a closed subspace of $p\Hc$ such 
that $b\Hc_B\subset \Hc_B$ for every $b\in B$. 
Also, 
$\Hc_B=p\Hc$ if $x_0$ can be chosen in $p\Hc$ and such that $B(x_0)$ is 
dense in $p\Hc$. 

The representation $\pi_A\colon a\mapsto\pi(a),\ A\to\Bc(\Hc_A)$ is 
the extension to
$\Hc_A$ of the left multiplication 
$\pi_A(a)\colon a'+N_0\mapsto(aa')+N_0,\ A/N_0\to A/N_0$. 
Thus it satisfies  
$\pi_A(a'+N_0)=(aa')+N_0\equiv a(a'x_0)=a(h)$, 
if $(a'+N_0)\leftrightarrow a(x_0)=h$. 
In other words, $\pi_A$ is the inclusion operator (by restriction) from $A$ 
into $\Bc(\Hc_A)$. Also, $\pi_B$ 
is in turn the inclusion operator from $B$ into $\Bc(\Hc_B)$. 

Since $E_p(N_0)\subseteq N_0$, the expectation $E_p$ induces 
a well-defined projection 
$P\colon A/N_0\to B/(N_0\cap B)$. On the other hand, 
$E_p(a^*a)-E_p(a)^*E_p(a)=pa^*\hat pap+\hat pa^*pap\ge0$ since $p,\hat p\ge0$. 
Hence $P$ extends once 
again as a bounded projection $P\colon\Hc_A\to\Hc_B$. Indeed, if $h=a(x_0)$ with $a\in A$, 
we have 
$$ 
P(h)\equiv P(a+N_0)=E(a)+(B\cap N_0)=E(a)(x_0)=(pa)(x_0)=p(h),
$$
that is, $P=p_{\mid\Hc_A}$.
\qed
\end{remark}

In the above setting, note that $\U_B=\U(p)$. 
Let $\Gamma(\U_A/\U(p),\U_A\times_{\U(p)}\Hc_B)$ be 
the section space of the bundle $\Pi_U$. 
The reproducing kernel associated
with $\Pi_U$ is given by 
$K_p(u_1\U(p),u_2\U(p))[(u_2,f_2)]:=[(u_1,pu_1^{-1}u_2f_2)]$, for every 
$u_1,u_2\in \U_A$ and $f_2\in\Hc_B$. 
The kernel $K_p$ generates a Hilbert subspace $\Hc^{K_p}$ 
of sections in $\Gamma(\U_A/\U(p),\U_A\times_{\U(p)}\Hc_B)$. 
Let $\gamma_p\colon\Hc_A\to\Gamma(\U_A/\U(p),\U_A\times_{\U(p)}\Hc_B)$ 
be the mapping defined
by $\gamma_p(h)(u\U(p)):=[(u,pu^{-1}h)]$ for every $h\in\Hc_A$ and $u\in \U_A$. 
Then $\gamma_p$ is injective 
and it intertwines the representation $\pi_A$ of $\U_A$ on $\Hc_A$ and 
the natural action of $\U_A$ 
on $\Hc^{K_p}$; that is, the diagram
\begin{equation}\label{intertwin}
\begin{CD}
\Hc_A @>{u}>> \Hc_A \\
@V{\gamma_p}VV @VV{\gamma_p}V \\
\Hc^{K_p} @>{\mu(u)}>> \Hc^{K_p},
\end{CD}
\end{equation}
is commutative for all $u\in \U_A$, where $\mu(u)F:=uF(u^{-1}\ \cdot\ )$ for every 
$F\in\Gamma(\U_A/\U(p),\U_A\times_{\U(p)}\Hc_B)$.  
In fact $\gamma(uh)(v\U(p)):=[(v,pv^{-1}uh)]=u[(u^{-1}v,pv^{-1}uh)]=:u\{\gamma(h)(u^{-1}v\U(p))\}$ 
for all $u,v\in \U_A$. 
See Theorem~5.4 of \cite{BR07} for details in the general case. 
We next show that $\Hc^{K_p}$ in fact consists of holomorphic sections.  

\begin{proposition}\label{cbundle}
Let $A$ be a unital $C^*$-algebra, $p=p^*\in\Pc(A)$, and $B:=\{p\}'$. 
In the above notation, the homogeneous Hermitian vector bundle 
$\Pi_U\colon \U_A\times_{\U(p)}\Hc_B\to \U_A/\U(p)$ is holomorphic, 
and the image of $\gamma_p$ consists of holomorphic sections. 
Thus  
$\Hc^{K_p}$ is a Hilbert space of holomorphic sections of~$\Pi_U$.   
\end{proposition}

\begin{proof}
Let $u_0\in \U_A$. 
Then $\Omega_G:=\{u_0g\mid g\in \G_A,\ \Vert\1-g^{-1}\Vert<1\}$ 
is open in $\G_A$ and contains~$u_0$, 
and similarly with $\Omega_U:=\Omega_G\cap \U_A$ in $\U_A$. 

It is readily seen that the mapping 
$\psi_0\colon$ $[(u,f)]\mapsto(u\U(p),E_p(u^{-1}u_0^{-1})^{-1}f),\ 
\Pi_U^{-1}(\Omega_U)\to\Omega_U\times\Hc_B$ is a diffeomorphism, with inverse map 
$(u\U(p),h)\mapsto[(uE_p(u^{-1}),h)]$ (this shows the local triviality of $\Pi_U$). 
Thus
every point in the manifold $\U_A\times_{\U(p)}\Hc_B$ has an open neighborhood which
is diffeomorphic to the manifold product $W\times\Hc_B$, 
where $W$ is an open subset of $\U_A/\U(p)$. 
By Proposition~\ref{complexes}, $\U_A/\U(p)$ is 
a complex homogeneous manifold 
and therefore the manifold $\U_A\times_{\U(p)}\Hc_B$ is locally complex, 
i.e., holomorphic. 
Also the bundle map $\Pi_U$ is holomorphic.

On the other hand, for fixed $h\in\Hc_A$, the mapping   
$\sigma_0\colon g\G_A([p])\mapsto E_p(g^{-1}u_0^{-1})^{-1}pg^{-1}h,\ 
\Omega_G\to\Hc_B$ 
is holomorphic on $\Omega_G$, so it defines a holomorphic function 
$\tilde\sigma_0\colon \Omega_G\G_A([p])\to\Hc_B$. 
By Proposition~\ref{complexes} the injection 
$j\colon \U_A/\U(p)\hookrightarrow \G_A/\G_A([p])$ is holomorphic, 
and so the restriction map 
$r:=\tilde\sigma_0\circ j$ is holomorphic around $u_0\U(p)$. Since 
$\gamma(h)=\psi_0^{-1}\circ(I_{\Omega_U}\times r)$ around $u_0\U(p)$, 
it follows that 
$\gamma(h)$ is (locally) holomorphic.

Finally, by applying Theorem~4.2 in \cite{BR07} we obtain that $K_p$ is holomorphic. 
\end{proof}

The starting point for the holomorphic picture given in Proposition~\ref{cbundle} has been 
the fact that $\U_A/\U(p)$ enjoys a holomorphic structure induced by the one of $\G_A/\G([p])$, 
see Proposition~\ref{complexes}. Such a picture can be made even more explicit 
if we have a global diffeomorphism $\U_A/\U_A(p)\simeq \G_A/\G_A([p])$. 
The prototypical example is to be found when $A$ is the algebra of bounded 
operators on a complex Hilbert space.  Let us recall the specific definition and 
some properties of the Grassmannian manifold in this case.

\begin{notation}\label{U1}
\normalfont
We shall use the standard notation 
$\Bc(\Hc)$ for the $C^*$-algebra of 
bounded linear operators on the complex Hilbert space $\Hc$ 
with the involution $T\mapsto T^*$.  
Let $\GL(\Hc)$ be the Banach-Lie group of all invertible elements of $\Bc(\Hc)$, and $\U(\Hc)$ 
its Banach-Lie subgroup of all unitary operators on $\Hc$. 
Also, 
\begin{itemize}
\item[$\bullet$] 
  $\Gr(\Hc):=\{\Sc\mid\Sc\text{ closed linear subspace of }\Hc\}$; 
\item[$\bullet$] 
  $\Tc(\Hc):=\{(\Sc,x)\in\Gr(\Hc)\times\Hc\mid x\in\Sc\}
   \subseteq\Gr(\Hc)\times\Hc$; 
\item[$\bullet$] 
  $\Pi_{\Hc}\colon\,(\Sc,x)\mapsto\Sc$, $\Tc(\Hc)\to\Gr(\Hc)$; 
\item[$\bullet$] 
  for every $\Sc\in\Gr(\Hc)$ we denote by $p_{\Sc}\colon\Hc\to\Sc$ 
  the corresponding orthogonal projection.
\end{itemize}
\qed 
\end{notation}

\begin{remark}\label{U2}
\normalfont
The objects introduced in Notation~\ref{U1} have the following well known properties:
\begin{itemize} 
\item[{\rm(a)}] Both $\Gr(\Hc)$ and $\Tc(\Hc)$ have structures of 
  complex Banach manifolds, and $\Gr(\Hc)$ carries 
 a natural (non-transitive) action of $\U(\Hc)$. 
(See Examples 3.11~and~6.20 in \cite{Up85}, 
or Chapter~2 in \cite{Do66}.) 
\item[{\rm(b)}] For every $\Sc_0\in\Gr(\Hc)$ 
the corresponding connected component 
of $\Gr(\Hc)$ is the $\GL(\Hc)$-orbit and is also the $\U(\Hc)$-orbit of $\Sc_0$, that is, 
 $$
\begin{aligned}
\Gr_{\Sc_0}(\Hc)
&=\{g\Sc_0\mid g\in\GL(\Hc)\}=\{u\Sc_0\mid u\in\U(\Hc)\} \\
&=\{\Sc\in\Gr(\Hc)\mid\dim\Sc=\dim\Sc_0\text{ and }\dim\Sc^\perp=\dim\Sc_0^\perp\}
\simeq\U(\Hc)/(\U(\Sc_0)\times\U(\Sc_0^\perp)).
\end{aligned}
 $$
(See Proposition~23.1 in \cite{Up85} or Lemma~\ref{U3} below, alternatively.) 
\item[{\rm(c)}] The mapping $\Pi_{\Hc}\colon\Tc(\Hc)\to\Gr(\Hc)$ 
is a holomorphic Hermitian vector bundle, and we call it the 
{\it universal (tautological) vector bundle} associated with  the Hilbert space $\Hc$. 
Set $\Tc_{\Sc_0}(\Hc):=\{(\Sc,x)\in\Tc(\Hc)\mid\Sc\in\Gr_{\Sc_0}(\Hc)\}$. The vector bundle 
$\Tc_{\Sc_0}(\Hc)\to\Gr_{\Sc_0}(\Hc)$ obtained by restriction of 
$\Pi_{\Hc}$ to $\Tc_{\Sc_0}(\Hc)$ 
will be called here the universal vector bundle {\it at $\Sc_0$}. It is also Hermitian and
holomorphic.
\end{itemize}
\qed
\end{remark}

Property (b) in Remark~\ref{U2} means that $U_A/U_A(p_{\Sc_0})\simeq G_A/G_A([p_{\Sc_0}])$ for
$A=\Bc(\Hc)$. 
For the sake of clarification we now connect Notation~\ref{set} and 
Notation~\ref{U1} in
more detail.  
For $A=\Bc(\Hc)$ we have $\Gr(A)=\Gr(\Hc)$, 
and with this identification the action $\beta$ of Lemma~\ref{beta} 
corresponds to the natural action (so-called {\it collineation} action) of 
the group of invertible operators on 
$\Hc$ upon the set of all closed linear subspaces of $\Hc$. 
The following lemma gives us the collineation orbits of $\Gr(\Hc)$ in terms of orbits of projections, 
and serves in particular to explain the property stated in Remark~\ref{U2}(b).
 
For short, denote $\Gc=\GL(\Hc)$ and $\Uc=\U(\Hc)$. 

\begin{lemma}\label{U3}
Let $\Sc_0\in\Gr(\Hc)$. Then the following assertions hold.
\begin{itemize}
\item[{\rm(i)}]$\Gc([p_{\Sc_0}])=\{g\in\Gc\mid g\Sc_0=\Sc_0\}$ and 
$\Uc([p_{\Sc_0}])=\Uc(p_{\Sc_0})=\{u\in\Uc\mid u\Sc_0=\Sc_0\}$.
\item[{\rm(ii)}] For every $g\in\Gc$ and $\Sc=g\Sc_0$ we have 
$\Sc^\perp=(g^*)^{-1}(\Sc_0^\perp)$.
\item[{\rm(iii)}] We have 
$$
\begin{aligned}
\Gr_{\Sc_0}(\Hc)
&=\{g\Sc_0\mid g\in\Gc\}\simeq\{[gp_{\Sc_0}g^{-1}]\mid g\in\Gc\} \\
&=\{u\Sc_0\mid u\in\Uc\}\simeq\{up_{\Sc_0}u^{-1}\mid u\in\Uc\}.
\end{aligned}
$$
\item[{\rm(iv)}]  We have
$$  
\Uc/\Uc(p_{\Sc_0})\simeq\Gc/\Gc([p_{\Sc_0}])\simeq\Gr_{\Sc_0}(\Hc),  
$$
where the symbol \lq\lq\ $\simeq$ " means diffeomorphism between the respective 
differentiable structures, and that the differentiable structure of the quotient spaces 
is the one associated with the corresponding quotient topologies. 
\item[{\rm(v)}] 
$\Gc/\Gc(p_{\Sc_0})\simeq\{(a\Sc_0,(a^*)^{-1}\Sc_0)\mid a\in\Gc\}$ and the map 
$(a\Sc_0,(a^*)^{-1}\Sc_0)\mapsto((a^*)^{-1}\Sc_0,a\Sc_0)$ is 
an involutive diffeomorphism on $\Gc/\Gc(p_{\Sc_0})$. 
Its set of fixed points is $\Gr_{\Sc_0}(\Hc)\equiv \{(u\Sc_0,u\Sc_0)\mid u\in\Uc\}$.
\end{itemize}
\end{lemma}

\begin{proof}
(i) As shown in Proposition~\ref{isotropy}, an element $g$ of $\Gc$ belongs to $\Gc([p_{\Sc_0}])$ 
if and only if $p_{\Sc_0}g^{-1}p_{\Sc_0}=g^{-1}p_{\Sc_0}$ and 
$p_{\Sc_0}g\ p_{\Sc_0}=g\ p_{\Sc_0}$. 
From this, it follows easily that 
$g(\Sc_0)\subset\Sc_0$ and $g^{-1}(\Sc_0)\subset\Sc_0$, that is, $g(\Sc_0)=\Sc_0$. 
Conversely, if $g(\Sc_0)\subset\Sc_0$ then $(g\ p_{\Sc_0})(\Hc)\subset p_{\Sc_0}(\Hc)$ 
whence $p_{\Sc_0}g\ p_{\Sc_0}=g\ p_{\Sc_0}$; similarly, $g^{-1}(\Sc_0)\subset\Sc_0$ 
implies that $p_{\Sc_0}g^{-1}p_{\Sc_0}=g^{-1}p_{\Sc_0}$. 
In conclusion, 
$\Gc([p_{\Sc_0}])=\{g\in\Gc\mid g\Sc_0=\Sc_0\}$.

Now, the above equality and Remark~\ref{fixed} imply that 
$\Uc([p_{\Sc_0}])=\Uc(p_{\Sc_0})=\{u\in\Uc\mid u\Sc_0=\Sc_0\}$.

(ii) Let $x\in\Sc_0^\perp$, $y\in\Sc$. 
Then 
$((g^*)^{-1}(x)\mid y)=((g^{-1})^*(x)\mid y)=(x\mid g^{-1}(y))=0$, so 
$(g^*)^{-1}(\Sc_0^\perp)\subset\Sc^\perp$. 
Take now $y\in\Sc^\perp$, $x=g^*(y)$ and $z\in\Sc_0$. 
Then 
$(x\mid z)=(g^*(y)\mid z)=(y\mid g(z))=0$, whence $x\in\Sc_0^\perp$ and therefore 
$y=(g^*)^{-1}(g^*y)=(g^*)^{-1}(x)\in(g^*)^{-1}(\Sc_0^\perp)$. In conclusion, 
$\Sc^\perp=(g^*)^{-1}(\Sc_0^\perp)$.

(iii) By (ii), we have $u(\Sc_0^\perp)=u(\Sc_0)^\perp$ for $u\in\Uc$. 
Thus $\Sc=u(\Sc_0)$ if and only if 
$\dim\Sc=\dim\Sc_0$ and $\dim\Sc^\perp=\dim\Sc_0^\perp$.
Also, if $\Sc=u(\Sc_0)$ and $\Sc^\perp=u(\Sc_0^\perp)$, 
then $up_{\Sc_0}=p_{\Sc}u$, 
that is,
$p_{\Sc}=up_{\Sc_0}u^{-1}$. Hence
$\Gr_{\Sc_0}(\Hc)=\{u\Sc_0\mid u\in\Uc\}=\{\Sc\in\Gr(\Hc)\mid\dim\Sc=\dim\Sc_0\text
{ and }\dim\Sc^\perp=\dim\Sc_0^\perp\}\simeq\{up_{\Sc_0}u^{-1}\mid u\in\Uc\}$.

Suppose now that $\Sc=g\Sc_0$ with $g\in\Gc$. Then $\dim\Sc=\dim\Sc_0$. 
By (ii) again,
$\Sc^\perp=(g^*)^{-1}(\Sc_0^\perp)$ and so $\dim\Sc^\perp=\dim\Sc_0^\perp$. 
Hence $\Sc\in\Gr_{\Sc_0}(\Hc)$. 
Finally, the bijective correspondence between $g\Sc_0$ 
and $g[p_{\Sc_0}]g^{-1}$ is straightforward.

(iv) This is clearly a consequence of parts (iii) and (i) from above, 
and Theorem 2.2 in \cite{DG02}. 

(v) For every $a\in\Gc$, the pairs $(a\Sc_0,(a^*)^{-1}\Sc_0)$ and 
$(a[p]a^{-1},(a^*)^{-1}[p]a^*)$ are in a
one-to-one correspondence, by part (iii) from above. Hence, this part (v) is 
a consequence of 
Remark~\ref{picture}. 
\end{proof}

Parts (iv) and (v) of Lemma~\ref{U3} tell us that the Grassmannian orbit
$\Gr_{\Sc_0}(\Hc)$ is a complex manifold which in turn admits a complexification, 
namely the orbit
$\Gc/\Gc(p_{\Sc_0})$. 

\begin{remark}\label{connected}
\normalfont
As said in Remark~\ref{U2}(b), every $\GL(\Hc)$-orbit (and so every $\U(\Hc)$-orbit) is 
a connected component of $\Gr(\Hc)$. Let us briefly discuss the connected components of 
$\Gr(A)$ when $A$ is an arbitrary unital $C^*$-algebra. 
Every element $g\in \G_A$ has a unique polar
decomposition $g=ua$ with $u\in\U_A$ and $0\le a\in \G_A$, 
hence there exists a continuous path $t\mapsto u\cdot((1-t)\1+ta)$ in $\G_A$ 
that connects $u=u\cdot\1$ to $g=u\cdot a$. 
Thus every connected component of the $\G_A$-orbit of $[p]\in\Gr(A)$ 
contains at least one connected component 
of the $\U_A$-orbit of $[p]\in\Gr(A)$ for any idempotent $p\in\Pc(A)$. 
(Loosely speaking, the $\U_A$-orbit of $[p]$ has more connected components 
than the $\G_A$-orbit of $[p]$.) 
Example~7.13 in \cite{PR87} shows that the $C^*$-algebra 
$A$ of the continuous functions $S^3\to M_2({\mathbb C})$ 
has the property that there indeed exist $\G_A$-orbits of elements 
$[p]\in\Pc(A)$ which are nonconnected. 

If the unitary group $\U_A$ is connected
(so that the invertible group $\G_A$ is connected), 
then all the $\U_A$-orbits and the $\G_A$-orbits in $\Gr(A)$ 
are connected since continuous images of connected sets are always connected. 
On the other hand, as said formerly, the Grassmannian $\Gr(A)$ is the discrete union 
of these $\G_A$-orbits. Thus if the unitary group $\U_A$ is connected, 
then the connected components of $\Gr(A)$ are precisely the $\G_A$-orbits in $\Gr(A)$. 
One important case of connected unitary group $\U_A$ is when $A$ is 
a $W^*$-algebra (since every $u\in \U_A$ can be written as $u=\exp(\ie a)$ 
for some $a=a^*\in A$ by the Borel functional calculus, 
hence the continuous path $t\mapsto\exp(\ie ta)$ 
connects $\1$ to $u$ in $\U_A$). 
For $W^*$-algebras such that $\Gr(A)$ is the discrete union of
$\U_A$-orbits, it is then clear that the $\G_A$-orbits and 
the $\U_A$-orbits coincide.   
This is the case if $A$ is the algebra of bounded operators on 
a complex Hilbert space, as we have seen before.
\qed
\end{remark}

The universal bundle $\Tc_{\Sc_0}(\Hc)\to\Gr_{\Sc_0}(\Hc)$ can be expressed 
as a vector bundle 
obtained from the so-called (principal) {\it Stiefel bundle} associated to 
$p_{\Sc_0}\leftrightarrow\Sc_0$, see \cite{DG02}. A similar result holds, 
by replacing the Stiefel bundle with certain, suitable, of its sub-bundles. 
To see this, let us now introduce several mappings.

Put $p:=p_{\Sc_0}$. 
We consider $\Gc\times_{\Gc([p])}\Sc_0$ and 
$\Uc\times_{\Uc(p)}\Sc_0$ as in Remark~\ref{quasitauto}. 
Note that $g_1\Sc_0=g_2\Sc_0$ and
$g_1(h_1)=g_2(h_2)$ ($g_1,g_2\in\Gc$, $h_1,h_2\in\Sc_0$) if and only if $(g_1,h_1)\sim(g_2,h_2)$, 
via $w=g_1^{-1}g_2\in\Gc([p])$, in $\Gc\times\Sc_0$. 
Hence, the mapping 
$\upsilon_{\Gc}\colon\Gc\times\Sc_0\to\Tc_{\Sc_0}(\Hc)$ defined by 
$\upsilon_{\Gc}((g,h))=(g\Sc_0,g(h))$ for every $(g,h)\in\Gc\times\Sc_0$, induces the 
usual (canonical) quotient map 
$\tilde\upsilon_{\Gc}\colon\Gc\times_{\Gc([p])}\Sc_0\to\Tc_{\Sc_0}(\Hc)$. 
We denote by $\upsilon_{\Uc}$ the restriction map of $\upsilon_{\Gc}$ on $\Gc\times\Sc_0$. 
As above, the quotient mapping 
$\tilde\upsilon_{\Uc}\colon\Uc\times_{\Uc(p)}\Sc_0\to\Tc_{\Sc_0}(\Hc)$ 
is well defined.

Since $\Uc(p)=\Uc\cap\Gc([p])$, the inclusion mapping 
$j\colon\Uc\times_{\Uc(p)}\Sc_0\to\Gc\times_{\Gc([p])}\Sc_0$ is well defined. 
Note that 
$j=(\tilde\upsilon_{\Gc})^{-1}\circ\tilde\upsilon_{\Uc}$.

Finally, let $P_\Gc\colon\Gc\times_{\Gc([p])}\Sc_0\to\Gc/\Gc([p])$ and 
$P_\Uc\colon\Uc\times_{\Uc(p)}\Sc_0\to\Uc/\Uc(p)$ denote the vector bundles built in 
the standard way from the Stiefel sub-bundles 
$g\mapsto g\Gc([p])\simeq g(\Sc_0),\ \Gc\to\Gc/\Gc([p])\simeq\Gr_{\Sc_0}(\Hc)$ and 
$u\mapsto u\Uc(p)\simeq u(\Sc_0),\ \Uc\to\Uc/\Uc(p)\simeq\Gr_{\Sc_0}(\Hc)$ respectively.

\begin{proposition}\label{tautoquo}
The following diagram is commutative in both sides, and the horizontal arrows are diffeomorphisms
between the  corresponding differentiable structures
$$
\begin{CD}
\Tc_{\Sc_0}(\Hc) @>{(\tilde\upsilon_{\Uc})^{-1}}>> 
\Uc\times_{\Uc(p)}\Sc_0 @>{j}>>
\Gc\times_{\Gc([p])}\Sc_0 \\ 
@V{\Pi_{\Hc}}VV @VV{P_{\Uc}}V @VV{P_{\Gc}}V \\
\Gr_{\Sc_0}(\Hc) @>{\simeq}>> \Uc/\Uc(p) @>{\simeq}>> \Gc/\Gc([p])
\end{CD}
$$
\end{proposition}

\begin{proof}
By construction, the mapping $\tilde\upsilon_{\Uc}$ is clearly one-to-one. 
Now we show that it is onto. Let $(\Sc,h)\in\Tc_{\Sc_0}(\Hc)$. This means that $h\in\Sc$ and 
that $\Sc=u\Sc_0$ for some $u\in\Uc$. Then $f:=u^{-1}(h)\in\Sc_0$ and $h=u(f)$, whence 
$\tilde\upsilon_{\Uc}([(u,f)])=(\Sc,h)$, where $[(u,f)]$ is the equivalence class of $(u,f)$ 
in $\Uc\times_{\Uc(p)}\Sc_0$. Hence $\tilde\upsilon_{\Uc}$ is a bijective map. 

Analogously, we have that $\tilde\upsilon_{\Gc}$ is bijective from 
$\Gc\times_{\Gc([p])}\Sc_0$ onto $\Tc_{\Sc_0}(\Hc)$ as well. 
As a consequence, 
$j=(\tilde\upsilon_{\Gc})^{-1}\circ\tilde\upsilon_{\Uc}$ is also bijective.
It is straightforward to check that all the maps involved in the diagram above are smooth.
\end{proof}

\begin{example}\label{universal}
\normalfont
By Proposition~\ref{tautoquo} one can show that 
the universal, tautological bundle 
$\Pi_{\Hc}\colon\Tc_{\Sc_0}(\Hc)\to\Gr_{\Sc_0}(\Hc)$ enters, as a canonical example, 
the framework outlined in Theorem 5.4 and Theorem 5.8 of \cite{BR07}. 

To see this in terms of the bundle $\Pi_{\Hc}$ itself, 
first note that the commutant algebra 
$\{p_{\Sc_0}\}'$ of $p_{\Sc_0}$ coincides with the Banach subalgebra $B$ of $A$ formed by 
the operators $T$ such that $T(\Sc_0)\subset\Sc_0$, $T(\Sc_0^\perp)\subset\Sc_0^\perp$. 
(It is straightforward to check directly on $B$ that it is stable
under the adjoint operation, so that $B$ is a $C^*$-subalgebra of
$A$, as it had to be.) Put $p=p_{\Sc_0}$. 
From Lemma~\ref{U3}, $u\in\Uc([p])$ if and only if
$u\Sc_0=\Sc_0$. Hence $u\in\Uc(p)=\Uc([p])\cap\Uc([1-p])$ if and only if 
$u\Sc_0=\Sc_0$ and  $u\Sc_0^\perp=\Sc_0^\perp$, that is, 
$\Uc(p)=\Uc_A\cap B=\Uc_B$.

Similarly to what has been done in Remark~\ref{quasitauto}, let $E_p\colon A\to B$ denote the
canonical expectation associated to the  tautological bundle at $\Sc_0$; that is,
$E_p(T):=pTp+\hat{p}T\hat{p}$ for every
$T\in A$. Also, for a fixed $x_0\in\Sc_0$ such that
$\Vert x_0\Vert=1$, let $\varphi\colon A\to{\mathbb C}$ 
be the state of $A$ given by 
$\varphi_0(T):=(Tx_0\mid x_0)_{\Hc}$. Then $\varphi_0\circ E_p=\varphi_0$. 
Since the mappings $T\mapsto T(x_0),\ \Bc(\Hc)\to\Hc$ and 
$T\mapsto T(x_0),\ B\to\Sc_0$ are surjective, we obtain that $\Hc_A=\Hc$ and $\Hc_B=\Sc_0$ 
in the GNS construction associated with $A=\Bc(\Hc)$, $B$ and $\varphi_0$. 
Moreover, in this case, $\pi_A$ coincides with the identity operator and the extension 
$P\colon\Hc_A\to\Hc_B$ of $E_p$ is $P=p$.
Denote by $p_1,p_2\colon\Gr(\Hc)\times\Gr(\Hc)\to\Gr(\Hc)$ 
the natural projections and define 
$$
Q_{\Hc}\colon\Gr(\Hc)\times\Gr(\Hc)\to
\Hom(p_2^*(\Pi_{\Hc}),p_1^*(\Pi_{\Hc}))
$$
by 
$$
Q_{\Hc}(\Sc_1,\Sc_2)=(p_{\Sc_1})|_{\Sc_2}\colon\Sc_2\to\Sc_1
$$
whenever $\Sc_1,\Sc_2\in\Gr(\Hc)$. 
This mapping $\Qc_{\Hc}$ is called the {\it universal reproducing kernel} 
associated with the Hilbert space $\Hc$. In fact, 
for $\Sc_1,\dots,\Sc_n\in\Gr(\Hc)$ and $x_j\in\Sc_j$ ($j=1,\dots,n$), 
$$
\sum_{j,l=1}^n(Q_{\Hc}(\Sc_l,\Sc_j)x_j\mid x_l)_{\Hc}
=\sum_{j,l=1}^n(p_{\Sc_l}x_j\mid x_l)_{\Hc}
=\sum_{j,l=1}^n(x_j\mid x_l)_{\Hc}
=(\sum_{j=1}^nx_j\mid\sum_{l=1}^nx_l)_{\Hc}\ge0, 
$$
so $\Qc_{\Hc}$ is certainly a reproducing kernel in the sense of \cite{BR07}. 
\qed
\end{example}

Using Example~\ref{universal} we get the following special case of 
Theorem~5.8 in \cite{BR07}.  

\begin{corollary}\label{realiztauto} 
For a complex Hilbert space $\Hc$, the action of $\Uc$ on 
$\Hc$ can be realized as the natural action of $\Uc$ on a Hilbert space of holomorphic sections 
from $\Gr_{\Sc_0}(\Hc)$ into $\Hc$, such a realization being implemented by 
$\gamma(uh)=u\ \gamma(h)u^{-1}$, 
for every $h\in\Hc$, $u\in\Uc$. 
\end{corollary}

\begin{proof}
If $\Sc\in\Gr_{\Sc_0}(\Hc)$, there exists 
$u\in\Uc$ such that $u\Sc_0=\Sc$ and then $p_{\Sc}=up_{\Sc_0}u^{-1}$. 
Thus for all $u_1,u_2\in\Uc$ and $x_1,x_2\in\Sc_0$ 
we have 
$Q_{\Hc}(u_1\Sc_0,u_2\Sc_0)(u_2x_2)
=p_{u_1\Sc_0}(u_2x_2)=u_1p_{\Sc_0}(u_1^{-1}u_2x_2)$. 
This formula shows that for every connected component $\Gr_{\Sc_0}(\Hc)$ 
the restriction of $Q_{\Hc}$ to $\Gr_{\Sc_0}(\Hc)\times\Gr_{\Sc_0}(\Hc)$ 
is indeed a special case of the reproducing kernels considered in Remark~\ref{quasitauto}.
For every $h\in\Hc$, the mapping 
$\gamma_{p_{\Sc_0}}(h)\colon\Gr_{\Sc_0}(\Hc)\to\Tc_{\Sc_0}(\Hc)$ 
which corresponds to $Q_{\Hc}$ can be identified to the  holomorphic map
$u\Sc_0\mapsto upu^{-1}h,\ \Gr_{\Sc_0}(\Hc)\to \Hc$. 
Then the conclusion follows by using the diffeomorphism 
$\Uc/\Uc(p)\simeq\Gc/\Gc([p])\simeq\Gr_{\Sc_0}(\Hc)$ of Lemma~\ref{U3}, together 
with Proposition~\ref{tautoquo}. 
\end{proof}

\begin{remark}\label{motiv}
\normalfont
Assume again the situation where $A$ and $B$ are arbitrary $C^*$-algebras, 
$B$ is a $C^*$-subalgebra of $A$, with unit $\1\in B\subseteq A$, 
$E\colon A\to B$ is a conditional expectation, and $\varphi\colon A\to{\mathbb C}$ 
is a state such that 
$\varphi\circ E=\varphi$. With the same notations as in Remark~\ref{quasitauto}, 
take $x_0:=\1+N_B\in B/N_B\subset A/N_A$. It is well known that $x_0$ is 
a cyclic vector of $\pi_X$, 
for $X\in\{A;B\}$: let $h\in\Hc_X$ such that 
$0=(\pi(c)x_0\mid h)_{\Hc_X}\equiv(c+N_X\mid h)_{\Hc_X}$ for all $c\in X$; 
since $X/N_X$ is dense in $\Hc_X$ we get $0=(h\mid h)_{\Hc_X}=\Vert h\Vert^2$, 
that is, $h=0$. Thus $\pi_X(X)x_0$ is dense in $\Hc_X$. 

Inspired by \cite{AS94}, we now consider the $C^*$-subalgebra $\Ag$ of $\Bc(\Hc_A)$ generated by 
$\pi_A(A)$ and $p$, where $p$ is the orthogonal projection from $\Hc_A$ onto 
$\Hc_B$. Set $\Bg:=\Ag\cap\{p\}'$. Clearly, the GNS procedure is applicable to 
$\Bg\subset\Ag\subset\Bc(\Hc_A)$, for 
the expectation $E_p\colon\Ag\to\Bg$ and state $\varphi_0$ defined by $x_0$, 
as we have done in Remark~\ref{quasitauto}. Then 
$\pi_A(A)(x_0)\subset\Ag(x_0)\subset\Hc_A$ and 
$\pi_A(B)(x_0)\subset\Bg(x_0)\subset\Hc_B$, whence, by the choice of $x_0$, we obtain that 
$\overline{\Ag(x_0)}=\Hc_A$ and $\overline{\Bg(x_0)}=\Hc_B$. Thus we have that 
$\Hc_{\Ag}=\Hc_A$ and $\Hc_{\Bg}=\Hc_B$.

According to former discussions there are two (composed) commutative diagrams, namely
\begin{equation}\label{basic}
\begin{CD} 
\G_A\times_{\G_B}\Hc_B @>{\pi_A\tilde\times I}>> 
\G_{\Ag}\times_{\G_{\Ag}(p)}\Hc_B @>{ }>> \G_{\Ag}\times_{\G_{\Ag}([p])}\Hc_B 
@>{j\tilde\times I}>> \Gc\times_{\Gc([p])}\Hc_B \\
@V{\Pi_G}VV @VV{ \Pi_{G_{\Ag}}}V @VV{  }V @VV{\Pi_{\Hc_B}}V \\
\G_A/\G_B @>{\tilde\pi_A}>> \G_{\Ag}/\G_{\Ag}(p) @>{ }>> \G_{\Ag}/\G_{\Ag}([p]) @>{\tilde j}>>
\Gc/\Gc([p])
\end{CD}
\end{equation}
and
\begin{equation}\label{basicu}
\begin{CD} 
\U_A\times_{\U_B}\Hc_B @>{\pi_A\tilde\times I}>> 
\U_{\Ag}\times_{\U_{\Ag}(p)}\Hc_B @>{j\tilde\times I}>> \Uc\times_{\Uc(p)}\Hc_B @>{\simeq}>>
\Tc_{\Hc_B}(\Hc_A) \\ @V{\Pi_U}VV @VV{ \Pi_{U_{\Ag}}}V @VV{\Pi_{\Uc}}V @VV{\Pi_{\Hc_B}}V \\
\U_A/\U_B @>{\tilde\pi_A}>> \U_{\Ag}/\U_{\Ag}(p) @>{\tilde j}>> \Uc/\Uc(p) @>{\simeq}>>
\Gr_{\Hc_B}(\Hc_A)
\end{CD}
\end{equation}
(where the meaning of the arrows is clear). 
We suggest to call 
$\Pi_G\colon \G_A\times_{\G_B}\Hc_B\to \G_A/\G_B$ and 
$\Pi_U\colon \U_A\times_{\U_B}\Hc_B\to \U_A/\U_B$ 
the {\it GNS vector bundle} and the {\it unitary GNS vector bundle}, respectively, 
for data 
$E\colon A\to B$ and $\varphi\colon A\to{\mathbb C}$. 
Following the terminology used in 
\cite{AS94}, \cite{ALRS97} for the maps 
$\G_A/\G_B\to \G_{\Ag}/\G_{\Ag}(p)$, $\U_A/\U_B\to \U_{\Ag}/\U_{\Ag}(p)$, 
we could refer to the left sub-diagrams of
\eqref{basic} and \eqref{basicu} as the {\it basic} vector bundle
representations of $\Pi_G$ and $\Pi_U$, respectively. Since 
$\Hc_{\Ag}=\Hc_A$ and $\Hc_{\Bg}=\Hc_B$, the process to construct such \lq\lq basic" 
objects, of Grassmannian type, is stationary. Also, 
since there is another way to associate Grassmannians to
the GNS and unitary GNS bundles, which is that one of considering the tautological bundle 
of $\Hc_A$ (see the right diagrams in \eqref{basic}, \eqref{basicu}), 
we  might call $\G_{\Ag}\times_{\G_{\Ag}(p)}\Hc_B\to \G_{\Ag}/\G_{\Ag}([p])$ 
the {\it minimal}
Grassmannian vector bundle,  and call $\Tc_{\Hc_B}(\Hc_A)\to\Gr_{\Hc_B}(\Hc_A)$ 
the {\it universal} Grassmannian vector bundle, associated with data 
$E\colon A\to B$ and $\varphi\colon A\to{\mathbb C}$. 
In the unitary case, we should add the adjective \lq\lq unitary" to both bundles.
 
Note that the vector bundles 
$\Gc\times_{\Gc([p])}\Hc_B\to\Gc/\Gc([p])$ and $\Tc_{\Hc_B}(\Hc_A)\to\Gr_{\Hc_B}(\Hc_A)$ 
are isomorphic. 
In this sense, both diagrams \eqref{basic} and \eqref{basicu} \lq\lq converge" towards 
the tautological bundle for
$\Hc_A$. Let us remark that \eqref{basic} is holomorphic, 
and everything in \eqref{basicu} is holomorphic with the only possible
exception of the bundle $\Pi_U$. On the other hand, we have that   
$\G_{\Ag}/\G_{\Ag}(p)$ and $\Gc/\Gc(p)$ are complexifications of 
$\U_{\Ag}/\U_{\Ag}(p)$ and $\Uc/\Uc(p)$ respectively, on account of Remark~\ref{picture} 
and Lemma~\ref{U3}. 
We shall see in Corollary~\ref{supplcor1} 
that $\G_A/\G_B$ is also 
a complexification of $\U_A/\U_B$ in general. 
Note in passing that the fact 
that $\G_A/\G_B$ is such a complexification implies interesting properties of metric nature 
in the differential geometry of $\U_A/\U_B$, see \cite{ALRS97}.     

The above considerations strongly suggest to investigate the relationships between 
\eqref{basic} and \eqref{basicu} in terms of holomorphy and geometric realizations. 
In this respect, note that the commutativity of \ref{basicu} corresponds, 
at the level of reproducing kernels, with the equality
$$   
(\pi_A\tilde\times I)\circ K(u_1\U_B,u_2\U_B)=
Q_{\Hc_B}(\pi_A(u_1)\Uc(p),\pi_A(u_2)\Uc(p))\circ(\pi_A\tilde\times I)
$$
for all $u_1, u_2\in \U_A$ 
(where the holomorphy supplied by $Q_{\Hc_B}$ appears explicitly). 
From this, a first candidate to reproducing kernel on $\G_A/\G_B$, 
in order to obtain a geometric realization of
$\pi_A$ on $G_A$, would be defined by 
$$
K(g_1G_B,g_2G_B)[(g_2,f)]:=[(g_1,p(\pi_A(g^{-1}_1)\pi_A(g_2)f))]
$$
for every $g_1, g_2\in G_A$ and $f\in\Hc_B$. 
Nevertheless, since the elements $g_1,g_2$ 
are not necessarily unitary, it is readily seen that the kernel $K$ 
so defined need not be definite-positive in general. 
There is also the problem of the existence of a suitable structure 
of Hermitian type in $\Pi_G$.

In the present paper, we propose a theory on bundles 
$\G_A\times_{\G_B}\Hc_B\to \G_A/\G_B$ and kernels
$K$ {\sl ad hoc}, based on the existence of 
suitable involutive diffeomorphisms in 
$\G_A/\G_B$, which allows us to
incorporate those bundles to a framework that contains as a special case 
the one established in \cite{BR07}.      
\qed
\end{remark}

\section{Like-Hermitian structures}\label{sect2}

We are going to introduce a variation of the notion of 
Hermitian vector bundle, 
which will turn out to provide the appropriate setting 
for the geometric representation theory 
of involutive Banach-Lie groups as developed in Section~\ref{sect4}. 

\begin{definition}\label{like}
\normalfont 
Assume that $Z$ is a real Banach manifold equipped with a diffeomorphism 
$z\mapsto z^{-*}$, $Z\to Z$, 
which is involutive in the sense that $(z^{-*})^{-*}=z$ for all $z\in Z$. 
Denote by $p_1,p_2\colon Z\times Z\to Z$ the natural projection maps.
Let $\Pi\colon D\to Z$ be a smooth vector bundle whose fibers 
are complex Banach spaces 
(see for instance \cite{AMR88} or \cite{La01} for details on 
infinite-dimensional vector bundles). 

We define a {\it like-Hermitian structure}
on the bundle $\Pi$ (with typical fiber the Banach space $\Ec$) as a family 
$\{(\cdot\mid\cdot)_{z,z^{-*}}\}_{z\in Z}$ 
with the following properties: 
\begin{itemize} 
\item[{\rm(a)}] 
For every $z\in Z$,  
$(\cdot\mid\cdot)_{z,z^{-*}}\colon D_z\times D_{z^{-*}}\to{\mathbb C}$
is  a sesquilinear strong duality pairing.
\item[{\rm(b)}] 
For all $z\in Z$, $\xi\in D_z$, and $\eta\in D_{z^{-*}}$ 
we have $\overline{(\xi\mid\eta)}_{z,z^{-*}}=(\eta\mid\xi)_{z^{-*},z}$. 
\item[{\rm(c)}] 
If $V$ is an arbitrary open subset of $Z$, 
and $\Psi_0\colon V\times\Ec\to\Pi^{-1}(V)$ 
and $\Psi_1\colon V^{-*}\times\Ec\to\Pi^{-1}(V^{-*})$ 
are trivializations of the vector bundle $\Pi$ 
over $V$ and $V^{-*}$ ($:=\{z^{-*}\mid z\in V\}$),  
respectively, 
then the function 
$(z,x,y)\mapsto(\Psi_0(z,x)\mid\Psi_1(z^{-*},y))_{z,z^{-*}}$, 
$V\times \Ec\times\Ec\to{\mathbb C}$, 
is smooth. 
\end{itemize}
We call {\it like-Hermitian vector bundle}
any vector bundle equipped 
with a  like-Hermitian structure.
\qed
\end{definition}

\begin{remark}\label{sesqui}
\normalfont
Here we explain the meaning of condition~(a) in Definition~\ref{like}. 
To this end let $\Xc$ and $\Yc$ be two complex Banach spaces. 
A functional 
$(\cdot\mid\cdot)\colon\Xc\times\Yc\to{\mathbb C}$
is said to be a {\it sesquilinear strong duality pairing}
if it is continuous, is linear in the first variable 
and antilinear in the second variable, and both the mappings 
$$x\mapsto(x\mid\cdot),\quad\Xc\to(\overline{\Yc})^*,
\quad\text{ and }\quad 
y\mapsto(\cdot\mid y),\quad \overline{\Yc}\to\Xc^*,$$
are
(not necessarily isometric) isomorphisms of complex Banach spaces.  

Here we denote, for any complex Banach space $\Zc$, 
by $\Zc^*$ its dual Banach space 
(i.e., the space of all continuous linear functionals $\Zc\to{\mathbb C}$) 
and by $\overline{\Zc}$ the complex-conjugate Banach space. 
That is, the {\it real} Banach spaces underlying 
$\Zc$ and $\overline{\Zc}$ coincide, and 
for any $z$ in the corresponding real Banach space 
and $\lambda\in{\mathbb C}$ 
we have 
$\lambda\cdot z\;\; \text{(in $\overline{\Zc}$)}
=\overline{\lambda}\cdot z\;\; \text{(in ${\Zc}$)}$.
\qed
\end{remark}

\begin{remark}\label{sesqui1}
\normalfont 
For later use we now record the following fact:
Assume that $\Xc$ and $\Yc$ are two Banach spaces over ${\mathbb C}$, and let  
$(\cdot\mid\cdot)\colon\Xc\times\Yc\to{\mathbb C}$
be a sesquilinear strong duality pairing. 
Now let $\Hc$ be a Hilbert space over ${\mathbb C}$ and 
let $T\colon\Hc\to\Xc$ be 
a continuous linear operator. 
Then there exists a unique linear operator 
$S\colon\Yc\to\Hc$ such that 
\begin{equation}\label{adj}
(\forall h\in\Hc,y\in\Yc)\quad 
(Th\mid y)=(h\mid Sy)_{\Hc}.
\end{equation} 
Conversely,  for every 
bounded linear operator $S\colon\Yc\to\Hc$ there exists 
a unique bounded linear operator $T\colon\Hc\to\Xc$ satisfying~\eqref{adj}, 
and we denote $S^{-*}:=T$ and $T^{-*}:=S$. 
\qed
\end{remark}

\begin{remark}\label{hermitian}
\normalfont
In Definition~\ref{like} if $z^{-*}=z$ 
and $(\xi\mid\xi)_{z,z}\ge0$ for all $z\in Z$ and $\xi\in D_z$, 
then we shall speak simply about {\it Hermitian} structures and bundles, 
since this is just the usual notion of Hermitian structure 
on a vector bundle. 
See for instance Definition~1.1 in Chapter~III of \cite{We80} 
for the classical case of finite-dimensional Hermitian vector bundles. 
\qed
\end{remark}

\begin{example}\label{rel}
\normalfont
Let $\Pi\colon D\to Z$ be a smooth vector bundle whose fibers 
are complex Banach spaces. 
Assume that there exist a complex Hilbert space~$\Hc$ and 
a smooth map $\Theta\colon D\to\Hc$ with the property that 
$\Theta|_{D_z}\colon D_z\to\Hc$ is a bounded linear operator
for all $z\in Z$. 
Then $\Theta$ determines a family of continuous sesquilinear functionals 
$$
(\cdot\mid\cdot)_{z,z^{-*}}\colon D_z\times D_{z^{-*}}\to {\mathbb C}, 
\quad 
(\eta_1\mid\eta_2)_{z,z^{-*}}=(\Theta(\eta_1)\mid\Theta(\eta_2))_\Hc.
$$ 
If in addition $\Theta|_{D_z}\colon D_z\to\Hc$ 
is injective and has closed range, and the scalar product 
of $\Hc$ determines a sesquilinear strong duality pairing 
between the subspaces $\Theta(D_z)$ and $\Theta(D_{z^{-*}})$ 
whenever $z\in Z$, 
then it is easy to see that we get 
a like-Hermitian structure on the vector bundle $\Pi$. 
\qed
\end{example}

\begin{definition}\label{invol}
\normalfont
An {\it involutive} Banach-Lie group is a (real or complex) 
Banach-Lie group $G$ equipped with a diffeomorphism $u\mapsto u^*$ satisfying 
$(uv)^*=v^*u^*$ and $(u^*)^*=u$ for all $u,v\in G$. 
In this case we denote 
$$(\forall u\in G)\quad u^{-*}:=(u^{-1})^*$$ 
and 
$$G^{+}:=\{u^*u\mid u\in G\}$$ 
and the elements of $G^{+}$ are called the {\it positive} elements of $G$. 

If in addition $H$ is a Banach-Lie subgroup of $G$, 
then we say that $H$ is an {\it involutive} Banach-Lie subgroup if $u^*\in H$ 
whenever $u\in H$. 
\qed
\end{definition}

\begin{remark}\label{invol1}
\normalfont
If $G$ is an involutive Banach-Lie group then for every $u\in G$ 
we have $(u^{-1})^*=(u^*)^{-1}$ and moreover $\1^*=\1$. 
To see this, just note that the mapping $u\mapsto(u^*)^{-1}$ 
is an automorphism of $G$, hence it commutes with the inversion mapping 
and leaves $\1$ fixed. 
\qed
\end{remark}

\begin{example}\label{invol2}
\normalfont
Every Banach-Lie group $G$ has a trivial structure of 
involutive Banach-Lie group defined by $u^*:=u^{-1}$ for all $u\in G$. 
In this case the set of positive elements is $G^{+}=\{\1\}$. 
\qed
\end{example}

\begin{example}\label{invol3}
\normalfont
Let $A$ be a unital $C^*$-algebra with the group of invertible elements 
denoted by $\G_A$. 
Then $\G_A$ has a natural structure of involutive complex Banach-Lie group 
defined by the involution of $A$. 
If $B$ is any $C^*$-subalgebra of $A$ such that there exists 
a conditional expectation $E\colon A\to B$, 
then $\G_B$ is an involutive complex Banach-Lie subgroup of $\G_A$. 
\qed
\end{example}

\begin{definition}\label{homog}
\normalfont
Assume that we have the following data: 
\begin{itemize}
\item[$\bullet$] $G_A$ is an involutive real 
(respectively, complex) Banach-Lie group and $G_B$ 
is an involutive real (respectively, complex) Banach-Lie subgroup of $G_A$. 
\item[$\bullet$] For $X=A$ or $X=B$, assume $\Hc_X$ is a complex Hilbert space 
with $\Hc_B$ closed subspace in $\Hc_A$, and 
$\pi_X\colon G_X\to \Bc(\Hc_X)$ is a uniformly continuous 
(respectively, holomorphic) 
$*$-representation such that 
$\pi_B(u)=\pi_A(u)|_{\Hc_B}$ for all $u\in G_B$. 
By $*$-representation we mean that $\pi_A(u^*)=\pi_A(u)^*$ 
for all $u\in G_A$. 
\item[$\bullet$] We denote by $P\colon\Hc_A\to\Hc_B$ 
the orthogonal projection. 
\end{itemize}
We define an equivalence relation on 
$G_A\times\Hc_B$ by 
$$
(u,f)\sim(u',f')\quad\text{ whenever there exists } 
w\in G_B\quad\text{ such that }  u'=uw
\quad\text{ and } f'=\pi_B(w^{-1})f.
$$
For every pair $(u,f)\in G_A\times\Hc_B$ we define its equivalence class by 
$[(u,f)]$ and let $D=G_A\times_{G_B}\Hc_B$ 
denote the corresponding set of equivalence classes. 
Then there exists a natural onto map 
$$
\Pi\colon\,[(u,f)]\mapsto s:=u\ G_B,\quad D\to G_A/G_B.
$$
For $s\in G_A/G_B$, let $D_s:=\Pi^{-1}(s)$ denote the {\it fiber} on $s$.
Note that $(u,f)\sim(u',f')$ implies that $\pi_A(u)f=\pi_A(u')f'$ 
so that the correspondence 
$[(u,f)]\mapsto\pi_A(u)f$, $D_s\to\pi_A(u)\Hc_B$, 
gives rise to a complex linear structure on $D_s$.
Moreover, 
$$
\Vert[(u,f)]\Vert_{D_s}:=\Vert\pi_A(u)f\Vert_{\Hc_A}
$$
where $[(u,f)]\in D_s$, defines on $D_s$ a Hilbertian norm.

Clearly, this structure does not depend on the choice of $u$. 
Nevertheless, note that the natural
bijection from 
$\Hc_B$ onto the fiber $\Pi^{-1}(s)$ defined by 
$$
\Theta_u\colon\,f\mapsto[(u,f)],\quad\Hc_B\to\Pi^{-1}(s),
$$
is a topological isomorphism but it {\it need not be} an isometry. 
In other words, the fiberwise maps
$$
\Theta_v\Theta^{-1}_u\colon\,[(u,f)]\mapsto f\mapsto[(v,f)],
\quad  D_s\to\Hc_B\to D_t,
$$ 
where $s=uG_B$, $t=vG_B$ and $f\in\Hc_B$, 
are topological isomorphisms but they are not unitary
transformations in general. 
As a complex Hilbert space, $D_s$ has so many realizations 
of the topological dual or predual. We next consider
the following ones.
For $\xi=[(u,f)]$, $\eta=[(v,g)]$ in $D$, and $s=uG_B$, $t=vG_B$, 
we set as in Example~\ref{rel},
$$
\bigl(\xi\mid\eta\bigr)_D\equiv\bigl(\xi\mid\eta\bigr)_{s,t}
:=\bigl(\pi_A(u)f\mid\pi_A(v)g\bigr)_{\Hc_A}.
$$
where $\bigl(\cdot\mid\cdot\bigr)_{\Hc_A}$ is 
the inner product which defines the complex Hilbert structure on
$\Hc_A$ and, by restriction, on $\Hc_B$. 
This is a well-defined, non-negative sesquilinear form on $D$.
In particular 
$\bigl(\cdot\mid\cdot\bigr)_{s,t}
=\overline{\bigl(\cdot\mid\cdot\bigr)}_{t,s}$.
We are mainly interested in forms $\bigl(\cdot\mid\cdot\bigr)_{s,t}$ 
with $t=s^{-*}\in G_A/G_B$. 
In this case
\begin{equation}\label{strong}
\bigl([(u,f)]\mid[(u^{-*},g)]\bigr)_{s,s^{-*}}
=\bigl(\pi_A(u)f\mid\pi_A(u^{-*})g\bigr)_{\Hc_A}
=\bigl(\pi_A(u^{-1})\pi_A(u)f\mid g\bigr)_{\Hc_A}
=\bigl(f\mid g\bigr)_{\Hc_B},
\end{equation}
whenever $[(u,f)]\in D_s$ and [$(u^{-*},g)]\in D_{s^{-*}}$. 
Thus Example~\ref{rel} shows that the {\it basic mapping} 
$$
\Theta\colon\,[(u,f)]\mapsto\pi_A(u)f,\quad D\to\Hc_A,
$$
gives rise to a like-Hermitian structure on the vector bundle $\Pi$.
 
We shall say that $\Pi\colon D\to G_A/G_B$ is the 
(holomorphic) {\it homogeneous like-Hermitian vector bundle} 
associated with the data~$(\pi_A,\pi_B,P)$.
\qed
\end{definition}

\begin{remark}\label{ws}
\normalfont
Let us see that Definition~\ref{homog} is correct, 
that is, condition~(a) of Definition~\ref{like} 
is satisfied.  
In fact, let $u\in G_A$ arbitrary, denote 
$z=uG_B\in G_A/G_B$,  
and let $\varphi$ be any bounded linear functional on $\overline{D}_{z^{-*}}$. 
Then the mapping 
$
\tilde\varphi\colon\,g\mapsto[(u^{-*},g)]\mapsto\varphi([(u^{-*},g)])$, 
$\Hc_B\to\overline{D}_{z^{-*}}\to{\mathbb C}$,  
is antilinear and bounded. 
By the Riesz' theorem there exists $f\in\Hc_B$ such that 
$$
\varphi([(u^{-*},g)])=\tilde\varphi(g)=\bigl(f\mid g\bigr)_{\Hc_B}
\stackrel{\scriptstyle\eqref{strong}}{=}
\bigl([(u,f)]\mid[(u^{-*},g)]\bigr)_{z,z^{-*}}
$$
and so $\bigl(\cdot\mid\cdot\bigr)_{z,z^{-*}}$ is 
a sesquilinear strong duality pairing between 
$D_z$ and $D_{z^{-*}}$.
\qed
\end{remark}

In order to get a better understanding 
of the structures introduced in Definition~\ref{homog}, 
we shall need the following notion. 

\begin{definition}\label{relationship}
\normalfont
Assume that we have the following objects: 
a complex involutive Banach-Lie group~$G$,  
a complex Banach manifold~$Z$ equipped 
with an involutive diffeomorphism $z\mapsto z^{-*}$, 
and a holomorphic like-Hermitian vector bundle 
$\Pi\colon D\to Z$, 
such that $\Pi\circ\alpha=\beta\circ(\id_G\times\Pi)$,
where $\alpha$ and $\beta$ are holomorphic actions of $G$ 
on $D$ and $Z$ and for all $u\in G$ and $z\in Z$ the mapping 
$\alpha(u,\cdot)|_{D_z}\colon D_z\to D_{\beta(u,z)} $
is a bounded linear operator.  
In addition we assume that 
$\beta(u^{-*},z^{-*})=\beta(u,z)^{-*}$ 
whenever $u\in G$ and $z\in Z$  
and we let  $\pi\colon G\to\Bc(\Hc)$ be a holomorphic $*$-representation.  

We say that a holomorphic mapping 
$\Theta\colon D\to\Hc$ \textit{relates} $\Pi$ to $\pi$ 
if it has the following properties: 
\begin{itemize}
\item[{\rm(i)}] for each $z\in Z$ the mapping 
$\Theta_z:=\Theta|_{D_z}\colon D_z\to\Hc$ 
is an injective bounded linear operator and 
in addition we have  
$(\xi\mid\eta)_{z,z^{-*}}=(\Theta(\xi)\mid\Theta(\eta))_{\Hc}$
whenever $\xi\in D_z$ and $\eta\in D_{z^{-*}}$; 
\item[{\rm(ii)}] for every $u\in G$ and $z\in Z$ 
we have 
$\Theta_{\beta(u,z)}\circ\alpha(u,\cdot)|_{D_z}
=\pi(u)\circ\Theta_z\colon D_z\to\Hc$. 
\end{itemize}
\qed
\end{definition}

Now we turn to a result (Theorem~\ref{compare}) 
which points out that the basic mapping $\Theta$ introduced in 
Definition~\ref{homog} indeed plays a central role in the whole picture. 
In this statement we denote by $i_0^*(\cdot)$ the pull-back 
of a bundle by the mapping $i_0$; 
see for instance \cite{La01} for some details. 

\begin{theorem}\label{compare}
In the setting of {\rm Definition~\ref{relationship}},  
let $z_0\in Z$ such that $z_0^{-*}=z_0$,  
and assume that the isotropy group 
$G_0:=\{u\in G\mid\beta(u,z_0)=z_0\}$
is a Banach-Lie subgroup of $G$. 
In addition assume that the orbit 
of $z_0$, that is, 
$\Oc_{z_0}=\{\beta(u,z_0)\mid u\in G\}$, 
is a submanifold of $Z$, 
and denote by $i_0\colon\Oc_{z_0}\hookrightarrow Z$ the corresponding 
embedding map. 
Then there exists a closed subspace $\Hc_0$ of $\Hc$ 
such that the following assertions hold: 
\begin{itemize}
\item[{\rm(i)}] For every $u\in G_0$ we have $\pi(u)\Hc_0\subseteq\Hc_0$. 
\item[{\rm(ii)}] Denote by $\pi_0\colon\,u\mapsto\pi(u)|_{\Hc_0}$,
$G_0\to\Bc(\Hc_0)$, 
the corresponding representation of $G_0$ on $\Hc_0$, 
by $\Pi_0\colon D_0\to G/G_0$ the like-Hermitian vector bundle 
associated with 
the data $(\pi,\pi_0,P_{\Hc_0})$, 
and by $\Theta_0\colon D_0\to\Hc$ the basic mapping associated with 
the data $(\pi,\pi_0,P_{\Hc_0})$. 
Then there exists a biholmorphic bijective $G$-equivariant map 
$\theta\colon D_0\to i_0^*(D)$ such that 
$\theta$ sets up an isometric isomorphism of like-Hermitian 
vector bundles over $G/G_0\simeq\Oc_{z_0}$ and the diagram 
$$
\begin{CD}
D_0 @>{\theta}>> i_0^*(D) \\
@V{\Theta_0}VV @VV{\Theta|_{i_0^*(D)}}V \\
\Hc @>{\id_{\Hc}}>> \Hc
\end{CD}
$$
is commutative.
\end{itemize} 
\end{theorem}

\begin{proof}
For $\eta,\xi\in D_{z_0}$ we have 
$$
(\Theta_{z_0}^{-*}(\Theta_{z_0}\xi)\mid\eta)_{z_0,z_0}
=(\Theta_{z_0}\xi \mid \Theta_{z_0}\eta)_{\Hc}=(\xi \mid \eta)_{z_0,z_0},
$$
where the first equality is derived from Remark~\ref{sesqui1} and the second one is the
hypothesis of Definition~\ref{relationship} (i), since $z_0=z_0^{-*}$. As 
$(\cdot\mid\cdot)_{z_0,z_0}$ is a strong duality pairing we have 
$\Theta_{z_0}^{-*}(\Theta_{z_0}\xi)=\xi$. 
This implies that $\Ran(\Theta_{z_0})$ is closed in 
$\Hc$ since $\Theta_{z_0}^{-*}$ is bounded. Put $\Hc_0:=\Ran(\Theta_{z_0})$. 

For arbitrary $u\in G_0$ we have $\beta(u,z_0)=z_0$. 
Then property~(ii) in Definition~\ref{relationship} shows that 
we have a commutative diagram 
$$
\begin{CD}
D_{z_0} @>{\alpha(u,\cdot)|_{D_{z_0}}}>> D_{z_0} \\
@V{\Theta_{z_0}}VV @VV{\Theta_{z_0}}V \\
\Hc @>{\pi(u)}>> \Hc
\end{CD}
$$
whence $\pi(u)(\Theta_{z_0}(D_{z_0}))\subseteq \Theta_{z_0}(D_{z_0})$, 
that is, $\pi(u)\Hc_0\subseteq\Hc_0$. 
Thus $\Hc_0$ has the desired property~(i). 

To prove~(ii) we first note that, since $G_0$ is a Banach-Lie subgroup of $Z$, 
it follows that the $G$-orbit $\Oc_{z_0}\simeq G/G_0$ 
has a natural structure of Banach homogeneous space of $G$ 
(in the sense of \cite{Rae77}) such that the inclusion map 
$i_0\colon\Oc_{z_0}\hookrightarrow Z$ is an embedding. 

Next define 
\begin{equation}\label{tilde}
\widetilde{\theta}\colon G\times\Hc_0\to D,\quad 
\widetilde{\theta}(u,f):=
\alpha(u,\Theta_{z_0}^{-1}(f))=
\Theta_{\beta(u,z_0)}^{-1}(\pi(u)f)
\in D_{\beta(u,z_0)}\subseteq D,
\end{equation}
where the equality follows by property~(ii) in Definition~\ref{relationship}. 
Then for all $u\in G$, $u_0\in G_0$, and $f\in\Hc_0$ 
we have 
$\beta(uu_0^{-1},z_0)=\beta(u,z_0)$ and 
$$
\widetilde{\theta}(uu_0^{-1},\pi(u_0)f)
=\Theta_{\beta(uu_0^{-1},z_0)}^{-1}(\pi(uu_0^{-1})\pi(u_0)f) 
=\Theta_{\beta(u,z_0)}^{-1}(\pi(u)f)
=\widetilde{\theta}(u,f).
$$
In particular there exists a well defined map 
$$
\theta\colon\,[(u,f)]\mapsto\Theta_{\beta(u,z_0)}^{-1}(\pi(u)f),
\quad  
 G\times_{G_0}\Hc_0\to D.
$$
This mapping is $G$-equivariant with respect to 
the actions of $G$ on $G\times_{G_0}\Hc$ and on $D$ 
since $\widetilde{\theta}$ is $G$-equivariant:
for all $u,v\in G$ and $f\in\Hc_0$ 
we have 
$$
\widetilde{\theta}(uv,f)
=\Theta_{\beta(uv,z_0)}^{-1}(\pi(uv)f)  
=\Theta_{\beta(u,\beta(v,z_0))}^{-1}(\pi(u)\pi(v)f) 
=\alpha\bigl(u,\Theta_{\beta(v,z_0)}^{-1}(\pi(v)f)\bigr) 
=\alpha(u,\widetilde{\theta}(v,f)), 
$$
where the second equality follows since 
$\beta\colon G\times Z\to Z$ is a group action, 
while the third equality is a consequence of 
property~(ii) in Definition~\ref{relationship}. 
Besides, it is clear that $\theta$ is a bijection onto $i_0^*(D)$ and 
a fiberwise isomorphism. 
Also it is clear from the above construction of $\theta$ 
and from the definition 
of the basic mapping $\Theta_0\colon D_0\to\Hc$ associated with 
the data $(\pi,\pi_0,P_{\Hc_0})$
(see Definition~\ref{homog}) that $\Theta\circ\theta=\Theta_0$, 
that is, the diagram in the statement is indeed commutative. 
In addition, since both mappings $\Theta$ and $\Theta_0$ 
are fiberwise ``isometric'' 
(see property~(i) in Definition~\ref{relationship} 
above and Definition~\ref{homog}), 
it follows by $\Theta\circ\theta=\Theta_0$ that 
$\theta$ gives us an isometric morphism of like-Hermitian 
bundles over $G/G_0\simeq\Oc_{z_0}$. 

Now we still have to prove that the map 
$\theta\colon D_0=G\times_{G_0}\Hc_0\to i_0^*(D)\subseteq D$ is biholomorphic. 
We first show that it is holomorphic. 
Since $\Oc_{z_0}$ is a submanifold of $Z$, it 
follows that $i_0^*(D)$ is a submanifold of $D$ 
(see for instance the comments after Proposition~1.4 
in Chapter~III of \cite{La01}). 

Thus it will be enough to show that $\theta\colon G\times_{G_0}\Hc_0\to D$ 
is holomorphic. 
And this property is equivalent (by Corollary~8.3(ii) in \cite{Up85}) 
to the fact that 
the mapping $\widetilde{\theta}\colon G\times\Hc_0\to D$ 
is holomorphic, 
since the natural projection $G\times\Hc_0\to G\times_{G_0}\Hc_0$ 
is a holomorphic submersion. 
Now the fact that $\widetilde{\theta}\colon G\times\Hc_0\to D$ 
is a holomorphic map follows by the first formula 
in its definition~\eqref{tilde}, since the group action 
$\alpha\colon G\times D\to D$ is holomorphic. 

Consequently the mapping $\theta\colon G\times_{G_0}\Hc_0\to i_0^*(D)$ 
is holomorphic. 
Then the fact that the inverse 
$\theta^{-1}\colon i_0^*(D)\to G\times_{G_0}\Hc_0$ 
is also holomorphic follows by 
general arguments in view of the following facts
(the first and the second of them have been already established, 
and the third one is well-known): 
Both $G\times_{G_0}\Hc_0$ and 
$i_0^*(D)$ are locally trivial holomorphic vector bundles; 
we have a commutative diagram 
$$
\begin{CD}
G\times_{G_0}\Hc_0 @>{\theta}>> i_0^*(D) \\
@VVV @VVV \\
G/G_0 @>>> \Oc_{z_0}
\end{CD}
$$ 
where the bottom arrow is the biholomorphic map 
$G/G_0\simeq\Oc_{z_0}$ induced by the action $\beta\colon G\times Z\to Z$, 
and the vertical arrows are the projections of the corresponding holomorphic 
like-Hermitian vector bundles; 
the inversion mapping is holomorphic on the open set of invertible 
operators on a complex Hilbert space. 

The proof is completed.   
\end{proof}

\begin{remark}\label{nonhomog}
\normalfont
The significance of Theorem~\ref{compare} is the following one: 
In the setting of Definition~\ref{relationship}, the special situation 
of Definition~\ref{homog} is met precisely when the action 
$\beta\colon G\times Z\to Z$ is transitive, 
and in this case the basic mapping  
is essentially the unique mapping that relates the bundle~$\Pi$ 
to the representation of the bigger group~$G$. 

On the other hand, by considering direct products of 
homogeneous Hermitian vector bundles, we can construct obvious examples 
of other maps relating bundles to representations as in 
Definition~\ref{relationship}. 
\qed
\end{remark}

\section{Reproducing $(-*)$-kernels}\label{sect3}

\begin{definition}\label{kernel}
\normalfont 
Let $\Pi\colon D\to Z$ be a like-Hermitian bundle, with involution $-*$ in $Z$. 
A {\it reproducing $(-*)$-kernel} on $\Pi$ is a section  
$$K\in\Gamma(Z\times Z,\Hom(p_2^*\Pi,p_1^*\Pi))$$
(whence $K(s,t)\colon D_t\to D_s$ for all $s,t\in Z$) 
which is ($-*$)-positive definite in the following sense: 
For every $n\ge1$ and $t_j\in Z$,
$\eta_j^{-*}\in D_{t_j^{-*}}$ ($j= 1,\dots, n$),
$$
\sum_{j,l=1}^n
\bigl(\eta_j^{-*}\mid K(t_j,t_l^{-*})\eta_l^{-*}\bigr)_{t_j^{-*},t_j}
=\sum_{j,l=1}^n
\bigl(K(t_l,t_j^{-*})\eta_j^{-*}\mid \eta_l^{-*}\bigr)_{t_l,t_l^{-*}}\ge0.
$$ 
Here $p_1,p_2\colon Z\times Z\to Z$ are the natural projection mappings. 
If in addition $\Pi\colon D\to Z$ is 
a holomorphic like-Hermitian vector bundle  
and  $K(\cdot,t)\eta\in\Oc (Z,D)$ for all $\eta\in D_t$ and 
$t\in Z$, 
then we say that $K$ is a {\it holomorphic} reproducing $(-*)$-kernel.
\qed
\end{definition}

\begin{remark}\label{obs}
\normalfont
In Definition~\ref{kernel}, the symbol $\eta_j^{-*}$ 
is just a way to refer to elements of $D_{t_j^{-*}}$, 
that is, 
$\eta_j^{-*}$ is not associated to any element $\eta_j$ of 
$D_{t_j}$ necessarily.
From the definition we have that $K(s,s^{-*})\ge0$ in the sense that 
$\bigl(K(s,s^{-*})\xi^{-*}\mid\xi^{-*}\bigr)_{s,s^{-*}}\ge0$ 
for all $\xi^{-*}\in D_{s^{-*}}$.
\qed
\end{remark}

The following results are related to the extension of 
Theorem~4.2 in \cite{BR07} to reproducing kernels 
on like-Hermitian vector bundles. 

\begin{proposition}\label{rkhs1}
Let $\Pi\colon D\to Z$ be a smooth like-Hermitian
vector bundle and, as usually,  
denote by 
$p_1,p_2\colon Z\times Z\to Z$ the projections. 
Next consider a section 
$K\in\Gamma(Z\times Z,\Hom(p_2^*\Pi,p_1^*\Pi))$ 
and for all $s\in Z$ and $\xi\in D_s$ denote 
$K_\xi=K(\cdot,s)\xi\in\Gamma(Z,D)$. 
Also denote 
$$\Hc^K_0:=\spann_{\mathbb C}\{K_\xi\mid\xi\in D\}\subseteq\Gamma(Z,D).$$ 
Then $K$ is a reproducing $(-*)$-kernel on $\Pi$ if and only
if 
there exists a complex Hilbert space $\Hc$ such that 
$\Hc^K_0$ is a dense linear subspace of $\Hc$ and 
\begin{equation}\label{prods}
(K_\eta\mid K_\xi)_{\Hc}=(K(s^{-*},t)\eta\mid\xi)_{s^{-*},s} 
\end{equation}
whenever $s,t\in Z$, $\xi\in D_{s}$, and $\eta\in D_t$.
\end{proposition}

\begin{proof}
First assume that there exists a Hilbert space as in the statement. 
Then for all $n\ge1$ and $t_j\in Z$,
$\eta_j\in D_{t_j^{-*}}$ ($j= 1,\dots, n$), 
it follows by \eqref{prods} that 
$$
\sum_{j,l=1}^n
\bigl(K(t_l,t_j^{-*})\eta_j\mid \eta_l\bigr)_{t_l,t_l^{-*}}
=\sum_{j,l=1}^n (K_{\eta_j}\mid K_{\eta_l})_{\Hc}
=\Bigl(\sum_{j=1}^n K_{\eta_j}\mid\sum_{l=1}^n K_{\eta_l}\Bigr)_{\Hc}
\ge0.
$$ 
In addition, for arbitrary $s,t\in Z$, $\xi\in D_s$ and $\eta\in D_t$, 
we have  
$$(\eta\mid K(t^{-*},s)\xi)_{t,t^{-*}}
=\overline{(K(t^{-*},s)\xi\mid\eta)}_{t^{-*},t}
=\overline{(K_\xi\mid K_\eta)}_{\Hc}
=(K_\eta\mid K_\xi)_{\Hc}
=(K(s^{-*},t)\eta\mid\xi)_{s^{-*},s},  
$$
where the second equality and the fourth one follow 
by \eqref{prods}. 

Conversely, let us assume that $K$ is a reproducing $(-*)$-kernel. 
We are going to define a positive Hermitian sesquilinear form 
on $\Hc^K_0$ by 
\begin{equation}\label{hidden}
(\Theta\mid\Delta)_{\Hc}=\sum_{j,l=1}^n
(K(s_l^{-*},t_j)\eta_j\mid\xi_l)_{s_l^{-*},s_l}
\end{equation}
for $\Theta,\Delta\in\Hc^K_0$ of the form 
$\Theta=\sum\limits_{j=1}^n K_{\eta_j}$ and 
$\Delta=\sum\limits_{l=1}^n K_{\xi_l}$, 
where $\eta_j\in D_{t_j}$, $\xi_l\in D_{s_l}$, 
and $t_j,s_l\in Z$ for $j,l=1,\dots,n$. 
The assumption that $K$ is a reproducing $(-*)$-kernel implies at once 
that for all $\Theta,\Delta\in\Hc^K_0$ 
we have $(\Theta\mid\Theta)_{\Hc}\ge0$ 
and $\overline{(\Theta\mid\Delta)}_{\Hc}=(\Delta\mid\Theta)_{\Hc}$. 
To see that $(\cdot\mid\cdot)_{\Hc}$ is well defined, 
note that 
it is clearly sesquilinear and \eqref{hidden} implies 
\begin{equation}\label{point}
\overline{(\Delta\mid\Theta)}_{\Hc}
=(\Theta\mid\Delta)_{\Hc}
=\sum_{l=1}^n (\Theta(s_l^{-*})\mid\xi_l)_{s_l^{-*},s_l},
\end{equation}
hence 
$(\Delta\mid\Theta)_{\Hc}=(\Theta\mid\Delta)_{\Hc}=0$ 
if it happens that $\Theta=0$. 
This implies that $(\cdot\mid\cdot)_{\Hc}$ is well defined, 
and the above remarks show that this is a 
nonnegative Hermitian sesquilinear form on $\Hc^K_0$. 

To check that $(\cdot\mid\cdot)_{\Hc}$ is also non-degenerate, 
let $\Theta\in\Hc^K_0$ such that $(\Theta\mid\Theta)_{\Hc}=0$. 
Since $(\cdot\mid\cdot)_{\Hc}$ is a nonnegative Hermitian sesquilinear form, 
it follows that it satisfies the Cauchy-Schwarz inequality, 
hence for all $\Delta\in\Hc^K_0$ we have
$\vert(\Theta\mid\Delta)_{\Hc}\vert 
\le(\Theta\mid\Theta)_{\Hc}^{1/2}(\Delta\mid\Delta)_{\Hc}^{1/2}=0$, 
whence
$(\Theta\mid\Delta)_{\Hc}=0$. 
It follows by this property along with the formula \eqref{point} 
that for arbitrary $s\in Z$ and $\xi\in D_{s^{-*}}$ we have 
$(\Theta(s)\mid\xi)_{s,s^{-*}}=(\Theta\mid K_{\xi})_{\Hc}=0$.
Since $\{(\cdot\mid\cdot)_{z,z^{-*}}\}_{z\in Z}$ 
is a like-Hermitian structure, it then follows that 
$\Theta(s)=0$ for all $s\in Z$, hence $\Theta=0$. 
Consequently $(\cdot\mid\cdot)_{\Hc}$ is a scalar product 
on $\Hc^K_0$, 
and then the completion of $\Hc^K_0$ with respect to 
this scalar product is a complex Hilbert space 
with the asserted properties. 
\end{proof}

\begin{definition}\label{rkhs2}
\normalfont
Let $\Pi\colon D\to Z$ be a smooth like-Hermitian
vector bundle, 
$p_1,p_2\colon Z\times Z\to Z$ the projections, and  
let 
$K\in\Gamma(Z\times Z,\Hom(p_2^*\Pi,p_1^*\Pi))$ 
be a reproducing $(-*)$-kernel. As above, for all $s\in Z$ and $\xi\in D_s$, put 
$K_\xi=K(\cdot,s)\xi\in\Gamma(Z,D)$. 
It is clear that the Hilbert space $\Hc$ given by Proposition~\ref{rkhs1} 
is uniquely determined. 
We shall denote it by $\Hc^K$ 
and we shall call it the {\it reproducing $(-*)$-kernel Hilbert space} 
associated with $K$. 

In the same framework we also define the mapping
\begin{equation}\label{coh}
\widehat{K}\colon D\to\Hc^K, \quad \widehat{K}(\xi)=K_\xi.
\end{equation}
It follows by Lemma~\ref{rkhs3a} below that 
 for every $s\in Z$  there exists a bounded linear operator 
$\theta_s\colon\Hc^K\to D_{s^{-*}}$ such that 
\begin{equation}\label{2}
(\forall\xi\in D_s,h\in\Hc^K)\quad 
(\widehat{K}(\xi)\mid h)_{\Hc^K}=(\xi\mid\theta_sh)_{s,s^{-*}}.
\end{equation}
Note that the operator $\theta_s$ is uniquely determined since 
$\{(\cdot\mid\cdot)_{z,z^{-*}}\}_{z\in Z}$ 
is a like-Hermitian structure, 
and in the notation of Remark~\ref{sesqui1} we have 
\begin{equation}\label{2and1/2}
(\theta_s)^{-*}=\widehat{K}|_{D_{s^{-*}}}
\end{equation}
\qed
\end{definition}

\begin{lemma}\label{rkhs3a}
Assume the setting of {\rm Definition~\ref{rkhs2}}. 
Then for every $s\in Z$ the operator 
$\widehat{K}|_{D_s}\colon D_s\to\Hc^K$ 
is bounded, linear and adjointable, in the sense 
that there exists a bounded linear operator 
$\theta_s\colon\Hc^K\to D_{s^{-*}}$ such that 
\eqref{2} is satisfied. 
\end{lemma}

\begin{proof}
Since at every point of $Z$ we have a sesquilinear strong duality pairing,  
it will be enough to show that for arbitrary $s\in Z$ 
the linear operator $\widehat{K}|_{D_s}\colon D_s\to\Hc^K$ 
is continuous.
(See Remark~\ref{sesqui1}.) 
To this end, let us denote by 
$\Vert\cdot\Vert_{D_s}$ 
any norm that defines the topology of the fiber $D_s$. 
Then for every $\xi\in D_s$ we have 
$\Vert\widehat{K}(\xi)\Vert_{\Hc^K}
=\Vert K_\xi\Vert_{\Hc^K}
=(K_\xi\mid K_\xi)_{\Hc^K}^{1/2}
\stackrel{\eqref{prods}}{=}(K(s^{-*},s)\xi\mid\xi)_{s^{-*},s}^{1/2}
\le M_s^{1/2}\Vert\xi\Vert_{D_s}$,  
where $M_s>0$ denotes the norm of the continuous sesquilinear functional 
$D_s\times D_s\to{\mathbb C}$ defined by  
$(\xi,\eta)\mapsto (K(s^{-*},s)\xi\mid\eta)_{s^{-*},s}$. 
So the operator $\widehat{K}|_{D_s}\colon D_s\to\Hc^K$  
is indeed bounded and
$\Vert\widehat{K}|_{D_s}\Vert\le M_s^{1/2}$. 
\end{proof}

\begin{example}\label{rkhs3}
\normalfont
Every reproducing kernel on a Hermitian vector bundle 
(see e.g., Section~4 in \cite{BR07}) 
provides an illustration for Definition~\ref{rkhs2}. 
In fact, this follows  
since every Hermitian vector bundle 
is like-Hermitian. 
\qed
\end{example}

\begin{proposition}\label{rkhs5}
Let $\Pi\colon D\to Z$ be a like-Hermitian bundle, 
and denote by $p_1,p_2\colon Z\times Z\to Z$ the natural projections. 
Then for every reproducing $(-*)$-kernel 
$K\in\Gamma(Z\times Z,\Hom(p_2^*\Pi,p_1^*\Pi))$ 
there exists a unique linear mapping 
$\iota\colon\Hc^K\to\Gamma(Z,D)$
with the following properties:
\begin{itemize}
\item[{\rm(a)}] The restriction of $\iota$ to the dense subspace 
$\Hc^K_0$ is the identity mapping. 
\item[{\rm(b)}] The mapping $\iota$ is injective. 
\item[{\rm(c)}] The evaluation operator 
$\ev^\iota_s\colon\,h\mapsto\bigl(\iota(h)\bigr)(s)$, 
$\Hc^K\to D_s$, 
is continuous linear for arbitrary $s\in Z$, and we have
$$(\forall s,t\in Z)\quad K(s,t^{-*})=\ev^\iota_s\circ(\ev^\iota_t)^{-*}.$$
\end{itemize}
\end{proposition}

\begin{definition}\label{realiz_oper}
\normalfont
In the setting of Proposition~\ref{rkhs5} 
we shall say that $\iota$ is the \textit{realization operator} 
associated with the reproducing $(-*)$-kernel $K$.  
\qed
\end{definition}

\begin{proof}[Proof of Proposition~\ref{rkhs5}]  
The uniqueness of $\iota$ is clear. 
To prove the existence of $\iota$, 
note that for every $s\in Z$ there exists 
a bounded linear operator $\theta_s\colon\Hc^K\to D_{s^{-*}}$ 
such that 
\begin{equation}\label{3}
(\forall\xi\in D_s,h\in\Hc^K)\quad 
(K_\xi\mid h)_{\Hc^K}=(\xi\mid\theta_sh)_{s,s^{-*}}
\end{equation}
(see Lemma~\ref{rkhs3a}). 
We shall define the wished-for mapping $\iota$ by 
\begin{equation}\label{4}
\iota\colon\Hc^K\to\Gamma(Z,D),\quad
\bigl(\iota(h)\bigr)(s):=\theta_{s^{-*}}h
\end{equation}
whenever $h\in\Hc^K$ and $s\in Z$. 
In particular we have 
\begin{equation}\label{5}
(\forall s\in Z)\quad \ev^\iota_s=\theta_{s^{-*}},
\end{equation}
and in addition equation~\eqref{2and1/2} holds.

It is also clear that the mapping $\iota$ defined by \eqref{4} is linear. 
To prove that it is injective, let $h\in\Hc^K$ with $\iota(h)=0$. 
Then $(\iota(h))(s^{-*})=0$ for all $s\in Z$, 
so that $\theta_sh=0$ for all $s\in Z$, according to \eqref{4}. 
Now~\eqref{3} shows that 
$(K_\xi\mid h)_{\Hc^K}=0$ for all $\xi\in D$, 
whence $h\perp\Hc^K_0$ in $\Hc^K$. 
Since $\Hc^K_0$ is a dense subspace of $\Hc^K$, 
it then follows that $h=0$. 

We shall check that the restriction of $\iota$ 
to $\Hc^K_0$ is the identity mapping. 
To this end it will be enough to see that 
for all $t\in Z$ and $\eta\in D_t$ we have $\iota(K_\eta)=K_\eta$. 
In fact, at any point $s\in Z$ we have 
$(\iota(K_\eta))(s)=\theta_{s^{-*}}(K_\eta)$ 
by \eqref{4}. 
Hence for all $\xi\in D_{s^{-*}}$ we get 
$$\begin{aligned}
(\xi\mid(\iota(K_\eta))(s))_{s^{-*},s}
 &=(\xi\mid\theta_{s^{-*}}(K_\eta))_{s^{-*},s} 
 \stackrel{\eqref{3}}{=}(K_\xi\mid K_\eta)_{\Hc^K} 
 \stackrel{\eqref{prods}}{=}(K(t^{-*},s^{-*})\xi\mid\eta)_{t^{-*},t} \\
 &=(\xi\mid K(s,t)\eta)_{s^{-*},s} 
  =(\xi\mid K_\eta(s))_{s^{-*},s}.
\end{aligned}$$
Since $\xi\in D_{s^{-*}}$ is arbitrary and 
$\{(\cdot\mid\cdot)_{z,z^{-*}}\}_{z\in Z}$ is 
a like-Hermitian structure, 
it then follows that $(\iota(K_\eta))(s)=K_\eta(s)$ for all $s\in Z$, 
whence $\iota(K_\eta)=K_\eta$, as desired. 

Next we shall prove that $\iota$ has the asserted property~(c). 
To this end, let $s,t\in Z$, $\eta\in D_{t^{-*}}$, and $\xi\in D_{s^{-*}}$ 
arbitrary. 
Then 
$$
\begin{aligned}
((\ev^\iota_s\circ(\ev^\iota_t)^{-*})\eta\mid\xi)_{s,s^{-*}}
 &\stackrel{\eqref{5}}{=}
  ((\theta_{s^{-*}}\circ(\theta_{t^{-*}})^{-*})\eta\mid\xi)_{s,s^{-*}}
  \stackrel{\eqref{2}}{=}
  (((\theta_{t^{-*}})^{-*})\eta\mid K_\xi)_{\Hc^K}
  \stackrel{\eqref{2and1/2}}{=}
   (K_\eta\mid K_\xi)_{\Hc^K}\\
 &\stackrel{\eqref{prods}}{=}
  (K(s,t^{-*})\eta\mid \xi)_{s,s^{-*}}.
\end{aligned}
$$
Since $\eta\in D_{t^{-*}}$ and $\xi\in D_{s^{-*}}$ 
are arbitrary and $\{(\cdot\mid\cdot)_{z,z^{-*}}\}_{z\in Z}$ is 
a like-Hermitian structure, 
it follows that 
$\ev^\iota_s\circ(\ev^\iota_t)^{-*}=K(s,t^{-*})$ 
for arbitrary $s,t\in Z$, as desired. 
\end{proof}

We now extend to our framework some basic properties of 
the classical reproducing kernels 
(see for instance the first chapter of \cite{Ne00}). 

\begin{proposition}\label{cont}
Assume that $\Pi\colon D\to Z$ is a like-Hermitian vector bundle, 
and $K$ is a continuous reproducing $(-*)$-kernel on $\Pi$ 
with the realization operator $\iota\colon\Hc^K\to\Gamma(Z,D)$.
Then the following assertions hold:
\begin{enumerate}
\item[{\rm(a)}] We have $\Ran\iota\subseteq\Cc(Z,D)$ and 
the mapping $\iota$ is continuous with respect to 
the topology of $\Cc(Z,D)$ defined by the uniform convergence on 
the compact subsets of $Z$. 
\item[{\rm(b)}] If $\Pi$ is a holomorphic bundle and 
$K$ is a holomorphic reproducing $(-*)$-kernel 
then we have $\Ran\iota\subseteq\Oc(Z,D)$.
\end{enumerate}
\end{proposition}

\begin{proof}
The proof has two stages. 

$1^\circ$ At this stage we prove that every $s\in Z$ 
has an open neighborhood $V_s$ such that 
for every sequence $\{h_n\}_{n\in{\mathbb N}}$ in $\Hc^K$ 
convergent to some $h\in\Hc^K$ we have 
$\lim\limits_{n\in{\mathbb N}}\bigl(\iota(h_n)\bigr)(z)
=\bigl(\iota(h)\bigr)(z)$ 
uniformly for $z\in V_s$. 

In fact, since the vector bundle $\Pi$ is locally trivial, 
there exists an open neighborhood $V$ of $s$ 
such that $\Pi$ is trivial over both 
$V$ and $V^{-*}:=\{z^{-*}\mid z\in V\}$. 
Let 
$\Psi_0\colon V\times \Ec\to\Pi^{-1}(V)$ 
and $\Psi_1\colon V^{-*}\times \Ec\to\Pi^{-1}(V^{-*})$ 
be trivializations of the vector bundle $\Pi$ 
over $V$ and $V^{-*}$ 
respectively, 
where the Banach space $\Ec$ is the typical fiber of $\Pi$. 
In particular, these trivializations 
allow us to endow each fiber $D_z$  with a norm 
(constructed out of the norm of $\Ec$)
for $z\in V\cup V^{-*}$.  
On the other hand, property~(c) in Definition~\ref{like} shows that 
the function 
$$
B\colon\,
(z,x,y)\mapsto(\Psi_0(z,x)\mid\Psi_1(z^{-*},y))_{z,z^{-*}}, 
\quad V\times\Ec\times\Ec\to{\mathbb C} 
$$
is smooth. 
Then by property~(c) in Definition~\ref{like} we get a well-defined mapping
$$
\widetilde{B}\colon V\to\Iso(\Ec,\overline{\Ec}^*),
\quad \widetilde{B}(z)x:=B(z,x,\cdot) \text{ for }z\in V\text{ and }x\in\Ec, 
$$
and it is straightforward to prove that $\widetilde{B}$ is continuous 
since $B$ is so.  
Here $\Iso(\Ec,\overline{\Ec}^*)$ stands for the 
set of all topological isomorphisms $\Ec\to\overline{\Ec}^*$, 
which is an open subset of the complex Banach space $\Bc(\Ec,\overline{\Ec}^*)$.  
Then by shrinking the open neighborhood $V$ of $s$ we may 
assume that there exists $M>0$ such that 
$\max\{\Vert\widetilde{B}(z)\Vert,\Vert\widetilde{B}(z)^{-1}\Vert\}<M$ 
whenever $z\in V$. 
In particular, for such $z$  and $x\in\Ec$ we have 
$\Vert x\Vert<M\Vert\widetilde{B}(z)x\Vert$, 
and then the definition of the norm of $\widetilde{B}(z)x\in\overline{\Ec}^*$ 
implies the following fact:
$$
(\forall z\in V)(\forall x\in\Ec)(\exists y\in\Ec,\;\Vert y\Vert=1)\quad 
\Vert x\Vert\le M\vert B(z,x,y)\vert.
$$
In view of the fact that the norms of the fibers $D_z$ and $D_{z^{-*}}$ 
are defined such that the operators 
$\Psi_0(z,\cdot)\colon\Ec\to D_z$ and 
$\Psi_1(z^{-*},\cdot)\colon\Ec\to D_{z^{-*}}$ 
are isometries whenever $z\in V$, it then follows that 
\begin{equation}\label{star}
(\forall z\in V)(\forall\eta\in D_z)(\exists\xi\in D_{z^{-*}},\;
\Vert\xi\Vert=1)\quad 
\Vert\eta\Vert\le M\vert(\xi\mid\eta)_{z^{-*},z}\vert.
\end{equation}
On the other hand, it follows by \eqref{4} that 
$\Vert\bigl(\iota(h)\bigr)(z)\Vert_{D_z}
=\Vert\theta_{z^{-*}}(h)\Vert_{D_z}$ for arbitrary 
$z\in V$ and $h\in\Hc^K$. 
Then by~\eqref{star} there exists $\xi\in D_{z^{-*}}$ such that 
$\Vert\xi\Vert=1$ and 
$$\Vert\bigl(\iota(h)\bigr)(z)\Vert_{D_z}
\le M \vert(\xi\mid\theta_{z^{-*}}(h))_{z^{-*},z}\vert
\mathop{=}\limits^{\eqref{3}}
M\vert(K_\xi\mid h)_{\Hc^K}\vert
\le M\Vert K_\xi\Vert_{\Hc^K}\Vert h\Vert_{\Hc^K}.$$
On the other hand, 
since $K\colon Z\times Z\to\Hom(p_2^*\Pi,p_1^*\Pi)$ 
is continuous, 
it follows that after shrinking again the neighborhood $V$ 
of $s$ we may suppose that $m:=\sup\limits_{z\in V}M_z<\infty$, 
where $M_z$ denotes the norm of 
the bounded sesquilinear functional 
$D_z\times D_z\to{\mathbb C}$ defined by  
$(\eta_1,\eta_2)\mapsto (K(z^{-*},z)\eta_1\mid\eta_2)_{z^{-*},z}$ 
whenever $z\in V$. 
Then the computation from the proof of Lemma~\ref{rkhs3a} 
shows that 
$\Vert K_\xi\Vert_{\Hc^K}\le m^{1/2}\Vert\xi\Vert_{D_z}=m^{1/2}$. 
It then follows by the above inequalities that 
we end up with an open neighborhood $V$ of $s$ with the following property: 
$$(\forall h\in\Hc^K)(\forall z\in V)\quad 
\Vert\bigl(\iota(h)\bigr)(z)\Vert_{D_z}
\le m^{1/2}M\Vert h\Vert_{\Hc^K}, $$
which clearly implies the claim from the beginning of the present stage 
of the proof. 

$2^\circ$ At this stage we come back to the proof of the assertions (a)--(b). 
Assertion~(a) follows by what we proved at stage~$1^\circ$ 
by means of a straightforward compactness reasoning  
and by what we proved at stage~$1^\circ$,  
since $K_\xi\in\Cc(Z,D)$ whenever $\xi\in D$ 
and $\spann_{\mathbb C}\{K_\xi\mid\xi\in D\}=\Hc_0^K$. 
Finally, assertion~(b) follows by the assertion (a)  
in a similar manner,  
since $\Oc(Z,D)$ is a closed subspace of $\Cc(Z,D)$ 
with respect to the topology of uniform convergence on 
the compact subsets of $Z$ 
(see Corollary~1.14 in \cite{Up85}). 
\end{proof}

\begin{remark}\label{schwartz}
\normalfont
It follows by Proposition~\ref{cont} (a) that every reproducing $(-*)$-kernel 
Hilbert space $\Hc^K$ is a Hilbert subspace of $\Cc(Z,D)$ in the sense of 
\cite{Sc64}. 
Thus the theory of reproducing $(-*)$-kernels developed in the present section 
provides a new class of examples of reproducing kernels in the sense of 
Laurent Schwartz. 
\qed
\end{remark}

\section{Homogeneous Like-Hermitian vector bundles and kernels}\label{sect4}

We develop here some spects of the theory of kernels introduced in the previous section, when the
manifold $Z$ is assumed to be a homogeneous manifold arising from the (smooth) action of a
Banach-Lie group, as in Definition~\ref{homog}. 
Specifically, we shall construct realizations
of $*$-representations, of Banach-Lie groups, on spaces of analytic sections in 
like-Hermitian vector bundles. 
A critical role in this connection will be played by 
the following class of examples of reproducing $(-*)$-kernels 
(compare the special case discussed in Example~\ref{universal}). 

\begin{example}\label{equivariant}
\normalfont
Assume the setting of Definition~\ref{homog}: Let $G_A$ be an involutive real 
(respectively, complex) Banach-Lie group and $G_B$ 
an involutive real (respectively, complex) Banach-Lie subgroup of $G_A$. 
For $X=A$ or $X=B$, let $\Hc_X$ be a complex Hilbert space 
with $\Hc_B$ closed subspace in $\Hc_A$ and 
$P\colon\Hc_A\to\Hc_B$ the corresponding orthogonal projection, 
and let 
$\pi_X\colon G_X\to \Bc(\Hc_X)$ be a uniformly continuous 
(respectively, holomorphic) 
$*$-representations such that 
$\pi_B(u)=\pi_A(u)|_{\Hc_B}$ for all $u\in G_B$. 
In addition, denote by 
$\Pi\colon D=G_A\times_{G_B}\Hc_B\to G_A/G_B$ 
the homogeneous like-Hermitian vector bundle 
associated with 
the data~$(\pi_A,\pi_B,P)$, and let $p_1,p_2\colon G_A/G_B\times G_A/G_B\to G_A/G_B$ 
be the natural projections.
Set
$$
K(s,t)\eta=[(u,P(\pi_A(u^{-1})\pi_A(v)f))]
$$
for $s,t\in G_A/G_B$, $s=u G_B$, $t=v G_B$, 
and $\eta=[(v,f)]\in D_t\subset D$. 
Then {\it $K$ is a reproducing $(-*)$-kernel, for which 
the corresponding reproducing $(-*)$-kernel Hilbert space $\Hc^K\subset \Ci(G_A/G_B,D)$
(respectively 
$\Hc^K\subset \Oc(G_A/G_B,D)$) consists of sections of the form 
$F_h:=[(\ \cdot\ ,P(\pi_A(\ \cdot\ )^{-1}h))]$, 
$h\in\overline{\spann}_{\mathbb C}(\pi_A(G_A)\Hc_B)$ 
in $\Hc_A$}.

To see this, let $s_j=u_jG_B\in G_A/G_B$ 
and $\xi_j=[(u_j^{-*},f_j)]\in D_{s_j^{-*}}$ for $j=1,\dots,n$. 
We have 
$$
\begin{aligned}
\sum_{j,l=1}^n(K(s_l,s_j^{-*})\xi_j\mid\xi_l)_{s_l,s_l^{-*}}
&=\sum_{j,l=1}^n(P(\pi_A(u_l^{-1})\pi_A(u_j^{-*})f_j)\mid f_l)_{\Hc_B} 
=\sum_{j,l=1}^n(\pi_A(u_l^{-1})\pi_A(u_j^{-*})f_j\mid f_l)_{\Hc_A} \\
&=\Bigl(\sum_{j=1}^n\pi_A(u_j^{-*})f_j\mid
 \sum_{l=1}^n\pi_A(u_l^{-*})f_l\Bigr)_{\Hc_A}\ge0.
\end{aligned}
$$
On the other hand, by the above calculation we get 
$$
\begin{aligned}
(K(s_l,s_j^{-*})\xi_j\mid\xi_l)_{s_l,s_l^{-*}}
&=(\pi_A(u_j^{-*})f_j\mid\pi_A(u_l^{-*})f_l)_{\Hc_A} 
=\overline{(\pi_A(u_l^{-*})f_l\mid\pi_A(u_j^{-*})f_j)}_{\Hc_A} \\
&=\overline{(K(s_j,s_l^{-*})\xi_l\mid\xi_j)}_{s_j,s_j^{-*}} 
=(\xi_j\mid K(s_j,s_l^{-*})\xi_l)_{s_j^{-*},s_j}. 
\end{aligned}
$$
Thus $K$ is a reproducing $(-*)$-kernel. 
Again by the above calculation it follows that 
\begin{equation}\label{2and3/4}
(K_{\xi_j}\mid K_{\xi_l})_{\Hc^K}
=(K(s_l,s_j^{-*})\xi_j\mid\xi_l)_{s_l,s_l^{-*}}
=(\pi_A(u_j^{-*})f_j\mid\pi_A(u_l^{-*})f_l)_{\Hc_A}.
\end{equation}
Now, by Proposition~\ref{cont}, $\Hc^K\subset \Cc(G_A/G_B,D)$. 
Let $F$ be a section in $\Hc^K$. 
By definition $F$ is a limit, in the norm of $\Hc^K$, of a sequence of sections of the form 
$\sum_{j=1}^{n(m)}K_{\xi_j^m}$, where $\xi_j^m=[(v_j^m,f_j^m)]\in D$, $j=1,\dots,n(m)$,
$m=1,2,\dots$ By \ref{2and3/4}, $\sum_{j=1}^{n(m)}\pi_A(v_j^m)f_j^m$ is a Cauchy sequence 
in $\Hc_A$, so that there exists $h:=\lim_{m\to\infty}\sum_{j=1}^{n(m)}\pi_A(v_j^m)f_j^m$ 
in $\Hc_A$. 
Now, by the proof of Proposition~\ref{cont}, convergence in $\Hc^K$ implies
(locally uniform) convergence in $\Cc(G_A/G_B,D)$ whence, for every $s=uG_B$ in $G_A/G_B$, we get 
$$  
\begin{aligned}
F(s)&=\lim_{m\to\infty}\sum_{j=1}^{n(m)}K_{\xi_j^m}(s)
=\lim_{m\to\infty}\sum_{j=1}^{n(m)}[(u,P(\pi_A(u)^{-1}\pi_A(v_j^m)f_j^m))] \\
    &=\lim_{m\to\infty}[(u,P(\pi_A(u)^{-1}\sum_{j=1}^{n(m)}\pi_A(v_j^m)f_j^m))]
\end{aligned}
$$
in $D_s$.
On the other hand, since the norm in $D_s$ is the copy of the norm in $\Hc_A$, 
through the action of
the basic mapping $\Phi$ associated with data $(\pi_A,\pi_B,P)$ 
(see Example~\ref{rel} and the bottom
of Definition~\ref{homog}), we also have 
$\lim_{m\to\infty}[(u,P(\pi_A(u)^{-1}\sum_{j=1}^{n(m)}\pi_A(v_j^m)f_j^m))]
=[(u,P(\pi_A(u)^{-1}h))]$.
Thus we have shown that $F=F_h$.
Also, for arbitrary $h\in\Hc_A$, 
$$
F_h=0\iff(\forall u\in G_A)\quad P(\pi_A(u^{-1})h)=0
\iff (\forall u\in G_A)\quad \pi_A(u^{-1})h\perp\Hc_B=0. 
$$
Since $\pi_A$ is a $*$-representation, 
it then follows that $F_h=0$ if and only if $h\perp\spann_{\mathbb C}(\pi_A(G_A)\Hc_B)$. Hence  
$\Hc_A/([\spann_{\mathbb C}(\pi_A(G_A)\Hc_B)]^\perp)=\Hc^K
=\{F_h\mid h\in\overline{\spann}_{\mathbb C}(\pi_A(G_A)\Hc_B)\}$.

Finally, note that $\Hc^K\subset\Ci(G_A/G_B,D)$ indeed, by definition of $F_h$ ($h\in\Hc_A$).
In the case where $G_A$ and $G_B$ are complex Banach-Lie groups then
$\Hc^K\subset\Oc(G_A/G_B,D)$, by the definition of $F_h$ as well.

Clearly, the mapping 
$F_h\mapsto h,\ \Hc^K\to\overline{\spann}_{\mathbb C}(\pi_A(G_A)\Hc_B)\}\subset\Hc_A$ is an isometry,
which we denote by $W$, such that $W(K_\eta)=\pi_A(v)f$ for $\eta=[(v,f)]\in D$. In addition, 
if $\overline{\spann}_{\mathbb C}\pi_A(\G_A)\Hc_B=\Hc_A$ then the operator $W$ is unitary. Recall
the mapping $\widehat{K}\colon D\to\Hc^K$ given by $\widehat{K}(\xi)=K_\xi$ if $\xi\in D$, as in 
\eqref{coh}. 
Clearly $W\circ\widehat{K}=\Theta$, where $\Theta$ is the basic mapping for the data 
$(\pi_A,\pi_B,P)$ 
(see Definition~\ref{homog}). 
\qed
\end{example}

The following result is an extension of Theorem~5.4 in \cite{BR07} and provides 
geometric realizations for $*$-representations of involutive Banach-Lie groups. 

\begin{theorem}\label{induction}
In the preceding setting, the following assertions hold: 
\begin{itemize}
\item[\rm(a)] The linear operator 
$$
\gamma\colon\Hc_A\to\Hc^K\subset\Ci(G_A/G_B,D),\quad 
(\gamma(h))(u G_B)=[(u,P(\pi_A(u^{-1})h))],$$
satisfies 
$\Ker\gamma=(\spann_{\mathbb C}(\pi_A(G_B)\Hc_B))^\perp$ 
and the operator 
$\iota:=\gamma\circ W$ is the canonical inclusion 
$\Hc^K\hookrightarrow\Ci(G_A/G_B,D)$. 
Moreover, $\gamma\circ\Theta=\widehat{K}$. 
\item[\rm(b)] For every point $t\in G_A/G_B$ the evaluation map
$\ev^\iota_t=\iota(\cdot)(t)\colon\Hc^K\to D_t$ 
is 
a continuous linear operator such that  
$$
(\forall s,t\in G_A/G_B)\quad 
K(s,t^{-*})=\ev^\iota_s\circ(\ev^\iota_t)^{-*}.
$$
\item[\rm(c)] The mapping $\gamma$ is a realization operator in the sense 
that it is an intertwiner 
between the $*$-representation $\pi_A\colon G_A\to\Bc(\Hc_A)$ 
and the natural representation of $G_A$ on the space of cross sections 
$\Ci(G_A/G_B,D)$. 
\end{itemize}
\end{theorem}

\begin{proof}

(a) 
This part is just a reformulation of what has been shown prior to the statement of the theorem.
The equality $\gamma\circ\Theta=\widehat{K}$ is obvious.

(b) 
Let $t\in G_A/G_B$ arbitrary and 
then pick $u\in G_A$ such that $t=uG_B$. 
In particular, once the element $u$ is chosen, we get a norm 
on the fiber $D_t$ (see Definition~\ref{homog}) 
and then for every $F=F_h\in\Hc^K$, where $h\in\overline{\spann}_{\mathbb C}\pi_A(G_A)\Hc_B$, 
we have 
$$
\begin{aligned}
\Vert\ev^\iota_t(F_h)\Vert_{D_t}
&=\Vert\iota F_h(t)\Vert_{D_t}=\Vert[(u,P(\pi_A(u)^{-1}h))]\Vert_{D_t} \\
&=\Vert\pi_A(u)P(\pi_A(u)^{-1}h)\Vert_{\Hc_A} 
\le\Vert\pi_A(u)\Vert\cdot\Vert\pi_A(u^{-1})\Vert\cdot\Vert h\Vert_{\Hc_A}
=C_u\Vert F_h\Vert_{\Hc^K},
\end{aligned}
$$
so that the evaluation map
$\ev^\iota_t\colon\Hc_A\to D_t$
is continuous. 
 
Let us keep $s=u G_B$ fixed for the moment. 
We first prove that 
\begin{equation}\label{comp_adj}
(\widehat{K}|_{D_s})^{-*}
=\ev^\iota_{s^{-*}}\colon\Hc^K\to D_{s^{-*}}.
\end{equation}
To this end we check that condition~\eqref{2} 
in Definition~\ref{rkhs2} 
is satisfied with 
$\theta_s=\ev^\iota_{s^{-*}}\colon\Hc^K\to D_{s^{-*}}$. 
In fact, 
let $\xi=[(u,f)]\in D_s$ arbitrary. 

Then for all $h\in\overline{\spann}_{\mathbb C}(\pi_A(G_A)\Hc_B)$ we have 
$$
\begin{aligned}
(\xi\mid\theta_s F_h)_{s,s^{-*}}
&=([(u,f)]\mid[(u^{-*},P(\pi_A((u^{-*})^{-1})h))])_{s,s^{-*}} \\
&=(\pi_A(u)f\mid\pi_A(u^{-*})P(\pi_A((u^{-*})^{-1})h))_{\Hc_A} \\
&=(f\mid P(\pi_A(u^*)h))_{\Hc_A}=(\pi_A(u)f\mid h)_{\Hc_A} \\
&=(W(\gamma\circ\Theta)(\xi)\mid W(\gamma(h)))_{\Hc_A}=((\gamma\Theta)(\xi)\mid F_h)_{\Hc^K} \\
&=(\widehat K(\xi)\mid F_h)_{\Hc^K}=(K_\xi\mid F_h)_{\Hc^K}. 
\end{aligned}
$$ 
Now let $s,t\in G_A/G_B$ arbitrary 
and $u,v\in G_A$ such that $s=uG_B$ and $t=vG_B$. 
It follows by~\eqref{comp_adj} that 
$(\ev^\iota_{t})^{-*}=\widehat{K}_{D_{t^{-*}}}$, hence 
for every 
$\eta=[(v^{-*},f)]\in D_{t^{-*}}$ we have 
$$
\ev^\iota_s\circ(\ev^\iota_t)^{-*}\eta
=\ev^\iota_s(\widehat{K}(\eta))
=(\iota(K_\eta))(s) 
=[(u,P(\pi_A(u^{-1})\pi_A(v^{-*})f))] 
=K(s,t^{-*})\eta.
$$
(c) Let $h\in\Hc_A$ and $v\in G_A$ arbitrary. 
Then at every point $t=uG_B\in G_A/G_B$ we have 
$$
\begin{aligned}
\bigl(\gamma(\pi_A(v)h)\bigr)(t)
&=[(u,P(\pi_A(u^{-1})\pi_A(v)h))]
=[(u,P(\pi_A((v^{-1}u)^{-1})h))]  \\
&=v\cdot[(v^{-1}u,P(\pi_A((v^{-1}u)^{-1})h))]
=v\cdot\bigl(\gamma(h)\bigr)(v^{-1}t)
\end{aligned}
$$
and the proof ends. 
\end{proof}

Part (c) of the above theorem tells us that it is possible to 
realize representations like 
$\pi_A\colon G_A\to\Bc(\Hc_A)$ as natural actions on spaces of analytic sections. 
We next take advantage of these
geometric realizations to point out some phenomena of holomorphic extension in bundle vectors 
and sections of them.
Firstly, we record some auxiliary facts in the form of a lemma. 

\begin{lemma}\label{invol4}
Let $G_A$ be an involutive Banach-Lie group and $G_B$ an involutive 
Banach-Lie subgroup of $G_A$, 
and denote by $\beta\colon\,(v,uG_B)\mapsto vuG_B$, 
$G_A\times G_A/G_B\to G_A/G_B$,  
the corresponding transitive action. 
Also denote $U_X=\{u\in G_X\mid u^{-*}=u\}$ for $X\in\{A,B\}$. 
Then the following assertions hold:
\begin{itemize}
\item[{\rm(a)}]
There exists a correctly defined involutive diffeomorphism 
$$z\mapsto z^{-*},\quad G_A/G_B\to G_A/G_B,$$
defined by $uG_B\mapsto u^{-*}G_B$. 
This diffeomorphism has the property 
$\beta(v^{-*},z^{-*})=\beta(v,z)^{-*}$ 
whenever $v\in G_A$ and $z\in G_A/G_B$. 
\item[{\rm(b)}] 
The group $U_X$ is a Banach-Lie subgroup of $G_X$ for $X\in\{A,B\}$ 
and $U_B$ is a Banach-Lie subgroup of~$U_A$.
\item[{\rm(c)}] 
If $G_B^{+}=G_A^{+}\cap G_B$, 
then the mapping 
$$\lambda\colon\,u U_B\mapsto uG_B,\quad
 U_A/U_B\to G_A/G_B, 
$$ 
is a diffeomorphism of $U_A/U_B$ onto the fixed-point submanifold  
of the involutive diffeomorphism of $G_A/G_B$ introduced above in 
assertion~{\rm(a)}.
\end{itemize}
\end{lemma}

\begin{proof}
Assertion~(a) follows since the mapping 
$u\mapsto u^{-*}$ is an automorphism of $G$ 
(Remark~\ref{invol1}). 
The proof of assertion~(b) is straightforward. 
 
As regards~(c), 
what we really have to prove is the equality 
$\lambda(U_A/U_B)=\{z\in G_A/G_B\mid z^{-*}=z\}$.  
The inclusion $\subseteq$ is obvious. 
Conversely, let $z\in G_A/G_B$ with $z^{-*}=z$. 
Pick $u\in G_A$ arbitrary such that $z=uG_B$. 
Since $z^{-*}=z$, it follows that 
$u^{-1}u^{-*}\in G_B$. 
On the other hand $u^{-1}u^{-*}\in G_A^{+}$, 
hence the hypothesis $G_B^{+}=G_A^{+}\cap G_B$
implies that $u^{-1}u^{-*}\in G_B^{+}$. 
That is, there exists $w\in G_B$ such that 
$u^{-1}u^{-*}=ww^*$. 
Hence $uw=u^{-*}(w^*)^{-1}$, 
so that $uw=(uw)^{-*}$. 
Consequently $uw\in U_A$, and in addition 
$z=uG_B=uwG_B=\lambda(uw U_B)$. 
\end{proof}

The next theorem gives a holomorphic extension of the Hermitian vector bundles 
and kernels introduced in \cite{BR07}.

\begin{theorem}\label{realization_hom}
For $X\in\{A,B\}$, let $G_X$ be a complex Banach-Lie group and $G_B$ a Banach-Lie 
subgroup of $G_A$. As above, 
set $U_X=\{u\in G_X\mid u^{-*}=u\}$. Let 
$\pi_X\colon X\to\Bc(\Hc_X)$ be a holomorphic $*$-representation such that 
$\pi_B(u)=\pi_A(u)|_{\Hc_B}$ 
for all $u\in G_B$. Denote by 
$\Pi\colon D\to G_A/G_B$ the like-Hermitian vector bundle, 
$K$ the reproducing $(-*)$-kernel,
and $W\colon\Hc^K\to\Hc_A$ the isometry and 
$\gamma\colon\Hc_A\to\Ci(G_A/G_B,D)$ the realization operator
associated with the data~$(\pi_A,\pi_B,P)$, 
where $P\colon\Hc_A\to\Hc_B$ is the orthogonal projection.  

Also denote by $\Pi^{U}\colon D^{U}\to U_A/U_B$ 
the like-Hermitian vector bundle, 
$K^{U}$ the reproducing $(-*)$-kernel,
and $W^{U}\colon\Hc^{K^{U}}\to\Hc_A$ the isometry and 
$\gamma^{U}\colon\Hc_A\to\Ci(U_A/U_B,D)$ the operators
associated with the data~$(\pi_A|_{U_A},\pi_B|_{U_B},P)$

Let assume in addition that $G_B^+=G_A^+\cap G_B$. Then the following assertions hold: 
\begin{itemize}
\item[\rm(a)] The inclusion $\iota:=\gamma\circ W\colon\Hc^K\to\Oc(G_A/G_B,D)$ 
is the realization operator associated with the reproducing $(-*)$-kernel~$K$. 
Moreover, $\gamma$ intertwines 
the $*$-representation $\pi_A\colon G_A\to\Bc(\Hc_A)$ 
and the natural representation of $G_A$ on the space of cross sections 
$\Oc(G_A/G_B,D)$.  
\item[\rm(b)] The like-Hermitian vector bundle 
$\Pi^{U}\colon D^{U}\to U_A/U_B$ 
is actually a Hermitian vector bundle. 
The mapping 
$\lambda\colon\,uU_B\mapsto uG_B$, $U_A/U_B\hookrightarrow G_A/G_B$, 
is a diffeomorphism of $U_A/U_B$ onto a  
submanifold of $G_A/G_B$ and we have 
\begin{equation}\label{xif}
\lambda(U_A/U_B)=\{z\in G_A/G_B\mid z^{-*}=z\}.
\end{equation}
In addition, there exists an $U_A$-equivariant real analytic embedding
$\Lambda\colon D^{U}\to D$ 
such that the diagrams 
$$
\begin{CD}
D^{U} @>{\Lambda}>> D \\
@V{\Pi^{U}}VV @VV{\Pi}V \\
U_A/U_B @>{\lambda}>> G_A/G_B
\end{CD}
\qquad\text{and}\qquad
\begin{CD}
D^{U} @>{\Lambda}>> D \\
@A{\gamma^{U}(h)}AA @AA{\gamma(h)}A \\
U_A/U_B @>{\lambda}>> G_A/G_B
\end{CD}
$$
for arbitrary $h\in\overline{\spann}_{\mathbb C}(\pi_A(G_A)\Hc_B)$
are commutative, the mapping $\Lambda$ is a fiberwise isomorphism, and 
$\Lambda(D^{U})=\Pi^{-1}(\lambda(U_A/U_B))$. 
\item[\rm(c)] 
The inclusion  
$\iota^{U}:=\gamma^{U}\circ W\colon\Hc^K\to\Co(U_A/U_B,D^{U})$ 
is the realization operator 
associated with the reproducing kernel $K^{U}$, 
where $\Co(U_A/U_B,D^{U})$ is the subspace of
$\Ci(U_A/U_B,D^{U})$ of real analytic sections. In addition, 
$\gamma^{U}$ is an intertwiner between 
the unitary representation $\pi_A\colon U_A\to\Bc(\Hc_A)$ 
and the natural representation of $U_A$ on the space of cross sections 
$\Co(U_A/U_B,D^{U})$.  
\end{itemize}
\end{theorem}

\begin{proof}
(a) This follows by Theorem~\ref{induction}  
applied to the data $(\pi_A,\pi_B, P)$. 
The fact that the range of the realization operator $\iota$ 
consists only of holomorphic sections follows 
either by Proposition~\ref{cont}(c) or directly by 
the definition of $\gamma$ (see Theorem~\ref{induction}(a)). 

(b) The fact that $\Pi^{U}$ is a Hermitian vector bundle 
(Remark~\ref{hermitian}) 
follows by \eqref{strong}. 
Moreover, the asserted properties of $\lambda$ 
follow by Lemma~\ref{invol4}(c) 
as we are assuming that $G_B^{+}=G_B\cap G_A^{+}$. 
As regards the $U_A$-equivariant embedding $\Lambda\colon D^{U}\to D$, 
it can be defined as the mapping 
that takes every equivalence class 
$[(u,f)]\in D^{U}$ 
into the equivalence class $[(u,f)]\in D$. 
Then $\Lambda$ clearly has the wished-for properties. 

(c) Use again Theorem~\ref{induction}  
for the data $(\pi_A|_{U_A},\pi_B|_{U_B},P)$. 
It is clear from definitions (see also \cite{BR07}) 
that $K^{U}$ is a reproducing {\it kernel} indeed.
The fact that the range of $\gamma^{U}$, or $\iota^{U}$, consists of 
real analytic sections follows by the definition of $\gamma^{U}$ 
(see Theorem~\ref{induction}(a) again). 
Alternatively, one can use assertions (a)~and~(b) above 
to see that for arbitrary $h\in\Hc_A$ 
the mapping 
$\Lambda\circ\gamma^{U}(h)\circ\lambda^{-1}\colon\lambda(U_A/U_B)\to D$ 
is real analytic since it is a section of $\Pi$ over 
$\lambda(U_A/U_B)$ which extends to the \textit{holomorphic} section 
$\gamma(h)\colon G_A/G_B\to D$. 
\end{proof}

\begin{remark}\label{that_is_complexif}
\normalfont
Theorem~\ref{realization_hom} (b) says that 
the image of $\Lambda$ is precisely 
the restriction of $\Pi$ to the fixed-point set of 
the involution on the base $G_A/G_B$, 
and this restriction is a Hermitian vector bundle.   
This remark along with the alternate proof of assertion~(c) 
show that there exists a close relationship 
between the setting of Theorem~\ref{realization_hom} 
and the circle of ideas related to complexifications of 
real analytic manifolds, 
and in particular complexifications of compact homogeneous spaces 
(see for instance \cite{On60}, \cite{IS66}, 
\cite{Sz04} and the references therein 
for the case of finite-dimensional manifolds). 
Specifically, the manifold $U_A/U_B$ can be identified with 
the fixed-point set of the antiholomorphic involution 
$z\mapsto z^{-*}$ of $G_A/G_B$. 
Thus we can view $G_A/G_B$ as a \textit{complexification} of $U_A/U_B$. 
By means of this identification, 
we can say that for arbitrary $h\in\Hc_A$ 
the real analytic section $\gamma^{U}(h)\colon U_A/U_B\to D^{\U}$ 
can be holomorphically extended to the section 
$\gamma(h)\colon G_A/G_B\to D$.
\qed
\end{remark}

The complex structure of $G_A/G_B$ (with self-conjugate space $U_A/U_B$) can be 
suitably displayed in certain cases where $G_X$ is the group of invertibles of a 
$C^*$-algebra $X=A$ or $B$.

\begin{theorem}\label{suppl1}
Assume that $\1\in B\subseteq A$ are two $C^*$-algebras such that 
there exists a conditional expectation $E\colon A\to B$ from $A$ onto $B$. 
Denote the groups of invertible elements in $A$~and~$B$ 
by $\G_A$~and~$\G_B$, respectively, 
and consider the quotient map $q\colon\,a\mapsto a\G_B$, $\G_A\to \G_A/\G_B$. 

Let $\pg:=(\Ker E)\cap\ug_A$, 
which is a real Banach space acted on by 
$\U_B$ by means of the adjoint action $(u,X)\mapsto uXu^{-1}$. 
Consider the corresponding quotient map 
$\kappa\colon\,(u,X)\mapsto[(u,X)]$, $\U_A\times\pg\to\U_A\times_{\U_B}\pg$, 
and define the mapping 
$\Psi^E_0\colon\,(u,X)\mapsto u\exp(\ie X)$, $\U_A\times\pg\to \G_A$. 
Then there is a unique $\U_A$-equivariant, real analytic diffeomorphism 
$\Psi^E\colon\U_A\times_{\U_B}\pg\to \G_A/\G_B$ such that the diagram 
$$
\begin{CD}
\U_A\times\pg @>{}\Psi^E_0>> \G_A \\
@V{\kappa}VV @VV{q}V \\
\U_A\times_{\U_B}\pg @>{\Psi^E}>> \G_A/\G_B
\end{CD}
$$
is commutative.
Thus the complex homogeneous space $\G_A/\G_B$ has the structure of 
a $\U_A$-equivariant real vector bundle over its real form $\U_A/\U_B$, 
the corresponding projection being given by the composition 
(depending on the conditional expectation~$E$)
$$
\G_A/\G_B\mathop{\longrightarrow}\limits^{(\Psi^E)^{-1}}
\U_A\times_{\U_B}\pg
\mathop{\longrightarrow}\limits^{\Xi}\U_A/\U_B, 
$$ 
where the typical fiber of the vector bundle 
$\Xi\circ(\Psi^E)^{-1}$ 
is the real Banach space $\pg=(\Ker E)\cap\ug_A$. 
\end{theorem}

\begin{proof}
The uniqueness of $\Psi^E$ follows since the mapping 
$\kappa$ is surjective. 
For the existence of $\Psi^E$, 
note that for all $u\in\U_A$, $v\in\U_B$, and $X\in\pg$ we have 
$$
\begin{aligned}
q(\Psi^E_0(uv,v^{-1}Xv))
 &=q(uv\cdot\exp(\ie v^{-1}Xv))
 =q(uv\cdot v^{-1}\exp(\ie X)v)
 =q(u\exp(\ie X)v) \\
 &=u\exp(\ie X)vG_B
 =u\exp(\ie X)G_B 
 =q(\Psi^E_0(u,X)). 
\end{aligned}
$$
This shows that the mapping 
\begin{equation}\label{Psi}
\Psi^E\colon\,[(u,X)]\mapsto u\exp(\ie X)\G_B, \quad 
\U_A\times_{\U_B}\pg\to \G_A/\G_B,
\end{equation}
is well defined, and it is clearly $\U_A$-equivariant. 
Moreover, since $\kappa$ is a submersion and 
$\Psi^E\circ\kappa$ ($=q\circ\Psi^E_0$) 
is a real analytic mapping, 
it follows by Corollary~8.4(i) in \cite{Up85} 
that $\Psi^E$ is real analytic. 

Now we prove that $\Psi^E$ is bijective. 
To this end we need the following fact: 
\begin{equation}\label{PR}
\text{for all }a\in \G_A\text{ there exist a unique }
(u,X,b)\in\U_A\times\pg\times \G_B^{+}
\text{ such that }
a=u\cdot\exp(\ie X)\cdot b
\end{equation}
(see Theorem~8 in \cite{PR94}). 
It follows by \eqref{Psi} and \eqref{PR}  that the mapping 
$\Psi^E\colon\U_A\times_{\U_B}\pg\to \G_A/\G_B$ is surjective. 
To see that it is also injective, assume that 
$u_1\exp(\ie X_1)\G_B=u_2\exp(\ie X_2)\G_B$, 
where $(u_j,X_j)\in\U_A\times\pg$ for $j=1,2$. 
Then there exists $b_1\in \G_B$ such that 
$u_1\exp(\ie X_1)b_1=u_2\exp(\ie X_2)$. 
Let $b_1=vb$ be the polar decomposition of $b_1\in \G_B$, 
where $v\in\U_B$ and $b\in \G_B^{+}$. 
Then 
$$
u_1\exp(\ie X_1)b_1=u_1\exp(\ie X_1)vb
=u_1v\exp(\ie v^{-1}X_1v)b.
$$
Note that $u_1v\in\U_A$ and $v^{-1}X_1v\in\pg$ 
since $E(v^{-1}X_1v)=v^{-1}E(X_1)v=0$. 
Since $u_1\exp(\ie X_1)b_1=u_2\exp(\ie X_2)$, 
it then follows by the uniqueness assertion in \eqref{PR} 
that $u_2=u_1v$ and $X_2=v^{-1}X_1v$. 
Hence $[(u_1,X_1)]=[(u_2,X_2)]$, 
and thus the mapping 
$\Psi^E\colon\U_A\times_{\U_B}\pg\to \G_A/\G_B$ is injective as well. 

Finally, we show that the inverse function 
$$
(\Psi^E)^{-1}\colon a\G_B=u\exp(\ie X)\G_B\mapsto[(u,X)],\ \G_A/\G_B\to\U_A\times_{\U_B}\pg
$$ 
is also smooth. 
For this,  note that $u$ and $X$ in \eqref{PR} depend on $a$ in a real analytic fashion 
(see \cite{PR94}). 
Hence, the mapping $\sigma\colon a\mapsto[(u,X)],\ \G_A\to\U_A\times_{\U_B}\pg$ is smooth. Since 
$\sigma=(\Psi^E)^{-1}\circ q$ and $q$ is a submersion, it follows again from Corollary~8.4(i) in \cite{Up85} 
that $(\Psi^E)^{-1}$ is smooth. 
In conclusion, $\Psi^E$ is a real analytic diffeomorphism (see page 268 in \cite{Be06}), 
as we wanted to show. 
\end{proof}

\begin{remark}\label{suppl2}
\normalfont
From the observation in the second part of the above statement, 
it follows that the mapping $\Xi\circ(\Psi^E)^{-1}$  
can be thought of as an infinite-dimensional version of Mostow fibration; 
see \cite{Mo55}, \cite{Mo05} and Section~3 in \cite{Bi04} 
for more details on the finite-dimensional setting. 
See also Theorem~1 in Section~3 of \cite{Las78} 
for a related property of complexifications of compact symmetric spaces. 

In fact the construction of the diffeomorphism $\Psi^E$ in Theorem~\ref{suppl1} relies on 
the representation \eqref{PR}, and so it depends on the decomposition of $A$ obtained in terms of 
the expectation $E$, see \cite{PR94}. It is interesting to see how $\Psi^E$ depends explicitly 
on $E$ at the level of tangent maps: We have 
$$
\begin{aligned}
T_{(\1,0)}\kappa 
 &\colon\,(Z,Y)\mapsto((\1-E)Z,Y),\quad \ug_A\times\pg\to 
 T_{[(\1,0)]}(\U_A\times_{\U_B}\pg)\simeq\pg\times\pg,\\
T_{(\1,0)}(\Psi^E_0) 
 &\colon\,(Z,Y)\mapsto Z+\ie Y,\quad\ug_A\times\pg\to A, \\
T_{\1}q
 &\colon\,Z\mapsto(\1-E)Z,\quad A\to\Ker E,
\end{aligned}
$$
hence $T_{[(\1,0)]}(\Psi^E)((\1-E)Z,Y)=(\1-E)(Z+\ie Y)=(\1-E)Z+\ie Y$ 
whenever $Z\in\ug_A$ and $Y\in\pg$. Thus 
$$
T_{[(\1,0)]}(\Psi^E)\colon\,(Y_1,Y_2)\mapsto Y_1+\ie Y_2,\quad
\pg\times\pg\to\Ker E,
$$
which is an isomorphism of real Banach spaces since 
$\Ker E=\pg\dotplus\ie \pg$. 
\qed
\end{remark}

\begin{corollary}\label{supplcor1}
Let $A$ and $B$ two $C^*$-algebras as in the preceding theorem. Then 
$\G_A/\G_B\simeq \U_A\times_{\U_B}\pg$ is a complexification of $\U_A/\U_B$ with respect to 
the anti-holomorphic involutive diffeomorphism 
$$
u\exp(\ie X)\G_B\mapsto u\exp(-\ie X)	\G_B,\ \G_A/\G_B\to \G_A/\G_B 
$$
where $u\in\U_A$, $X\in\pg$ (alternatively, $[(u,X)]\mapsto[(u,-X)]$).
\end{corollary}

\begin{proof}
First, note that $\G^+_B=\G_B\cap \G^+_A$. This is a direct consequence 
of the fact that the $C^*$-algebras are closed 
under taking square roots of positive elements. So Theorem~\ref{suppl1} applies to get 
$\U_A/\U_B$ as the set of fixed points of the mapping $a\G_B\mapsto a^{-*}\G_B$ on $\G_A/¼G_B$, 
where 
$a^{-*}:=(a^{-1})^*$ for $a\in A$, and $*$ is the involution in $A$. 
By \eqref{PR}, every element $a\G_B$ in $\G_A/\G_B$ is of the form 
$a\G_B=u\exp(\ie X)\G_B$ with 
$u\in\U_A$ and $X\in\pg$, and the correspondence $u\exp(\ie X)\G_B\mapsto [(u,X)]$ 
sis a bijection. 
But then $(u\exp(\ie X))^{-*}=u\exp(-\ie X)$ since $X^*=-X$, and the proof ends.
\end{proof}

To put Theorem~\ref{suppl1} and Corollary~\ref{supplcor1} in a proper perspective, 
we recall that for $X\in\{A,B\}$ 
the Banach-Lie group $\G_X$ is the universal complexification of $\U_X$ 
(see Example~VI.9 in \cite{Ne02}, and also \cite{GN03}). 
Besides this, we have seen in Theorem~\ref{realization_hom} that 
the homogeneous space $\G_A/\G_B$ 
is a complexification of $\U_A/\U_B$. Now Corollary~\ref{supplcor1} 
implements such a complexification
in the explicit terms of a sort of polar decomposition (if $X\in\pg$ then 
$\exp(\ie X)^*=\exp(-\ie(-X))=\exp(\ie X)$ whence 
$\exp(\ie X)=\exp(\ie X/2)\ \exp(\ie X/2)^*$ is positive). For the group case, see \cite{GN03}.

\begin{remark}\label{cutout}
\normalfont
It is to be noticed that there is an alternative way to express the involution mapping 
considered in this section as multiplication by positive elements. This representation 
was suggested by Axiom~4 for involutions of homogeneous reductive spaces 
as studied in the paper~\cite{MR92}.

Specifically, under the conditions assumed above the following condition is satisfied: 
\begin{equation}\label{MR}
(\forall a\in \G_A)(\exists a_{+}\in \G_A^{+},\, b_{+}\in \G_B^{+})\quad a^{-*}=a_{+}ab_{+}. 
\end{equation}

To see this first note that 
we can assume $\Vert a^*\Vert<\sqrt 2$.
Then, if $\Hc$ is a Hilbert space such that $A$ is canonically embedded in $\Bc(\Hc)$ and $x\in\Hc$,
$$
\begin{aligned}
\Vert a^*ax\Vert^2<2\Vert ax\Vert^2
\Leftrightarrow
(x-a^*ax\mid x-a^*ax\bigr)_{\Hc}<(x\mid x\bigr)_{\Hc}
\Leftrightarrow
\Vert (\1-a^*a)x\Vert^2<\Vert x\Vert^2.
\end{aligned}
$$
Thus $\Vert\1-a^*a\Vert<1$ and so $\Vert\1-E(a^*a)\Vert=\Vert E(\1)-E(a^*a)\Vert<1$, whence 
$b_{+}:=E(a^*a)\in\G_B^{+}$. Now it is clear that (5.6) holds with $a_{+}:=(a^{-*})b_{+}^{-1}a^{-1}\in\G_A^{+}$.

As a consequence of (5.6), we have that $a^{-*}\G_B=a_{+}\ a\G_B$ for every $a\in\G_B$. Let us see the 
correspondence of such an identity with the decomposition of $a^{-*}\G_B$ given in Theorem~\ref{suppl1}. 
Since $a=ue^{iX}b$ in (5.5), we have $a^*a=(be^{iX}u^{-1})(ue^{iX}b)=be^{2iX}b$ and so 
$E(a^*a)=bE(e^{2iX})b$. It follows that $a^{-*}=a_{+}ab_{+}$ where $b_{+}=bE(e^{2iX})b$ and
$a_{+}=ue^{-iX}b^{-2}E(e^{2iX})^{-1}b^{-2}e^{-iX}u^{-1}$.
\qed
\end{remark}

There is a natural identification between the vector bundle 
$\Xi\colon \U_A\times_{\U_B}\pg\to\U_A/\U_B$ and the tangent bundle $T(\U_A/\U_B)\to\U_A/\U_B$. 
In view of Theorem~\ref{suppl1}, we get an interesting interpretation of 
the homogeneous space $\G_A/\G_B$ as the tangent bundle of $\U_A/\U_B$. 

\begin{corollary}\label{supplcor2}
In the above notation, the vector bundle $\Xi\colon \U_A\times_{\U_B}\pg\to\U_A/\U_B$ 
is $\U_A$-equivariantly isomorphic to the tangent bundle $T(\U_A/\U_B)\to\U_A/\U_B$. 
Hence, the composition 
$$
\G_A/\G_B\mathop{\longrightarrow}\limits^{(\Psi^E)^{-1}}
\U_A\times_{\U_B}\pg
\mathop{\longrightarrow}\limits^{\simeq} T(\U_A/\U_B) 
$$
defines a $\U_A$-equivariant diffeomorphism between the complexification 
$\G_A/\G_B$ and the tangent bundle $T(\U_A/\U_B)$ of the homogeneous space $\U_A/\U_B$.
\end{corollary}

\begin{proof}
Let $\alpha\colon\,(u,v\U_B)\mapsto uv\U_B$, $\U_A\times\U_A/U_B\to\U_A/\U_B$. 
Then let $p_0=\1\U_B\in\U_A/\U_B$ 
and $\partial_2\alpha\colon\U_A\times T(\U_A/\U_B)\to T(\U_A/\U_B)$ 
the partial derivative of $\alpha$ with respect to 
the second variable. 
Since $T_{p_0}(\U_A/\U_B)\simeq\pg$, by restricting 
$\partial_2\alpha$ to $\U_A\times T_{p_0}(\U_A/U_B)$ 
we get a mapping 
$\alpha^E_0\colon\U_A\times\pg\to T(\U_A/\U_B)$. 
Then it is straightforward to show that there exists a unique 
$\U_A$-equivariant diffeomorphism 
$\alpha^E\colon\U_A\times_{\U_B}\pg\to T(\U_A/\U_B)$ 
such that $\alpha^E\circ\kappa=\alpha^E_0$. 

Now it follows by Theorem~\ref{suppl1} that the composition
$\G_A/\G_B\mathop{\longrightarrow}\limits^{(\Psi^E)^{-1}}
\U_A\times_{\U_B}\pg
\mathop{\longrightarrow}\limits^{\alpha^E} T(\U_A/\U_B)$
defines a $\U_A$-equivariant diffeomorphism between the complexification 
$\G_A/\G_B$ and the tangent bundle $T(\U_A/\U_B)$ of 
the homogeneous space $\U_A/\U_B$. 
\end{proof}

\begin{remark}\label{suppl3}
\normalfont
It is known that conditional expectations can be regarded as connection forms 
of principal bundles, 
see \cite{ACS95}, \cite{CG99}, and \cite{Ga06}. 
Thus Corollary~\ref{supplcor2} leads to numerous examples of real
analytic Banach manifolds whose tangent bundles have complex structures 
associated with certain connections. 
See for instance \cite{LS91}, \cite{Bi03}, and \cite{Sz04} 
for the case of finite-dimensional manifolds. 
\qed
\end{remark}

\section{Stinespring representations}\label{sect5}

In this section we are going to apply the preceding theory of reproducing $(-*)$-kernels, 
for homogeneous
like-Hermitian bundles, to explore the differential geometric background of 
completely positive maps. 
Thus we shall find geometric realizations of the Stinespring representations 
which will entail an unexpected
bearing on the Stinespring dilation theory. 
Specifically, it will follow that the classical constructions of 
extensions of representations and induced representations of $C^*$-algebras 
(see \cite{Di64} and \cite{Ri74}, respectively), which seemed to pass beyond 
the realm of geometric structures, actually have geometric interpretations 
in terms of reproducing kernels on vector bundles. 
See Remark~\ref{DR} below for some more details. 

\begin{notation}\label{matrices}
\normalfont
For every linear map $\Phi\colon X\to Y$ between two vector spaces and 
every integer $n\ge1$ we denote 
$\Phi_n=\Phi\otimes\id_{M_n({\mathbb C})}\colon M_n(X)\to M_n(Y)$, 
that is, $\Phi_n((x_{ij})_{1\le i,j\le n})=(\Phi(x_{ij}))_{1\le i,j\le n}$ 
for every matrix $(x_{ij})_{1\le i,j\le n}\in M_n(X)$. 
\qed
\end{notation}

\begin{definition}\label{cp}
\normalfont
Let $A_1$ and $A_2$ be two unital $C^*$-algebras and $\Phi\colon A_1\to A_2$ 
a linear map. 
We say that $\Phi$ is {\it completely positive} if 
for every integer $n\ge1$ the map $\Phi_n\colon M_n(A_1)\to M_n(A_2)$ 
is positive in the sense that it takes positive elements in the $C^*$-algebra 
$M_n(A_1)$ to positive ones in $M_n(A_2)$. 

If moreover $\Phi(\1)=\1$ then we say that $\Phi$ is unital and in this case 
we have $\Vert\Phi_n\Vert=1$ for every $n\ge1$ by the Russo-Dye theorem 
(see e.g., Corollary~2.9 in \cite{Pa02}). 
\qed
\end{definition}

\begin{definition}\label{stinespring}
\normalfont
Let $A$ be a unital $C^*$-algebra, $\Hc_0$ a complex Hilbert space and 
$\Phi\colon A\to\Bc(\Hc_0)$ a unital completely positive map. 
Define a nonnegative sesquilinear form on $A\otimes\Hc_0$ by the formula 
$$
\Bigl(\sum_{j=1}^nb_j\otimes\eta_j\mid\sum_{i=1}^n a_i\otimes\xi_i\Bigr)
=\sum_{i,j=1}^n(\Phi(a_i^*b_j)\eta_j\mid\xi_i)
$$
for all $a_1,\dots,a_n,b_1,\dots,b_n\in A$, 
$\xi_1,\dots,\xi_n,\eta_1,\dots,\eta_n\in\Hc_0$ 
and $n\ge1$. 
In particular 
\begin{equation}\label{st}
\Bigl(\sum_{j=1}^nb_j\otimes\eta_j\mid\sum_{j=1}^nb_j\otimes\eta_j\Bigr)
=(\Phi_n((b_i^*b_j)_{1\le i,j\le n})\eta\mid\eta),
\end{equation}
where 
$\eta=\begin{pmatrix}\eta_1 \\ \vdots \\ \eta_n \end{pmatrix}
\in M_{n,1}({\mathbb C})\otimes\Hc_0$. 
Consider the linear space $N=\{x\in A\otimes\Hc_0\mid (x\mid x)=0\}$ 
and denote by $\Kc_0$ the Hilbert space obtained as 
the completion of $(A\otimes\Hc_0)/N$ 
with respect to the scalar product defined by $(\cdot\mid\cdot)$ 
on this quotient space. 

On the other hand define a representation $\widetilde{\pi}$ of $A$ 
by linear maps 
on $A\otimes\Hc_0$ by 
$$(\forall a,b\in A)(\forall \eta\in\Hc_0)\qquad
\widetilde{\pi}(a)(b\otimes\eta)=ab\otimes\eta.$$
Then every linear map 
$\widetilde{\pi}(a)\colon A\otimes\Hc_0\to A\otimes\Hc_0$ 
induces a continuous map $(A\otimes\Hc_0)/N\to(A\otimes\Hc_0)/N$, 
whose extension by continuity will be denoted by $\pi_\Phi(a)\in\Bc(\Kc_0)$. 
We thus obtain a unital $*$-representation 
$\pi_\Phi\colon A\to\Bc(\Kc_0)$
which is called the {\it Stinespring representation associated with $\Phi$}. 

Additionally, denote by $V\colon\Hc_0\to\Kc_0$ 
the bounded linear map obtained as 
the composition 
$$V\colon\Hc_0 \to A\otimes\Hc_0\to(A\otimes\Hc_0)/N\hookrightarrow\Kc_0,$$
where the first map is defined by $A\ni h\mapsto\1\otimes h\in A\otimes\Hc_0$ 
and the second map is the natural quotient map. 
Then $V\colon\Hc_0\to\Kc_0$ is an isometry satisfying 
$\Phi(a)=V^*\pi(a)V$ for all $a\in A$.
\qed
\end{definition}

\begin{remark}\label{GNS}
\normalfont
The construction sketched in Definition~\ref{stinespring} 
essentially coincides with the proof of the Stinespring theorem on dilations 
of completely positive maps (\cite{St55}); 
see for instance Theorem~5.2.1 in \cite{ER00} 
or Theorem~4.1 in \cite{Pa02}. 
Minimal Stinespring representations are 
uniquely determined up to a unitary equivalence; 
see Proposition~4.2 in \cite{Pa02}. 

We also note that in the case when $\dim\Hc_0=1$, that is, 
$\Phi$ is a state of $A$, 
the Stinespring representation associated with $\Phi$ coincides with 
the Gelfand-Naimark-Segal (GNS) representation associated with $\Phi$. Thus in this case the isometry 
$V$ identifies with an element $h$ in $\Kc_0$ such that 
$\Phi(a)=(\pi(a)h\mid h)_{\Hc}$ for all $a\in A$.
\qed
\end{remark}

We now start the preparations necessary for obtaining 
the realization theorem for Stinespring representations 
(Theorem~\ref{stinerealiz}). 

\begin{lemma}\label{schwarz}
Let $\Phi\colon A\to B$ be a unital completely positive map 
between two $C^*$-algebras. 
Then for every $n\ge1$ and every $a\in M_n(A)$ we have 
$\Phi_n(a)^*\Phi_n(a)\le\Phi_n(a^*a)$. 
\end{lemma}

\begin{proof}
Note that $\Phi_n\colon M_n(A)\to M_n(B)$ is 
in turn a unital completely positive map, 
hence after replacing $A$ by $M_n(A)$, $B$ by $M_n(B)$, 
and $\Phi$ by $\Phi_n$, 
we may assume that $n=1$. 
In this case we may assume $B\subseteq\Bc(\Hc_0)$ 
for some complex Hilbert space $\Hc_0$ 
and then, using the notation in Definition~\ref{stinespring} we have  
$$
\Phi(a^*a)=V^*\pi_\Phi(a^*a)V=V^*\pi_\Phi(a)^*\id_{\Kc_0}\pi_\Phi(a)V
\ge V^*\pi_\Phi(a)^*VV^*\pi_\Phi(a)V=\Phi(a)^*\Phi(a),
$$
where the second equality follows since the Stinespring representation 
$\pi_\Phi\colon A\to\Bc(\Kc_0)$ is in particular a $*$-homomorphism. 
See for instance Corollary~5.2.2 in \cite{ER00} for more details. 
\end{proof}

For later use we now recall the theorem of Tomiyama on 
conditional expectations.

\begin{remark}\label{tomiyama}
\normalfont
Let $\1\in B\subseteq A$ be two $C^*$-algebras and such that 
there exists a {\it conditional expectation} $E\colon A\to B$, that is, 
$E$ is a linear map satisfying $E^2=E$, $\Vert E\Vert=1$ and $\Ran E=B$. 
Then for every $a\in A$ and $b_1,b_2\in B$ we have 
$E(a^*)=E(a)^*$, $0\le E(a)^*E(a)\le E(a^*a)$, and 
$E(b_1ab_2)=b_1E(a)b_2$. 
(See for instance \cite{To57} or \cite{Sak71}.) 
Additionally, $E$ is completely positive and $E(\1)=\1$, and this explains 
why $E$ has the Schwarz property stated in the previous Lemma~\ref{schwarz}. 
\qed
\end{remark}

\begin{lemma}\label{two_squares}
Assume that $\1\in B\subseteq A$ are $C^*$-algebras 
with a conditional expectation $E\colon A\to B$ 
and a unital completely positive map $\Phi\colon A\to\Bc(\Hc_0)$
satisfying $\Phi\circ E=\Phi$, where $\Hc_0$ is a complex Hilbert space.
Denote by $\pi_A\colon A\to\Bc(\Hc_A)$ and $\pi_B\colon B\to\Bc(\Hc_B)$ 
the Stinespring representations associated with 
the unital completely positive maps 
$\Phi$ and $\Phi|_B$, respectively. 
Then $\Hc_B\subseteq\Hc_A$, and for every
$h_0\in\Hc_0$ and 
$b\in B$ we have the commutative diagrams 
$$
\begin{CD}
A @>{\iota_{h_0}}>> \Hc_A @>{\pi_A(b)}>> \Hc_A \\
@V{E}VV      @V{P}VV              @VV{P}V \\
B @>{\iota_{h_0}}>> \Hc_B @>{\pi_B(b)}>> \Hc_B
\end{CD}
$$
where $P\colon\Hc_A\to\Hc_B$ is the orthogonal projection, and  
$\iota_{h_0}\colon A\to\Hc_A$ is the map induced by $a\mapsto a\otimes h_0$. 
\end{lemma}

\begin{proof}
We first check that the right-hand square is a commutative diagram. 
In fact, it is clear from the construction in Definition~\ref{stinespring} 
that $\Hc_B\subseteq\Hc_A$ and for every $b\in B$ we have 
$\pi_A(b^*)|_{\Hc_B}=\pi_B(b^*)$. 
In other words, 
if we denote by $I\colon\Hc_B\hookrightarrow\Hc_A$ the inclusion map, 
then $\pi_A(b^*)\circ I=I\circ\pi_B(b^*)$. 
Now note that $I^*=P$ and take the adjoints in the previous equation to get 
$P\circ \pi_A(b)=\pi_B(b)\circ P$. 

To check that the left-hand square is commutative, first note that 
$E\otimes\id_{\Hc_0}\colon A\otimes\Hc_0\to A\otimes\Hc_0$ 
is an idempotent mapping. 
To investigate the continuity of this map, 
let $x=\sum\limits_{i=1}^na_i\otimes\xi_i\in A\otimes\Hc_0$ and note that 
$\bigl((E\otimes\id_{\Hc_0})x\mid(E\otimes\id_{\Hc_0})x\bigr)
=\bigl(\Phi_n\bigl((E(a_i^*)E(a_j))_{1\le i,j\le n}\bigr)\xi\mid\xi\bigr),
$
where $\xi=\begin{pmatrix}\xi_1\\ \vdots \\ \xi_n\end{pmatrix}
\in M_{n,1}({\mathbb C})\otimes\Hc_0$.
On the other hand 
$(E(a_i^*)E(a_j))_{1\le i,j\le n}=E_n(a^*)E_n(a)\le E_n(a^*a)=
E_n\bigl((a_i^*a_j)_{1\le i,j\le n}\bigr)$, 
where 
$$
a=\begin{pmatrix}a_1&\dots&a_n\\ 0 &\dots & 0\\ 
\vdots & & \vdots \\ 0 &\dots & 0\end{pmatrix}\in M_n(A)
$$
and the above inequality follows by Lemma~\ref{schwarz}. 
Now, since $\Phi_n\colon M_n(A)\to M_n(B)$ is a positive map, we get 
$\Phi_n\bigl((E(a_i^*)E(a_j))_{1\le i,j\le n}\bigr)\le
\Phi_n\bigl(E_n\bigl((a_i^*a_j)_{1\le i,j\le n}\bigr)\bigr)$. 
Furthermore we have $\Phi_n\circ E_n=(\Phi\circ E)_n=\Phi_n$ by hypothesis, 
hence 
$\bigl((E\otimes\id_{\Hc_0})x\mid(E\otimes\id_{\Hc_0})x\bigr)
\le \bigl(\Phi_n\bigl((a_i^*a_j)_{1\le i,j\le n}\bigr)\xi\mid\xi\bigr)
=(x\mid x)$. 
Thus the linear map $E\otimes\id_{\Hc_0}\colon A\otimes\Hc_0\to A\otimes\Hc_0$ 
is continuous (actually contractive) with respect to the semi-scalar product 
$(\cdot\mid\cdot)$ and then it induces a bounded linear operator 
$\widetilde{E}\colon\Hc_A\to\Hc_A$. 
Moreover, since $E^2=E$ and $E(A)=B$, 
it follows that $\widetilde{E}^2=\widetilde{E}$ 
and $\widetilde{E}(\Hc_A)=\Hc_B$. 
On the other hand, it is obvious that for every $h_0\in\Hc$ we have 
$\widetilde{E}\circ\iota_{h_0}=\iota_{h_0}\circ E$. 
Hence it will be enough to prove that $\widetilde{E}=P$. 

To this end let $x=\sum\limits_{i=1}^n a_i\otimes\xi_i\in A\otimes\Hc_0$ 
and $y=\sum\limits_{j=1}^n b_j\otimes\eta_j\in B\otimes\Hc_0$ arbitrary. 
We have 
$$
\begin{aligned}
\bigl((E\otimes\id_{\Hc_0})x\mid y\bigr)
&=\bigl(\sum\limits_{i=1}^n E(a_i)\otimes\xi_i\mid 
     \sum\limits_{j=1}^n b_j\otimes\eta_j\bigr)    
=\sum_{i,j=1}^n\bigl(\Phi(b_j^*E(a_i))\mid \eta_j\bigr) 
=\sum_{i,j=1}^n\bigl(\Phi(E(b_j^*a_i))\mid \eta_j\bigr) \\ 
&=\sum_{i,j=1}^n\bigl(\Phi(b_j^*a_i)\mid \eta_j\bigr) 
=(x\mid y),
\end{aligned}
$$
where the third equality follows since $E(ba)=bE(a)$ for all $a\in A$ and $b\in B$, 
while the next-to-last equality follows by the hypothesis $\Phi\circ E=\Phi$. 
Since $y\in B\otimes\Hc_0$ is arbitrary, 
the above equality shows that $(E\otimes\id_{\Hc_0})x-x\perp B\otimes\Hc_0$. 
This implies that $\widetilde{E}(\tilde{x})-\tilde{x}\perp\Hc_B$, 
whence $\widetilde{E}(\tilde{x})=P(\tilde{x})$ for all $x\in A\otimes\Hc_0$, 
where $x\mapsto\tilde{x}$, $A\otimes\Hc_0\to\Hc_A$,  
is the canonical map obtained as the composition 
$A\otimes\Hc_0\to(A\otimes\Hc_0)/N\hookrightarrow\Hc_A$. 
(See Definition~\ref{stinespring}.) 
Since $\{\tilde{x}\mid x\in A\otimes\Hc_0\}$ is 
a dense linear subspace of $\Hc_A$, 
it follows that $\widetilde{E}=P$ throughout $\Hc_A$, and we are done.
\end{proof}

\begin{remark}\label{generation}
\normalfont
Under the assumptions of the previous lemma, we also obtain that 
$\Hc_A=\overline{\spann}\, \pi_A(\U_A)\Hc_B$: 
By standard arguments in $C^*$-algebras, 
we have that $A=\spann_{\mathbb C}\U_A$ or, equivalently, $A=\spann\, \U_A\cdot B$ 
since we have $\1\in B$. 
So $A\otimes\Hc_0=\spann\, \U_A\cdot(B\otimes\Hc_0)$ whence by quotienting and then
by passing to the completion we get $\Hc_A=\overline{\spann}\, \pi_A(\U_A)\Hc_B$.

Hence the mapping $\gamma$ is an isometry from $\Hc_A$ onto $\Hc^K$ and the inverse mapping 
$\gamma^{-1}$ coincides with $W$, see the remark prior to Theorem~\ref{induction}.
\qed
\end{remark}

\begin{remark}\label{DR}
\normalfont
In the setting of Lemma~\ref{two_squares}, if the restriction of $\Phi$ to $B$ 
happens to be a nondegenerate $*$-representation of $B$ on $\Hc_0$, then $\Hc_B=\Hc_0$ 
and $\pi_B=\Phi|_B$ by the uniqueness property of the minimal Stinespring dilation 
(see Remark~\ref{GNS}). 
In this special case our Lemma~\ref{two_squares} is related 
to the constructions of extensions of representations 
(see Proposition~2.10.2 in \cite{Di64}) and 
induced representations of $C^*$-algebras 
(see Lemma~1.7, Theorem~1.8, and Definition~1.9 in \cite{Ri74}). 
\qed
\end{remark}

In the following theorem we are using notation of Section~\ref{sect4}.

\begin{theorem}\label{stinerealiz}
Assume that $B\subseteq A$ are two unital $C^*$-algebras such that 
there exists a conditional expectation $E\colon A\to B$ from $A$ onto $B$, 
and let $\Phi\colon A\to\Bc(\Hc_0)$ be a unital completely positive map 
satisfying $\Phi\circ E=\Phi$, where $\Hc_0$ is a complex Hilbert space. 
Let $(\pi_A|_{\G_A},\pi_B|_{\G_B},P)$ be the Stinespring data associated with 
$E$ and $\Phi$. 
Set $\lambda\colon\,u\U_B\mapsto u\G_B$, $\U_A/\U_B\hookrightarrow\G_A/\G_B$.

Then the following assertions hold: 
\begin{itemize}
\item[\rm(a)] There exists a real analytic diffeomorphism
$
a\mapsto (u(a),X(a),b(a)),\ \G_A\to \U_A\times \pg \times \G_B^+
$ 
so that $a=u(a)\exp(\ie X(a))b(a)$ for all $a\in A$, which induces the polar decomposition 
in $\G_A/\G_B$, 
$$
a\G_B=u(a)\exp(\ie X(a))\G_B, \qquad{a\in G_A}.
$$
\item[\rm(b)] The mapping
$
-*:\ u\exp(\ie X)G_B\mapsto u\exp(-\ie X)G_B, \ \G_A/\G_B\to \G_A/\G_B 
$
is an anti-holomorphic involutive diffeomorphism of $\G_A/\G_B$ such that 
$$
\lambda(\U_A/\U_B)=\{s\in \G_A/\G_B\mid s=s^{-*}\}.
$$
\item[\rm(c)] The projection
$$
u\exp(\ie X)\G_B\mapsto u\U_B, \ \G_A/\G_B\to \U_A/\U_B 
$$
has the structure of a vector bundle isomorphic to the tangent bundle 
$\U_A\times_{\U_B}\pg\to\U_A/\U_B$ of the manifold $\U_A/\U_B$.  
The corresponding isomorphism is given by 
$u\exp(\ie X)\G_B\mapsto[(u,X)]$ for all $u\in\U_A$, $X\in\pg$.
\item[\rm(d)] 
Set $\Hc(E,\Phi):=\{\gamma(h)\mid h\in\Hc_A\}\subset\Oc(G_A/G_B,D)$ where  
$\gamma\colon\Hc_A\to\Oc(G_A/G_B,D)$ is 
the realization operator defined by  
$\gamma(h)(a\G_B)=[(a,P(\pi_A(a)^{-1}h))]$ for $a\in \G_A$ and $h\in\Hc_A$.
Put $\gamma^{\U}:=\gamma(\cdot)|_{\U_A/\U_B}\colon\Hc_A\to\Co(\U_A/\U_B,D^{\U})$ and 
$\Hc^{\U}(E,\Phi):=\{\gamma^{\U}(h)\mid h\in\Hc_A\}$. Denote by $\mu(a)$ the operator 
on the spaces $\Co(\U_A/\U_B,D^{\U})$ and $\Oc(G_A/G_B,D)$ defined by natural
multiplication by $a\in \G_A$. Then $\Hc(E,\Phi)$ and
$\Hc^{\U}(E,\Phi)$ are Hilbert spaces isometric with
$\Hc_A$. Moreover, for every $a\in \G_A$ the following diagram is commutative
$$
\begin{CD}
\Hc_A @>{\gamma^{\U}}>> \Hc^{\U}(E,\Phi) @>{\simeq}>> \Hc(E,\Phi) \\ 
@V{\pi(a)}VV @VV{\mu(a)}V @VV{\mu(a)}V \\
\Hc_A @>{\gamma^{\U}}>> \Hc^{\U}(E,\Phi) @>{\simeq}>> \Hc(E,\Phi),
\end{CD}
$$
that is, $\gamma\circ\pi(a)=\mu(a)\circ\gamma$. 
\item[\rm(e)]
There exists an isometry $V_{E,\Phi}\colon\Hc_0\to\Hc(E,\Phi)$ such that 
$$
\Phi(a)=V_{E,\Phi}^*(T_{\1}\mu)(a)V_{E,\Phi}\ , \qquad{a\in A},
$$
where $T_{\1}\mu$ is the tangent map of $\mu(\cdot)|_{\Hc(E,\Phi)}$ at $\1\in\G_A$. 
In fact, $T_{\1}\mu$ is a Banach algebra
representation of $A$ which extends $\mu$.
\end{itemize}
\end{theorem}

\begin{proof}
(a) Let $(\pi_A|_{\G_A},\pi_B|_{\G_B},P)$ be the Stinespring data introduced in  
Lemma~\ref{two_squares}, so that $\Hc_A=\overline{\spann}\ \pi_A(\G_A)\Hc_B$ according to 
Remark~\ref{generation}. We have that $\G_B^{+}=\G_B\cap \G_A^{+}$ as a direct consequence 
of the fact that the $C^*$-algebras are closed under taking square roots of positive elements. 
Then parts (a)-(d) of the theorem follows immediately by application of 
Theorem~\ref{suppl1}, Corollary~\ref{supplcor1}, Corollary~\ref{supplcor2}
and Theorem~\ref{realization_hom}.

As regards (e) note that for every $a\in A$ and $h\in\Hc_A$, 
$$
T_{\1}\mu(a)\gamma(h)=(d/dt)|_{t=0}\mu(e^{ta})\gamma(h)
=(d/dt)|_{t=0}e^{ta}[(e^{-ta}(\cdot)\ ,P(\pi_A(\cdot)^{-1}\pi_A(e^{ta})h))]
=\gamma(\pi_A(a)(h)).
$$

Since $\gamma$ is bijective (and isometric) we have that 
$T_{\1}\mu(a)=\gamma^{-1}\pi_A(a)\gamma$ for all $a\in A$, whence it is clear that
$T_{\1}\mu$ becomes a Banach algebra representation 
(and not only a Banach-Lie algebra representation).

Now take $V_{E,\Phi}:=\gamma\circ V$ where $V$ is the isometry $V\colon\Hc_0\to\Hc_A$ given in 
Definition~\ref{stinespring}. It is clear that $\gamma^*=\gamma^{-1}$ and then that 
$V_{E,\Phi}$ is the isometry we wanted to find.
\end{proof}

\begin{remark}\label{that_is_stine}
\normalfont
Theorem~\ref{stinerealiz} extends to the holomorphic setting, 
and for Stinespring representations, the geometric realization framework given in 
Theorem~5.4 
of \cite{BR07} for GNS representations. As part of such an extension we have found that the
real analytic sections obtained in \cite{BR07} are always restrictions of holomorphic sections 
of suitable (like-Hermitian) vector bundles on fairly natural complexifications.

Part (e) of the theorem provides us with a strong geometric view of the completely positive
mappings on $C^*$-algebras $A$: such a map is the compression of the \lq\lq natural action of $A$"
(in the sense that it is obtained by differentiating the non-ambiguous natural action of $\G_A$) on
a Hilbert space formed by holomorphic sections of a vector bundle of 
the formerly referred to type.
\qed
\end{remark}

\section{Further applications and examples}\label{sect6}

\medskip
1) {\bf Banach algebraic amenability} 

\begin{example}\label{amenable}
\normalfont
Let $\Ag$ be a Banach algebra. 
A {\it virtual diagonal} of $\Ag$ is by definition 
an element $M$ in the bidual 
$\Ag$-bimodule $(\Ag\hat{\otimes}\Ag)^{**}$ such that
$$
{\rm y}\cdot M=M\cdot {\rm y} \qquad\text{and}\qquad {\rm{m}}(M)\cdot {\rm y}={\rm y} \, 
\qquad({\rm y}\in\Ag)
$$
where $\rm{m}$ is the extension to $(\Ag\hat{\otimes}\Ag)^{**}$ of 
the multiplication map in
$\Ag$, ${\rm y}\otimes {\rm y}'\mapsto {\rm y}{\rm y}'$. 
The algebra $\Ag$ is called {\it amenable} 
when it possesses a
virtual diagonal as above. When $\Ag$ is a $C^*$-algebra, 
then $\Ag$ is amenable if and only 
it is nuclear. 
Analogously, a {\it dual} Banach algebra $\Mg$ is called 
{\it Connes-amenable} if $\Ag$ has a 
virtual diagonal which in addition is {\it normal}. 
Then a von Neumann algebra $\Ag$ is
Connes-amenable if and only it is injective. 
For all these concepts and results, see 
\cite{Ru02}. 

Let $\Ag$ be a $C^*$-algebra and let $\Bg$ be a von Neumann algebra. 
By 
$\Rep(\Ag,\Bg)$ we denote the set of bounded representations 
$\rho\colon\Ag\to\Bg$ such that $\overline{\rho(\Ag)\Bg_*}=\Bg_*$ 
where $\Bg_*$ is the 
(unique) predual of $\Bg$ (recall that $\Bg_*$ is 
a left Banach $\Bg$-module). In the case 
$\Bg=\Bc(\Hc)$, for a complex Hilbert space $\Hc$, the property that 
$\overline{\rho(\Ag)\Bg_*}=\Bg_*$ is equivalent to have 
$\overline{\rho(\Ag)\Hc}=\Hc$, that is, $\rho$ is nondegenerate. 
Let $\Rep_*(\Ag,\Bg)$ 
denote the subset of $*$-representations in $\Rep(\Ag,\Bg)$. 
For a von Neumann algebra $\Mg$, we
denote by
$\Rep^\omega(\Mg,\Bg)$ the subset of homomorphisms in 
$\Rep(\Mg,\Bg)$ which are ultraweakly
continuous, or {\it normal} for short. As above, the set of $*$-representations of 
$\Rep^\omega(\Mg,\Bg)$ is denoted by 
$\Rep_*^\omega(\Mg,\Bg)$.  

From now on, 
$\Ag$, $\Mg$ will denote a nuclear $C^*$-algebra and 
an injective von Neumann algebra
respectively. Fix $\rho\in\Rep(\Ag,\Bg)$. 
The existence of a virtual diagonal $M$ for 
$\Ag$ allows us to define an operator $E_\rho\colon\Bg\to\Bg$ by
$$
(E_\rho(T)x\mid x')
:=M\left({\rm y}\otimes {\rm y}'\mapsto(\rho({\rm y})T)\rho({\rm y}')x\mid x')\right)
\equiv\int_{\Ag\otimes\Ag}(\rho({\rm y})T)\rho({\rm y}')x\mid x')\ dM({\rm y},{\rm y}')
$$
where $x,x'$ belong to a Hilbert space $\Hc$ such that 
$\Bg\hookrightarrow\Bc(\Hc)$ 
canonically, and $T\in\Bg$. 
In the formula, the \lq\lq integral" corresponds to the 
Effros notation, see \cite{CG98}. 
The operator $E_\rho$ is a 
bounded projection such that
$$
E_\rho(\Bg)=\rho(\Ag)':=\{T\in\Bg\mid T\rho({\rm y})=\rho({\rm y})T,\ {\rm y}\in\Ag\}.
$$
In fact, it is readily seen that 
$\Vert E_\rho\Vert\le \Vert M\Vert\ \Vert\rho\Vert^2$, so that $E_\rho$ 
becomes a conditional expectation provided that 
$\Vert M\Vert=\Vert\rho\Vert=1$. For instance, 
if $\rho$ is a $*$-homomorphism then its norm is one, 
see page 7 in \cite{Pa02}. The existence
of (normal) virtual diagonals of norm one in (dual) Banach algebras 
is not a clear fact 
in general, but it is true, and not simple, 
that such (normal) virtual diagonals exist for 
(injective von Neumann) nuclear $C^*$-algebras, 
see page 188 in \cite{Ru02}.    

For $\rho\in\Rep(\Ag,\Bg)$ and $T\in\Bg$, let  
$T\rho T^{-1}\in\Rep(\Ag,\Bg)$ defined as 
$(T\rho T^{-1})({\rm y}):=T\rho({\rm y})T^{-1}$ (${\rm y}\in\Ag$). 
Put 
$$\Sg(\rho):=\{T\rho T^{-1}\mid T\in G_{\Bg}\}
\text{ and }
\Ug(\rho):=\{T\rho T^{-1}\mid T\in\U_{\Bg}\}.$$ 
The set $\Sg(\rho)$ is called 
the \emph{similarity orbit} of $\rho$, and $\Ug(\tau)$ is called 
the \emph{unitary orbit} of $\tau\in\Rep_*(\Ag,\Bg)$. 
It is known that
$\Rep(\Ag,\Bg)$, endowed with the norm topology, 
is the discrete union of orbits 
$\Sg(\rho)$. 
Moreover, each orbit $\Sg(\rho)$ is 
a homogeneous Banach manifold 
with a reductive structure induced by the connection form $E_\rho$. 
In the same way,
$\Rep_*(\Ag,\Bg)$ is the disjoint union of orbits $\Ug(\tau)$, 
and the restriction of
$E_\rho$ on $\ug_{\Bg}$ is a connection form which induces 
a homogeneous reductive 
structure on $\Ug(\tau)$ 
-- see \cite{ACS95}, \cite{CG98}, and \cite{Ga06}. 
We next compile some more information, 
about the similarity and unitary orbits, 
which is obtained on the basis of results of the present
section.

Let $\rho\in\Rep(\Ag,\Bg)$. 
As it was proved in \cite{Bu81} 
and independently in \cite{Ch81} (see also \cite{Ha83}), there exists 
$\tau\in\Rep_*(\Ag,\Bg)\cap\Sg(\rho)$, whence $\Sg(\rho)=\Sg(\tau)$. 
Hence, without loss of generality, $\rho$ can be
assumed to be a $*$-representation, so that $\Vert\rho\Vert=1$. 
Moreover, since we are assuming that 
$\Ag$ is nuclear, we can choose a virtual diagonal $M$ of $\Ag$ of norm one. 
Thus the operator 
$E_\rho$ is a conditional expectation. 
Set $A:=\Bg$, $B:=\rho(\Ag)'$. 
With this notation, 
$\Sg(\rho)=\G_A/\G_B$ and $\Ug(\rho)=\U_A/\U_B$ diffeomorphically.

For $X\in\pg_\rho:=\Ker E_\rho\cap\ug_A$, 
let $[X]$ denote the equivalence class of $X$ under 
the adjoint
action of $\U_B$ on $\pg_\rho$ considered in Theorem~\ref{suppl1}. Also, 
set $e^{\ie X}:=\exp(\ie X)$.
\end{example}

\begin{corollary}\label{simiunita}
Let $\Ag$ be a nuclear $C^*$-algebra and let $\Bg$ be a von Neumann algebra. 
The following assertions hold:
\begin{itemize}
\item[\rm(a)] 
Each connected component of $\Rep(\Ag,\Bg)$ is 
a similarity orbit $\Sg(\rho)$, for
some $\rho\in\Rep_*(\Ag,\Bg)$. 
Moreover, each orbit $\Sg(\rho)$ is the disjoint union
$$
\Sg(\rho)=\bigcup_{[X]\in\pg_\rho/\U_B}\Ug(e^{\ie X}\rho e^{-\ie X})
$$
where $\Ug(e^{\ie X}\rho e^{-\ie X})$ is connected, 
for all $[X]\in\pg_\rho/\U_B$.
\item[\rm(b)]
The similarity orbit $\Sg(\rho)$ is a complexification of 
the unitary orbit $\Ug(\rho)$ with respect to 
the involutive diffeomorphism 
$ue^{\ie X}\rho e^{-\ie X}u^{-1}\mapsto ue^{-\ie X}\rho e^{\ie X}u^{-1}$  
($u\in\U_\Bg$).
\item[\rm(c)] 
The mapping 
$ue^{\ie X}\rho e^{-\ie X}u^{-1}\mapsto u\rho u^{-1},\ \Sg(\rho)\to\Ug(\rho)$ 
is a continuous retraction 
which defines a vector bundle diffeomorphic to the tangent bundle 
$\U_A\times_{\U_B}\pg_\rho\to\Ug(\rho)$ of $\Ug(\rho)$.
\item[\rm(d)] 
Let $\Hc_0$ be a Hilbert space such that $\Bg\hookrightarrow\Bc(\Hc_0)$. 
For every
$\rho\in\Rep(\Ag,\Bg)$ there exists a Hilbert space $\Hc_0(\rho)$ isometric 
with $\Hc_0$, which is formed
by holomorphic sections of a like-Hermitian vector bundle with 
base $\Sg(\rho)$. 
Moreover, $\Bg$ acts
continuously by natural multiplication on $\Hc_0(\rho)$, 
and the representation
$R$ obtained by  transferring $\rho$ on $\Hc_0(\rho)$ coincides with 
multiplication by $\rho$;
that is, $R({\rm y})F=\rho({\rm y})\cdot F$ for all ${\rm y}\in\Ag$ and 
section $F\in\Hc_0(\rho)$.
\end{itemize}
\end{corollary}

\begin{proof} (a) As said before, 
every similarity orbit of $\Rep(\Ag,\Bg)$ is of the form 
$\Sg(\rho)$ for some $\rho\in\Rep_*(\Ag,\Bg)$. 
Since $A=\Bg$ is a von Neumann algebra, the set of unitaries
$\U_A=\U_\Bg$ is connected whence it follows (as in Remark~\ref{connected}) 
that the orbits $\Sg(\rho)$ and 
$\Ug(e^{\ie X}\rho e^{-\ie X})$ are connected for all $X\in\pg_\rho$. 
For $X,Y\in\pg_\rho$, we have 
$\Ug(e^{\ie X}\rho e^{-\ie X})=\Ug(e^{\ie Y}\rho e^{-\ie Y})$ 
if and only if there exists $u\in\U_A$ such
that $e^{\ie Y}\rho e^{-\ie Y}=ue^{\ie X}\rho e^{-\ie X}$, 
which means that $u\in\U_B$ and $Y=uXu^{-1}$
(see Theorem~\ref{suppl1}). 
Hence $[X]=[Y]$. 
Finally, by Theorem~\ref{suppl1} again we have 
$\Sg(\rho)=\bigcup_{[X]\in\pg_\rho/\U_B}\Ug(e^{\ie X}\rho e^{-\ie X})$.

(b) This is Theorem~\ref{stinerealiz} (b).

(c) This follows by Theorem~\ref{stinerealiz} (c).

(d) Given $\rho$ in $\Rep(\Ag,\Bg)$, 
there is $\tau=\tau(\rho)$ in $\Rep_*(\Ag,\Bg)$ such that 
$\Sg(\rho)=\Sg(\tau)$. Now we fix a virtual diagonal of $\Ag$ of norm one 
and then define 
the conditional expectation $E_\rho\equiv E_{\tau(\rho)}$ as prior to 
this corollary. So 
$E_\rho\colon\Bg\to\Bc(\Hc_0)$ is a completely positive mapping and 
one can apply 
Theorem~\ref{stinerealiz} (d).
As a result, one gets a Hilbert space $\Hc(\rho):=\Hc(E_\rho,E_\rho)$ 
of holomorphic sections of 
a like-Hermitian bundle on $\Sg(\rho)=\G_\Bg/\G_B$ (where $B=\rho(\Ag)'$), 
and an isometry
$V_\rho:=V_{E_\rho}\colon\Hc_0\to\Hc(\rho)$, satisfying 
$E_\rho(\rho({\rm y}))=V_\rho^*T_\1\mu(\rho({\rm y}))V_\rho$ for all ${\rm y}\in\Ag$, 
in the notations
of Theorem~\ref{stinerealiz}. Note that $V_\rho^*V_\rho=\1$ and 
therefore the correspondence 
${\rm y}\mapsto V_\rho\rho({\rm y})V_\rho^*$ defines a (bounded) representation of 
$\Hc(\rho)$. Now take 
$\Hc_0(\rho):=V(\Hc_0)$ and define $R({\rm y})$ as the restriction of 
$V_\rho\rho({\rm y})V_\rho^*$ on 
$\Hc_0(\rho)$ for every ${\rm y}\in\Ag$. 
Clearly, $R$ is the transferred representation of $\rho$ from 
$\Hc_0$ to $\Hc_0(\rho)$. 

Also, for every $F\in\Hc_0(\rho)$ there exists $h_0\in\Hc_0$ such that 
$F=V_\rho(\1\otimes h_0)$, that is, $F(aG_B)=[(a,P(a^{-1}\otimes h_0)]$ 
for all $a\in \G_\Bg$, 
where $P$ is as in Lemma~\ref{two_squares}. 
Then, for ${\rm y}\in\Ag$,
$$
\begin{aligned}
R({\rm y})F&=R({\rm y})V_\rho(\1\otimes h_0)=V_\rho\ \rho({\rm y})V_\rho^*V_\rho(\1\otimes h_0)
=V_\rho\ \rho({\rm y})(\1\otimes h_0) \\
&=V_\rho(\rho({\rm y})\otimes h_0)=T_\1\mu(\rho({\rm y}))\ V_\rho(\1\otimes h_0)      
=\rho({\rm y})\cdot F,
\end{aligned}
$$
as we wanted to show.
\end{proof}

\begin{remark}\label{final}
\normalfont
\begin{enumerate}
\item[{\rm(i)}] The first part of Corollary~\ref{simiunita} (a) was already well known 
(see for example \cite{ACS95}). 
In the decomposition of the second part, 
the orbit $\Ug(e^{iX}\rho e^{-iX})$ 
for $X=0$ corresponds to the unitary orbit of $\rho$. 
So the disjoint union supplies 
a sort of configuration of the similarity orbit $\Sg(\rho)$ by relation with 
the unitary orbit $\Ug(\rho)$.

\item[{\rm(ii)}]  Parts (a), (b), (c) of Corollary~\ref{simiunita} are consequences of 
the Porta-Recht decomposition given in \cite{PR94}, see \eqref{PR}. 
Such a decomposition has been 
considered previously in relation with similarity orbits of 
nuclear $C^*$-algebras, 
though in a different perspective, see Theorem 5.7 in \cite{ACS95}, 
for example.

\item[{\rm(iii)}] Corollary~\ref{simiunita} admits a version entirely analogous for 
injective 
von Neumann algebras $\Mg$ (replacing the nuclear $C^*$ algebra $\Ag$ of 
the statement) 
and representations in $\Rep^\omega(\Ag,\Bg)$ and $\Rep_*^\omega(\Ag,\Bg)$. 
Proofs are similar to the nuclear, $C^*$, case. 
For the analog of (d) one needs 
to take a normal virtual diagonal of $\Mg$ of norm one.

\item[{\rm(iv)}] Corollary~\ref{simiunita} applies in particular to locally compact groups $G$ 
for which the group 
$C^*$-algebra $C^*(G)$ is amenable, see \cite{ACS95} and \cite{CG99}. 
When the group is compact the method 
to define the expectation $E_\rho$ works for every representation $\rho$ 
taking values in any Banach algebra $A$. 
We shall see a particular example of this below, involving Cuntz algebras. 
\end{enumerate}
\qed
\end{remark}

\medskip
2) {\bf Completely positive mappings}

Let $A$ be a complex unital $C^*$-algebra, with unit $\1$, 
included in the algebra
$\Bc(\Hc)$ of bounded operators on a Hilbert space $\Hc$. 
Assume that $\Phi\colon A\to\Bc(\Hc)$
is a unital, completely bounded mapping. 
(In the following we shall assume freely that $\Hc$ is separable, 
when necessary.)

\begin{lemma}\label{cpositive} 
Given $\Phi$ as above and $u\in\G_A$, let $\Phi_u$ denote the mapping 
$\Phi_u:=u\Phi(u^{-1}\ \cdot\ u)u^{-1}$. Then 
\begin{itemize}
\item[{\rm(i)}] For every $u\in\G_A$, $\Phi_u$ is completely bounded and 
$\Vert\Phi_u\Vert_{cb}
\le\Vert\Phi\Vert_{cb}\Vert u^*\Vert\ \Vert u\Vert\ \Vert u^{-1}\Vert^2$.
\item[{\rm(ii)}] If $\Phi$ is completely positive then $\Phi_u$ is completely positive 
for every $u\in\U_A$. 
\end{itemize}
\end{lemma}

\begin{proof} (i) Let $n$ be a natural number. Take 
$f=(f_1,\cdots,f_n), h=(h_1,\cdots,h_n)\in\Hc^n$ and 
$(a_{ij})_{ij}\in M_n(A)$ all of them of
respective norms less than or equal to 1. 
In the following we shall think of $f$ and $g$ in their
column version. Then, for $u\in\G_A$, we have
$$
\begin{aligned}
\vert(\Phi^{(n)}_u(a_{ij})_{ij}f\mid h)_{\Hc^n}\vert &
=\vert\sum_{i,j}(\Phi(u^{-1}a_{ij}u)(u^{-1}f_j)\mid u^*h_i)_\Hc\vert \\
&\le\Vert\Phi^{(n)}(u^{-1}a_{ij}u)_{ij}\Vert_{\Bc(\Hc^n)}\ 
\Vert(u^{-1}f_j)_j\Vert_{\Hc^n}\ \Vert(u^*h_i)_i\Vert_{\Hc^n} \\
&\le\Vert\Phi\Vert_{cb}\ \Vert (u^{-1}I)(a_{ij})_{ij}(uI)\Vert_{\Bc(\Hc^n)}\ 
\Vert u^{-1}\Vert\ \Vert u^*\Vert
\le\Vert\Phi\Vert_{cb}\Vert u^{-1}\Vert^2\Vert u\Vert\ \Vert u^*\Vert.
\end{aligned}
$$ 
(ii) Assume now that $\Phi$ is completely positive. For natural $n$, 
take $(a_{ij})_{ij}\ge0$ in $M_n(A)$ and $h=(h_1,\cdots,h_n)\in\Hc^n$. Then
$$
\begin{aligned}
(\Phi^{(n)}_u(a_{ij})_{ij}h\mid h)_{\Hc^n}=&
((\sum_{j=1}^n\Phi_u(a_{1j})h_j,\dots,\sum_{j=1}^n\Phi_u(a_{nj})h_j)\mid h)_{\Hc^n} \\
&=\sum_{i,j=1}^n(\Phi_u(a_{ij})h_j\mid h_i)_\Hc
=\sum_{i,j=1}^n(\Phi^{(n)}(b_{ij}f\mid f)_{\Hc^n}
\end{aligned}
$$
where $f=u^{-1}h$, $b_{ij}=u^{-1}a_{ij}u$, $u\in\U_A$. So $(b_{ij})_{ij}\ge0$ in $M_n(A)$ 
and, since $\Phi$ is completely positive, 
we conclude that $\Phi^{(n)}(b_{ij})\ge0$ as we wanted
to show.
\end{proof}

Now let $\Phi\colon A\to\Bc(\Hc)$ be 
a fixed, unital \textit{completely positive} mapping. By
Proposition~3.5 in \cite{Pa02}, $\Phi$ is completely bounded. 
According to the preceding Lemma~\ref{cpositive}, if
$\Uc(\Phi)$ and
$\Sc(\Phi)$ denote, respectively, the {\it unitary} orbit
$\Uc(\Phi):=\{\Phi_u\mid u\in\U_A\}$ and the {\it similarity} orbit
$\Sc(\Phi):=\{\Phi_u\mid u\in\G_A\}$ of $\Phi$, 
there are natural actions of $\G_A$ on 
$\Sc(\Phi)$ and of $\U_A$ on $\Uc(\Phi)$, under usual conjugation.  
Note that the elements of the
orbit $\Sc(\Phi)$ are completely bounded maps but they 
do not need to be completely positive.

Put $\G(\Phi):=\{u\in\G_A\mid\Phi_u=\Phi\}$ and $\U(\Phi):=\G(\Phi)\cap\U_A$. 

\begin{corollary}\label{algebrly}
In the above notation, 
$\Sc(\Phi)=\G_A/\G(\Phi)$ and $\Uc(\Phi)=\U_A/\U(\Phi)$. 
\end{corollary}

\begin{proof} 
It is enough to observe that $\G(\Phi)$ and $\U(\Phi)$ are 
the isotropy subgroups 
of the actions of $\G_A$ on $\Sc(\Phi)$ and 
of $\U_A$ on $\Uc(\Phi)$, respectively.
\end{proof}

Note that $\G(\Phi)$ is defined by the family of polynomial equations
$$
\varphi(\Phi(axa^{-1})-a\Phi(x)a^{-1})=0,\qquad 
x\in A,\, \varphi\in\Bc(\Hc)_*,\, a\in\G_A
$$
on $\G_A\times\G_A$, so $\G(\Phi)$ is algebraic and a Banach-Lie group 
with respect to the relative topology of $A$
(see for instance the Harris-Kaup theorem in \cite{Up85}). 
To see when the isotropy groups $\G(\Phi)$ and $\U(\Phi)$ are  
Banach-Lie subgroups of $\G_A$, we need to  
compute their Lie algebras $\gg(\Phi)$ and $\ug(\Phi)$, respectively, 
and to see whether they are complemented subspaces of $A$. 

\begin{lemma}\label{differen} 
In the above notation we have 
$\gg(\Phi)=\{X\in A\mid(\forall a\in A)\quad\Phi([a,X])=[\Phi(a),X]\}$, and therefore 
$\ug(\Phi)=\{X\in\ug_A\mid(\forall a\in A)\quad\Phi([a,X])=[\Phi(a),X]\}$.
\end{lemma}

\begin{proof} 
To prove the inclusion ``$\subseteq$'' just note that 
if $X\in A$ and $e^{tX}:=\exp(tX)\in\G(\Phi)$, 
then for every $a\in A$ we get 
$\Phi(e^{tX}ae^{-tX})=e^{tX}\Phi(a)e^{-tX}$ for all $t\in\R$.  
Hence by differentiating in $t$ and
taking values at $t=0$ we obtain 
$\Phi(aX-Xa)=\Phi(a)X-X\Phi(a)$; that is, 
$\Phi([a,X])=[\Phi(a),X]$. 

Now let $X$ in the right-hand side of the first equality from the statement. 
Then $A$ is an invariant subspace for 
the mapping $\ad\,X=[X,\cdot]\colon\Bc(\Hc)\to\Bc(\Hc)$, 
since $X\in A$. 
In addition, $\Phi\circ(\ad\,X)|_A=(\ad\, X)\circ\Phi$. 
Hence for every $t\in{\mathbb R}$ and $n\ge0$ we get 
$\Phi\circ (t\ad\,X)^n|_A=(t\ad\, X)^n\circ\Phi$, 
whence $\Phi\circ\exp(t(\ad\,X)|_A)=\exp(t\ad\, X)\circ\Phi$. 
Since $\exp(t\ad\, X)b=\ee^{tX}b\ee^{-tX}$ for all $b\in\Bc(\Hc)$, 
it then follows that $\ee^{tX}\in\G(\Phi)$ for all $t\in{\mathbb R}$, 
whence $X\in\gg(\Phi)$. 
The remainder of the proof is now clear. 
\end{proof}

As regards the description of the isotropy Lie algebra $\gg(\Phi)$ 
in Lemma~\ref{differen}, let us note the following fact:  

\begin{proposition}\label{banalg} 
The isotropy Lie algebra $\gg(\Phi)$ is a closed involutive Lie subalgebra of $A$. 
If the range of $\Phi$ is contained in the commutant of $A$ then 
$\gg(\Phi)$ is actually a unital $C^*$-subalgebra of $A$, given by 
$\gg(\Phi)=\{X\in A\mid \Phi(aX)=\Phi(Xa)\text{ for all }a\in A\}$. 
In this case, $\G_{\gg(\Phi)}=\G(\Phi)$.
\end{proposition}

\begin{proof} 
It is clear from Lemma~\ref{differen} that $\gg(\Phi)$ is 
a closed linear subspace of $A$ which contains the unit $\1$. 
Moreover, since $\Phi(a^*)=\Phi(a)^*$ for all $a\in A$ 
(this is automatic by the Stinespring's dilation theorem, for instance), 
then for $X\in\gg(\Phi)$ and $a\in A$ we have
$
\Phi([X^*,a])=\Phi([a^*,X]^*)=\Phi([a^*,X])^*=[\Phi(a^*),X]^*=[X^*,\Phi(a)]
$
whence $\gg(\Phi)$ is stable under involution as well. The fact that 
$\gg(\Phi)$ is a Lie subalgebra of $A$ follows by Theorem 4.13 in \cite{Be06} 
(see the proof there).
 
If the range of $\Phi$ is contained in the commutant of $A$, 
then $\gg(\Phi)=\{X\in A\mid \Phi(aX)=\Phi(Xa)\text{ for all }a\in A\}$, 
and so $\gg(\Phi)$ is a $C^*$-subalgebra of $A$. Finally, note that
$u\in\G_{\gg(\Phi)}$ if and only if $u\in\G_A$ and $\Phi(uau^{-1})=\Phi(a)=u\Phi(a)u^{-1}$ 
(since $\Phi(A)\subseteq A'$), if and only if $u\in\G(\Phi)$. 	
\end{proof}

The condition in the above statement for $\Phi$ to be contained in the commutant of $A$ holds if 
for instance, $\Phi$ is a state of $A$. Next, we give another example suggested by 
Example~\ref{amenable}. For a $C^*$-algebra $\Ag$ and von Neumann algebra $A$ with predual $A_{*}$, 
let $\rho\colon\Ag\to A$ be a bounded $*$-homomorphism such that $\overline{\rho(\Ag)A_*}=A_{*}$. 
Denote $w^*$ the (generic) weak operator topology in a von Neumann algebra. 

\begin{corollary}\label{commut} 
Assume that $\Ag$ is a nuclear $C^*$-algebra or 
an injective von Neumann algebra (in the second case we assume in addition that $\rho$ is normal), 
and that $A=\overline{\rho(\Ag)}^{w^*}$. Let $\Phi=E_\rho\colon A\to A$ be a conditional expectation 
associated with $\rho$ as in Example~\ref{amenable}. Then $B:=\Phi(A)\subseteq A'$ and therefore 
$\gg(\Phi)$ is a von Neumann subalgebra of $A$. Also, $B$ is commutative and so 
it is isomorphic to an algebra of $L^\infty$ type.   
\end{corollary}

\begin{proof}
From $\Phi(A)=\rho(\Ag)'$ and $A=\overline{\rho(\Ag)}^{w^*}$ it is readily seen (recall that 
$A_{*}$ is an $A$-bimodule for the natural module operations) that $\Phi(a)$ 
commutes with every element of $A$ for all $a\in A$. The remainder of the corollary is clear.
\end{proof}

It is not difficult to find representations as those of the preceding corollary. If $\pi\colon\Ag\to\Bg$ is 
a representation as in Example~\ref{amenable}, then it is enough to take $A:=\overline{\pi(\Ag)}^{w^*}$ 
in $\Bg$, and 
$\rho\colon\Ag\to A$ defined by $\rho({\rm y}):=\pi({\rm y})$ (${\rm y}\in\Ag$), to obtain a representation satisfying
the hypotheses of Corollary~\ref{commut}. It is straightforward to check that $A$ is a $C^*$-algebra and, 
moreover, that $A$ is a dual Banach space. In effect, if $\Bg_{*}$ is the predual of $\Bg$ and $^\perp A$ is 
the pre-annihilator subspace of $A$ in $\Bg_{*}$, then the quotient $\Bg_{*}/^\perp A$ is a predual of $A$, and 
an $A$-submodule of $A^*$, such that $\overline{\rho(\Ag)(\Bg_{*}/^\perp A)}=\Bg_{*}/^\perp A$. Note that 
in the case when $\Bg=\Bc(\Hc)$ the von Neumann's bicommutant theorems says that $A=\pi(\Ag)''$.

\medskip
3) {\bf Conditional expectations}

As regards the isotropy group $\G(\Phi)$ of a completely positive map 
$\Phi\colon A\to\Bc(\Hc)$, 
we are going to see 
that it is actually a Banach-Lie subgroup of $\G_A$ 
(in the sense that the isotropy Lie algebra $\gg(\Phi)$ is 
a complemented subspace of $A$) 
in the important special case when $\Phi$ is a faithful normal conditional 
expectation. 
This will provide us with a wide class of completely
positive mappings whose similarity orbits 
illustrate the main results of the present paper.

Thus, let assume in this subsection that $\Phi=E$ is 
a faithful, normal, \textit{conditional expectation} 
$E\colon A\to B$, where $A$ and $B$ are von Neumann algebras with 
$B\subseteq A\subseteq\Bc(\Hc)$. 
In this case all of the elements in the unitary orbit $\Uc(E)$ are
conditional expectations, 
whereas all we can say about the elements in the similarity orbit $\Sc(E)$ is
that they are completely bounded quasi-expectations. 
We would like to present $\Sc(E)$ and $\Uc(E)$ as examples of
the theory given in 
the previous Theorems \ref{realization_hom} and/or \ref{suppl1}, 
or even Section~\ref{sect5}, of this paper.

Denote $A_E:=\{x\in B'\cap A\mid E(ax)=E(xa), a\in A\}$ 
and fix a faithful, normal state $\varphi$
on $B$. 
(Such a faithful state exists if the Hilbert space $\Hc$ is separable.) 
The set $A_E$ is a von Neumann subalgebra of $A$  
and, using the modular group of
$A$ induced by the gauge state $\psi:=\varphi\circ E$, 
it can be proven that there exists a faithful, normal, conditional
expectation $F\colon A\to A_E$ such that 
$E\circ F=F\circ E$ and $\psi\circ F=\psi$ (see
Proposition 4.5 in \cite{AS01}).
Set $\Delta=E+F-EF$. 
Then $\Delta$ is a bounded projection from $A$ onto 
$$
\Delta(A)=A_E+B=(A_E\cap\ker E)\oplus B.
$$
By considering  
the \emph{connected $\1$-component} $\G(E)^0=\G_{A_E}\cdot\G_B$ 
of the isotropy group $\G(E)$ 
(see Proposition~3.3 in \cite{AS01}), 
the existence of 
$\Delta$ implies that $\G(E)$ is in fact a Banach-Lie subgroup of $\G_A$, 
the orbits $\Sc(E)$ and 
$\Uc(E)$ are homogeneous Banach manifolds, and the quotient map 
$\G_A\to\Sc(E)\simeq\G_A/\G(E)$ is an analytic submersion, 
see Corollary~4.7 and Theorem~4.8
in \cite{AS01}. 
Also, the following assertions hold: 

\begin{proposition}\label{subalg} 
In the notations from above and from the first subsection,
$\Delta(A)=\gg(E)$. 
In particular, $A$ splits trough $\gg(E)$ 
and $\gg(E)$ is a $w^*$-closed Lie subalgebra of $A$. 
\end{proposition}

\begin{proof} 
By Theorem 4.8 in \cite{AS01} 
the quotient mapping $\G_A\to\Sc(E)=\G_A/\G(E)$ is an
analytic submersion. In fact the kernel of its differential is $\gg(E)$ 
(see Theorem 8.19 in 
\cite{Up85}).
Also, $\gg(E):=T_\1(\G(E))=\Delta(A)$ by Proposition 4.6 in \cite{AS01}. 
\end{proof}

Now let $\Phi\colon A\to\Bc(\Hc_0)$ be 
any unital completely positive map such that 
$\Phi\circ E=\Phi$ 
and apply Stinespring's dilation procedure to 
the mapping $\Phi$ and the conditional expectation $E\colon A\to B$. 
Thus, for $J=A,B$ there are the Hilbert spaces 
$\Hc_J(\Phi)$ 
and (Stinespring) representations 
$\pi_J\colon J\to\Bc(\Hc_J(\Phi))$ such that  
$\Hc_B(\Phi)\subseteq\Hc_A(\Phi)$ and  
$\pi_B(u)=\pi_A(u)|_{\Hc_B(\Phi)}$ for each $u\in B$, 
as given in Lemma~\ref{two_squares}. 
Denote by $P\colon\Hc_A(\Phi)\to\Hc_B(\Phi)$ 
the corresponding orthogonal projection. 

We are going to construct representations of 
the intermediate groups in the sequence 
$$\G_B\subseteq\G(E)^0\subseteq\G(E)\subseteq\G_A. $$
For this purpose set 
$\Hc_E(\Phi):=\overline{\spann}(\pi_A(\G(E))\Hc_B(\Phi))$ and 
$P_E$ the orthogonal projection from $\Hc_A(\Phi)$ onto
$\Hc_E(\Phi)$. 
We have that $\overline{\spann}(\pi_A(\G_A)\Hc_E(\Phi))=\Hc_A(\Phi)$, 
since
$\overline{\spann}(\pi_A(\G_A)\Hc_B(\Phi))=\Hc_A(\Phi)$ by 
Remark~\ref{generation}. 
For every $u\in\G(E)$, put 
$\pi_E(u):=\pi_A(u)|_{\Hc_E(\Phi)}$. 
Then $\pi_E(u)(\Hc_E(\Phi))\subseteq\Hc_E(\Phi)$ and so $(\pi_A,\pi_E,P_E)$
is a {\it data} in the sense of Definition~\ref{homog} 
(with holomorphic $\pi_A$ and $\pi_E$). 
Similarly, 
set $\Hc_E^0(\Phi):=\overline{\spann}(\pi_A(\G(E)^0)\Hc_B(\Phi))$ and 
$P_E^0$ the orthogonal projection from $\Hc_A(\Phi)$ onto
$\Hc_E^0(\Phi)$, 
and then for every $u\in\G(E)^0$, define
$\pi_E^0(u):=\pi_A(u)|_{\Hc_E^0(\Phi)}$.

Next set $D_B:=\G_A\times_{\G_B}\Hc_B(\Phi)$, 
$D_E^0:=\G_A\times_{\G(E)^0}\Hc_E(\Phi)$, and 
$D_E:=\G_A\times_{\G(E)}\Hc_E(\Phi)$. 
Let 
$\Hc_B(P,\Phi)$, $\Hc_E^0(P_E^0,\Phi)$ and $\Hc_E(P_E,\Phi)$
denote 
the (reproducing kernel) Hilbert spaces of holomorphic sections in 
these bundles, respectively, 
given by Theorems \ref{induction} and \ref{stinerealiz}(d). 

\begin{corollary}\label{sandwich} 
Let $B\subseteq A$ be unital von Neumann algebras, 
$E\colon A\to B$ be a faithful, normal, conditional expectation, 
and use the above notations. 
Then the inclusion maps 
$\Hc_B(\Phi)\hookrightarrow\Hc_E^0(\Phi)\hookrightarrow\Hc_E(\Phi)$ 
and $\G_B\hookrightarrow\G(E)^0\hookrightarrow\G(E)$ 
induce bundle homorphisms
$$
\begin{CD}
D_B @>>> D_E^0 @>{ }>> D_E \\ 
@VVV @VVV @VVV \\
\G_A/\G_B @>>> G_A/\G(E)^0 @>>> \Sc(E),
\end{CD}
$$
which leads to $\G_A$-equivariant isometric isomorphisms  
$\Hc_B(P,\Phi)\to\Hc_E^0(P_E^0,\Phi)\to\Hc_E(P_E,\Phi)$. 
In particular, the Stinespring 
representation $\pi_A|_{\G_A}\colon\G_A\to\Bc(\Hc_A(\Phi))$ can be realized 
as the natural 
representation $\mu\colon\G_A\to\Bc(\Hc_E(P_E,\Phi))$ 
on the vector bundle $D_E$ over the similarity orbit $\Sc(E)$. 
\end{corollary}

\begin{proof} 
Recall that $\Sc(E)\simeq\G_A/\G(E)$ and 
the elements or sections of 
the spaces 
$\Hc_B(P,\Phi)$, $\Hc_E^0(P_E^0,\Phi)$ and $\Hc_E(P_E,\Phi)$
are of the form 
$$
u\G_B\mapsto[(u,P(\pi_A(u^{-1})h))];\ 
u\G(E)^0\mapsto[(u,P_E^0(\pi_A(u^{-1})h))];\  
E_u\cong u\G(E)\mapsto[(u,P_E(\pi_A(u^{-1})h))],
$$
respectively, for $h$ running over $\Hc_A(\Phi)$. 
This gives us the quoted isometries. 
The fact that $\pi_A|_{\G_A}$ is
realized as $\mu$ acting on $\Hc_E(P_E,\Phi)$ is 
a consequence of Theorem~\ref{realization_hom}.
\end{proof}

Corollary~\ref{sandwich} admits a version in the unitary setting, that is, 
for the unitary groups $\U_A$, $\U_B$, $\U(E)^0$, $\U(E)$ and unitary orbit $\Uc(E)$ 
playing the role of 
the corresponding invertible groups and orbit. 
The following result answers in the affirmative the natural question of 
whether the similarity orbit $\Sc(E)\simeq\G_A/\G(E)$ 
endowed with the involutive diffeomorphism $a\G(E)\mapsto a^{-*}\G(E)$ 
is the complexification of the unitary orbit $\Uc(E)$ 
of the conditional expectation $E$. 

\begin{corollary}\label{similcomp} 
In the above situation, the similarity orbit $\Sc(E)$ of the conditional expectation $E$ is 
a complexification  of its unitary orbit $\Uc(E)$, and it is also 
$\U_A$-equivariantly diffeomorphic to the tangent bundle of $\Uc(E)$.
\end{corollary}

\begin{proof} 
Since the tangent bundles of $\Uc(E)$ 
and $\U_A/\U(E)$ coincide the assertion that 
the tangent bundle of $\Uc(E)$ is 
$\U_A$-equivariantly diffeomorphic to 
$\G_A/\G(E)$ is a consequence of 
Corollary~\ref{supplcor2}. 
On the other hand, as recalled above, 
to prove the fact that 
$\G_A/\G(E)$ is the complexification of $\U_A/\U(E)$ 
it will be enough to check that $\G(E)^+=\G_A^+\cap\G(E)$ 
(and then to apply Lemma~\ref{invol4}). 
The inclusion $\subseteq$ is obvious. 
Now let $c\in\G_A^+\cap\G(E)$. 
By Definition~\ref{invol}, there exists $g\in\G_A$ such that 
$c=g^*g\in\G(E)$. 
Then the reasoning from the proof of Theorem~3.5 in \cite{AS01} 
shows that $g^*g=ab$ with $a\in\G_{A_E}^+$ and $b\in\G_B^+$, 
whence $c=ab\in\G_{A_E}^+\cdot\G_B^+\subseteq\G(E)^+$. 
\end{proof}

\begin{remark}\label{weyl}
\normalfont
In connection with the commutative diagram 
of Corollary~\ref{sandwich}, note that   
since $\G(E)^0$ is the connected $\1$-component of $\G(E)$, 
it follows that the arrow $\G_A/\G(E)^0\to\G_A/\G(E)=\Sc(E)$ is actually 
a covering map whose fiber is the \emph{Weyl group} 
$\G(E)/\G(E)^0$ of the conditional expectation $E$ 
(cf. \cite{AS01} and the references therein). 
\qed
\end{remark}

\begin{remark}\label{staralgebra}
\normalfont
It is interesting to observe how Corollary~\ref{sandwich} looks in the case 
when $\gg(E)$ is an associative algebra, as in the second part of Proposition~\ref{banalg}.

Thus let assume that for a conditional expectation $E\colon A\to A$ as in former situations we have 
that $B:=E(A)\subseteq A'$. Then $B$ is commutative and $A\subseteq B'$ 
(note that $B'$ need not be commutative; in other words, $B$ is not maximal abelian). 
Hence, by Proposition~\ref{banalg}, 
$$ 
\gg(E)=\{X\in A\mid E(aX)=E(Xa), a\in A\}=\{X\in B'\cap A\mid E(aX)=E(Xa), a\in A\}:=A_E.
$$
By Proposition~\ref{subalg}, $A_E+B=\Delta(A)=\gg(E)=A_E$ whence $B\subseteq A_E$. 
Also, as regards to groups, Proposition~\ref{subalg} applies to give
$\G_{\gg(E)}=\G(E)$ whence we have 
$$
\G(E)=\G_{\gg(E)}=\G_{A_E}\subseteq\G_{A_E}\cdot\G_B=\G(E)^0\subseteq\G(E),
$$
and we obtain that $\G(E)^0=\G(E)$. 
This implies that the bundles 
$D^0_E\to\G_A/\G(E)^0$ and $D_E\to\Sc(E)$ of Corollary~\ref{sandwich} coincide.

Moreover, from the fact that $B\subseteq A_E=\gg(E)$ it follows that $F\circ E=E$ 
where $F$ is the conditional expectation given prior to Proposition~\ref{subalg}. 
In fact, for $a\in A$, $E(a)\in A_E=F(A)$ so there is some $a'\in A$ such that 
$E(a)=F(a')$. 
Then $(FE)(a)=F(F(a'))=(FF)(a')=F(a')=E(a)$ as required. 
Since 
$FE=EF$ we have eventually that $FE=EF=E$.

Suppose now that $\Phi\colon A\to\Bc(\Hc_0)$ is a completely positive mapping 
such that $\Phi\circ E=\Phi$. 
Then $\Phi\circ F=(\Phi\circ E)\circ F=\Phi\circ(EF)=\Phi\circ E=\Phi$, and one can use again the 
argument preceding Corollary~\ref{sandwich} to find vector bundles with corresponding Hilbert spaces 
(fibers) and representations $\pi_A\colon A\to\Bc(\Hc_A(\Phi))$, $\pi_{A_E}\colon A_E\to\Bc(\Hc_{A_E}(\Phi))$ 
and $\pi_B\colon B\to\Bc(\Hc_B(\Phi))$, and so on. 
In particular, from $E\vert_{A_E}\colon A_E\to B$ one gets 
$\Hc_{A_E}(\Phi)=\overline{\spann}\ \pi_{A_E}(\G(E))\Hc_B(\Phi)=\overline{\spann}\ \pi_A(\G(E))\Hc_B(\Phi)=\Hc_E(\Phi)$. 
Hence, in this case, the bundle $D_E\to\Sc(E)$ is a Stinespring bundle with respect to data 
${\pi_{A_E}}\vert_{\G(E)}$, ${\pi_B}\vert_{\G_B}$ 
(and the corresponding projection) to which 
Theorem~\ref{stinerealiz} can be applied. 

More precisely, part (c) of that theorem implies that $\Sc(E)$ is diffeomorphic to the tangent bundle 
$\U_A\times_{\U(E)}\pg^F$ of $\Uc(E)$, where $\pg^F=\ker F\cap\ug_A$, in the same way as $\G_A/\G_B$ is 
diffeomorphic to $\U_A\times_{\U_B}\pg^E$, $\pg^E=\ker E\cap\ug_A$.
\qed
\end{remark}

\medskip
4) {\bf Representations of Cuntz algebras}

We wish to illustrate the theorem on geometric realizations 
of Stinespring representations  
by an application to representations of Cuntz algebras. 
For the sake of simplicity we shall be working in the classical setting 
(\cite{Cu77}), although a part of what we are going to do 
can be extended to more general versions of these $C^*$-algebras 
(see \cite{CK80}, \cite{Pi97}, and also \cite{DPZ98}) 
or to more general $C^*$-dynamical systems.

\begin{example}\label{cuntz}
\normalfont
Let $N\in\{2,3,\dots\}\cup\{\infty\}$ 
and denote by $\Oc_N$ the $C^*$-algebra generated 
by a family of isometries $\{v_j\}_{0\le j<N}$ 
that act on the same Hilbert space and satisfy the condition 
that their ranges are mutually orthogonal,  
and in addition $v_0v_0^*+\cdots+v_{N-1}v_{N-1}^*=\1$ 
in the case $N\ne\infty$. 
The Cuntz algebra $\Oc_N$
has a canonical uniqueness property with respect to the choice of 
the generators $\{v_j\}_{1\le j<N}$ subject to the above conditions 
(see \cite{Cu77}). 
In particular, 
this implies that 
there exists a pointwise continuous \textit{gauge $*$-automorphism group} 
parameterized by the unit circle, 
$\lambda\mapsto\tau(\lambda)=\tau_\lambda$, 
${\mathbb T}\to\Aut(\Oc_N)$, 
such that $\tau_\lambda(v_j)=\lambda v_j$ if $0\le j<N$ 
and $\lambda\in{\mathbb T}$. 
For every $m\in{\mathbb Z}$ we denote  
\begin{equation}\label{note1}
\Oc_N^{(m)}=\{x\in\Oc_N\mid(\forall\lambda\in{\mathbb T})\quad 
\tau_\lambda(x)=\lambda^m x\}
\end{equation}
the spectral subspace associated with $m$, and then $\Fc_N:=\Oc_N^{(0)}$ (the fixed-point algebra 
of the gauge group). 
For every $m\in{\mathbb Z}$ we have
a contractive surjective linear idempotent mapping 
$E^{(m)}\colon\Oc_N\to\Oc_N^{(m)}$ defined by 
\begin{equation}\label{note2}
(\forall x\in\Oc_N)\qquad 
E^{(m)}(x)=\int\limits_{\mathbb T}\lambda^{-m}\tau_\lambda(x)\de\lambda,
\end{equation}
which is a faithful conditional expectation in the case $m=0$
(see for instance Theorem~V.4.3 in \cite{Dav96}).  
We shall denote $E^{(0)}=E$ for the sake of simplicity. 
\qed
\end{example}

The following statement is inspired by some remarks 
from Section~2 in \cite{Lac93}.  
It shows that, under specific hypothesis on the $C^*$-algebras 
from Theorem~\ref{stinerealiz},  
the corresponding reproducing kernel Hilbert space 
has a circular symmetry that resembles 
the one of the classical spaces of holomorphic functions on the unit disk. 
Thus we get series expansions and the natural setting of harmonic analysis 
in the spaces of bundle sections associated with completely positive maps. 

\begin{corollary}\label{fourier}
Let $N\in\{2,3,\dots\}\cup\{\infty\}$, a completely positive unital map 
$\Phi\colon\Oc_N\to\Bc(\Hc_0)$, 
and the corresponding Stinespring representation 
$\pi_\Phi\colon\Oc_N\to\Bc(\Kc_0)$, and isometry $V\colon\Hc_0\to\Kc_0$ such that 
$\Phi=V^*\pi V$. 
Put $\widetilde{\Hc}_0:=V(\Hc_0)$.  
Then the condition $\Phi\circ E=\Phi$ is satisfied 
if and only if $\Phi$ is gauge invariant, 
in the sense that for each $\lambda\in{\mathbb T}$ we have 
$\Phi\circ\tau_\lambda=\Phi$. 
In addition, if this is the case, then the following assertions hold: 
\begin{itemize}
\item[{\rm(a)}] Consider the geometric realization 
$\gamma\colon\Kc_0\to\Hc(E,\Phi)$ of 
the Stinespring representation $\pi_\Phi$  
and let $\Pi\colon D\to\G_{\Oc_N}/\G_{\Fc_N}$ 
be the corresponding homogeneous vector bundle. 
Then the gauge automorphism group of $\Oc_N$ induces 
smooth
actions $\widetilde{\tau}$ and $\overline{\tau}$ 
of the circle group ${\mathbb T}$ 
on the total space $D$ and the base $\G_{\Oc_N}/\G_{\Fc_N}$, 
respectively, 
such that the diagram 
$$
\begin{CD}
{\mathbb T}\times D @>{\widetilde{\tau}}>> D \\
@V{\id_{\mathbb T}\times\Pi}VV @VV{\Pi}V \\
\qquad{\mathbb T}\times(\G_{\Oc_N}/\G_{\Fc_N}) @>{\overline{\tau}}>> 
\G_{\Oc_N}/\G_{\Fc_N}
\end{CD}
$$
is commutative. 
The action on the base of the vector bundle also commutes with 
the natural involutive diffeomorphism thereof. 
\item[{\rm(b)}] If for all $m\in{\mathbb Z}$ we denote by $\Hc(E,\Phi)^{(m)}$ 
the closed linear subspace generated by $\gamma(\pi_\Phi(\Oc_N^{(m)})\widetilde{\Hc}_0)$, 
then we have the orthogonal direct sum decomposition 
$\Hc(E,\Phi)=\bigoplus\limits_{m\in\mathbb Z}\Hc(E,\Phi)^{(m)}$
and each term of this decomposition is $\G_{\Fc_N}$-invariant. 
\item[{\rm(c)}] For each $m\in{\mathbb Z}$, the orthogonal projection 
$P^{(m)}_{E,\Phi}\colon\Hc(E,\Phi)\to\Hc(E,\Phi)^{(m)}$ is given by the formula 
$$
(P^{(m)}_{E,\Phi}\Delta)(z)=\int\limits_{\mathbb T}\lambda^{-m}
(\widetilde{\tau}_\lambda\circ\Delta\circ
\overline{\tau}_\lambda^{-1})(z)\de\lambda
$$
whenever $z\in\G_{\Oc_N}/\G_{\Fc_N}$ and $\Delta\in\Hc(E,\Phi)$. 
\end{itemize}
\end{corollary}

\begin{proof}
Firstly note that \eqref{note2} implies that 
$E\circ\tau_\lambda=\tau_\lambda\circ E=E$ 
for all $\lambda\in{\mathbb T}$,  
because of the invariance property of the Haar measure $\de\lambda$ 
on the unit circle ${\mathbb T}$. 
Consequently, if we assume that $\Phi$ is gauge invariant, then 
for all $x\in\Oc_N$ we have 
$\Phi(E(x))
=\Phi\Bigl(\int\limits_{\mathbb T}\tau_\lambda(x)\de\lambda\Bigr)
=\int\limits_{\mathbb T}\Phi(\tau_\lambda(x))\de\lambda
=\int\limits_{\mathbb T}\Phi(x)\de\lambda
=\Phi(x)$.
Conversely, if $\Phi\circ E=\Phi$, 
then for every $\lambda\in{\mathbb T}$ we have 
$\Phi\circ\tau_\lambda=\Phi\circ E\circ\tau_\lambda 
=\Phi\circ E=\Phi$. 

(a) To define the action of ${\mathbb T}$ upon 
the base $\G_{\Oc_N}/\G_{\Fc_N}$ 
we use the fact that each gauge automorphism $\tau_\lambda$ 
leaves $\Fc_N$ pointwise invariant and therefore induces 
a mapping of $\G_{\Oc_N}/\G_{\Fc_N}$ onto itself. 
It is straightforward to show that in this way we get an action 
$(\lambda,a\G_{\Fc_N})\mapsto\tau_\lambda(a)\G_{\Fc_N}$ 
of ${\mathbb T}$ 
as claimed. 
The action of the circle group upon the total space $D$ 
can be defined by the formula 
$[(a,f)]\mapsto[(\tau_\lambda(a),f)]$ 
for all $[(a,f)]\in D$ and $\lambda\in{\mathbb T}$. 

(b) The realization operator $\gamma\colon\Kc_0\to\Hc(E,\Phi)$ is unitary, 
hence it will be enough to prove that 
$\Kc_0=\bigoplus\limits_{m\in\mathbb Z}\Kc_0^{(m)}$
and that each term of this decomposition is invariant 
under all of the operators in the $C^*$-algebra $\pi_\Phi(\Fc_N)$, 
where $\Kc_0^{(m)}$ is the closed linear subspace of $\Kc_0$ 
spanned by $\pi_\Phi(\Oc_N^{(m)})\widetilde{\Hc}_0$ for all $m\in{\mathbb Z}$. 
(Note that 
$\Kc_0=\overline{\spann}\ \pi_\Phi(\Fc_N)\widetilde{\Hc}_0$ by construction.)

The proof of this assertion 
follows the lines of Section~1 in \cite{Lac93} and 
relies on the fact that, 
as an easy consequence of \eqref{note1}, we have 
$\Oc_N^{(m)}\Oc_N^{(n)}\subseteq\Oc_N^{(m+n)}$ and 
$(\Oc_N^{(m)})^*\subseteq\Oc_N^{(-m)}$ for all $m,n\in{\mathbb Z}$
(where $(\cdot)^*$ stands for the image under the involution of~$\Oc_N$). 
Note that $VV^*\colon\Kc_0\to\widetilde{\Hc}_0$ is the orthogonal projection from 
$\Kc_0$ onto $\widetilde{\Hc}_0$. 
It follows that for all $m,n\in{\mathbb Z}$ with $m\ne n$, 
and $x\in\Oc_N^{(m)}$, $y\in\Oc_N^{(n)}$, $\xi,\eta\in\Hc_0$ we have 
$$
\begin{aligned}
(\pi_\Phi(x)V\xi\mid\pi_\Phi(y)V\eta)
&=(\pi_\Phi(y^*x)V\xi\mid V\eta)
=(VV^*(\pi_\Phi(y^*x)V\xi)\mid V\eta)\\
&=(\Phi(y^*x)\xi\mid\eta))_{\Hc_0}
=(\Phi(E(y^*x))\xi\mid\eta))_{\Hc_0}=0,
\end{aligned}
$$ 
where the latter equality follows since $y^*x\in\Oc_N^{(m-n)}$ 
with $m-n\ne0$, so that $E(y^*x)=0$ as an easy consequence 
of~\eqref{note2}. 
The above computation shows that $\Kc_0^{(m)}\perp\Kc_0^{(n)}$ 
whenever $m\ne n$. 

To see that $\bigcup\limits_{m\in{\mathbb Z}}\Kc_0^{(m)}$ 
spans the whole $\Kc_0$, 
just recall from \cite{Cu77} that the set 
$\bigcup\limits_{m\in{\mathbb Z}}\Oc_N^{(m)}$ 
spans a dense linear subspace of $\Oc_N$, 
and use the image of this set 
under the unital $*$-homomorphism 
$\pi_\Phi\colon\Oc_N\to\Bc(\Kc_0)$. 
Consequently the asserted orthogonal direct sum decomposition 
of $\Kc_0$ is proved. 
On the other hand, since $\Fc_N=\Oc_N^{(0)}$, 
it follows that for all $m\in{\mathbb Z}$ we have 
$\Fc_N\Oc_N^{(m)}=\Oc_N^{(0)}\Oc_N^{(m)}\subseteq\Oc_N^{(m)}$, 
so $\pi_\Phi(\Fc_N)\Kc_0^{(m)}\subseteq\Kc_0^{(m)}$ 
according to the construction of $\Kc_0^{(m)}$. 

(c) For every $\lambda\in{\mathbb T}$ denote by 
$a\otimes f\mapsto\overline{a\otimes f}$, $\Oc_N\otimes\Hc_0\to\Kc_0$, 
the canonical map. 
Since $\Phi\circ\tau_\lambda=\Phi$, it follows that the mapping 
$a\otimes f\mapsto\tau_\lambda(a)\otimes f$, $\Oc_N\otimes\Hc_0\to\Oc_N\otimes\Hc_0$, 
induces a unitary operator 
$V_\lambda\colon\overline{a\otimes f}\mapsto\overline{\tau_\lambda(a)\otimes f},\ \Kc_0\to\Kc_0$. 
As the closure of the image of $\Oc_N^{(m)}\otimes\Hc_0$ in $\Kc_0$ 
is equal to $\Kc_0^{(m)}$, it follows that for all $m\in{\mathbb Z}$ 
we have 
\begin{equation}\label{eigenspaces}
\Kc_0^{(m)}=\{h\in\Kc_0\mid(\forall\lambda\in{\mathbb T})
\quad V_\lambda h=\lambda^m h\}. 
\end{equation}
On the other hand, for all $c,a\in\Oc_N$ and $f\in\Hc_0$ we have 
$$
\begin{aligned}
\pi_\Phi(\tau_{\lambda}(c))\overline{a\otimes f}
&=\overline{\tau_{\lambda}(c)a\otimes f}
=\overline{\tau_{\lambda}(c\,\tau_{\lambda^{-1}}(a))\otimes f}
=V_{\lambda}\overline{c\,\tau_{\lambda^{-1}}(a)\otimes f}
=V_{\lambda}\pi_\Phi(c)\overline{\tau_{\lambda^{-1}}(a)\otimes f} \\
&=V_{\lambda}\pi_\Phi(c)V_{\lambda^{-1}}\overline{a\otimes f},  
\end{aligned}
$$
whence 
\begin{equation}\label{inter}
(\forall \lambda\in{\mathbb T},\,c\in\Oc_N)\qquad 
\pi_\Phi(\tau_{\lambda}(c))=V_\lambda\pi_\Phi(c)V_{\lambda^{-1}}. 
\end{equation}
It follows by \eqref{eigenspaces} and \eqref{inter} that  
\begin{equation}\label{weight}
(\forall \lambda\in{\mathbb T},\,n\in{\mathbb Z},\,
c\in\Oc_N,\,h\in\Kc_0^{(n)})\qquad 
\pi_\Phi(\tau_{\lambda}(c))h
=\lambda^{-n}V_{\lambda}\pi_\Phi(c)h.
\end{equation}
Now for $\lambda\in{\mathbb T}$, $\eta\in\Kc_0$, 
and $c\in\G_{\Oc_N}$ 
we get 
$$
\begin{aligned}
(\widetilde{\tau}_\lambda\circ
 \gamma(\eta)\circ\overline{\tau}_\lambda^{-1})(c\,\G_{\Fc_N})
 &=\widetilde{\tau}_\lambda(\gamma(\eta)(\tau_{\lambda^{-1}}(c)\G_{\Fc_N})) 
  =\widetilde{\tau}
 [(\tau_{\lambda^{-1}}(c),P(\pi_\Phi(\tau_{\lambda^{-1}}(c^{-1})\eta)))]\\
 &=[(c,P(\pi_\Phi(\tau_{\lambda^{-1}}(c^{-1})\eta)))], 
 \end{aligned}
$$
whence by \eqref{weight} it follows that for $h\in\Kc_0^{(n)}$ we have 
$(\widetilde{\tau}_\lambda\circ\gamma(h)\circ
\overline{\tau}_\lambda^{-1})(c\,\G_{\Fc_N}) 
=\lambda^n[(c,P(V_{\lambda^{-1}}\pi_\Phi(c^{-1})h))]$
where for $m\in{\mathbb Z}$, $P^{(m)}\colon\Kc_0\to\Kc_0^{(m)}$ 
is the orthogonal projection and $P=P^{(0)}$.

On the other hand, 
it follows by \eqref{eigenspaces} 
that 
$\int\limits_{\mathbb T}\lambda^{n-m}(V_\lambda h)\de\lambda
=P^{(m-n)}h$ for all $h\in\Kc_0$,  
and thus by the above equality we get for $h\in\Kc_0^{(n)}$, 
$$
\begin{aligned}
\int\limits_{\mathbb T}
\lambda^{-m}(\widetilde{\tau}_\lambda\circ\gamma(h)
\circ\overline{\tau}_\lambda^{-1})(c\,\G_{\Fc_N}) \de\lambda
&=\int\limits_{\mathbb T}\lambda^{n-m}
[(c,P(V_{\lambda^{-1}}\pi_\Phi(c^{-1})h))]\de\lambda \\
&=[(c,P^{(0)}P^{(n-m)}(\pi_\Phi(c^{-1})h))]  
=
 \begin{cases}\gamma(h)(c\,\G_{\Fc_N}) &\text{ if }n=m, \\
              0                      &\text{ if } n\ne m.                      
 \end{cases}
\end{aligned}
$$
Thus we have proved the asserted formula for $\Delta=\gamma(h)$ 
with $h\in\bigcup\limits_{n\in{\mathbb Z}}\Kc_0^{(n)}$. 
On the other hand, the reproducing $(-*)$-kernel associated with 
the Stinespring representation $\pi_\Phi$ is continuous, hence 
the realization operator 
$\gamma\colon\Kc_0\to\Oc(\G_{\Oc_N}/\G_{\Fc_N},D)$ 
is continuous with respect to the uniform convergence on 
the compact subsets of the base $\G_{\Oc_N}/\G_{\Fc_N}$. 
Since $\Kc_0=\bigoplus\limits_{n\in{\mathbb Z}}\Kc_0^{(n)}$ 
and the right hand side of the asserted formula 
is linear and continuous with respect to the latter topology, 
it follows that the corresponding equality extends by 
linearity and continuity to the whole space $\Hc(E,\Phi)=\gamma(\Kc_0)$. 
\end{proof}

\begin{remark}\label{wavelets}
\normalfont
It is noteworthy that orthogonal decompositions 
similar to the one of Corollary~\ref{fourier}(b) 
also show up in connection with representations of Cuntz algebras 
that do not necessarily occur as Stinespring dilations 
of gauge invariant maps;  
see for instance the representations studied in \cite{BJ99}. 
\qed 
\end{remark}

There is a close relationship between $*$-endomorphisms 
of algebras $\Bc(\Hc)$ and representations of Cuntz algebras, see for instance 
\cite{Lac93a}, \cite{Lac93}. 
In the remainder of this subsection, we point out that some of the notions underlying 
this relationship provide us with more examples of the theory proposed in the present paper.
Thus let $\Hc$ be a separable Hilbert space, let $A:=\Bc(\Hc)$ and let 
$\alpha\colon A\to A$ a unital $*$-endomorphism; 
then $\alpha$ is normal, as noted for instance in \cite{Lac93a}. 
By a celebrated result of W.~Arveson 
there exist
$N\in\{1,\dots,+\infty\}$ and a $*$-representation $\rho\colon\Oc_N\to A$, where $\Oc_N$ is the Cuntz algebra
generated by $N$ isometries $v_j$, on $\Hc$, such that 
\begin{equation}\label{arveson}
\alpha(T)=\sum_{j=1}^N \rho(v_j)T\rho(v_j^*),\ \qquad {T\in A}.
\end{equation}
(See Proposition~2.1 in \cite{Ar89}.) 
Since $\alpha$ is unital we have that
$\sum_{j=1}^N
\rho(v_j)\rho(v_j^*)=\1$ even for
$N=\infty$  
(in the strong operator topology, in the latter case). 
In the sequel we assume that $N\ge2$. 

We do observe that, in
order to make a link between the geometry of $\rho$ and the one of $\alpha$, the endomorphism $\alpha$ can be
regarded  either as a $*$-representation of $A$ of the injective von Neumann algebra $A$ on $\Hc$ or 
as a completely
positive mapping from $A$ into $\Bc(\Hc)$. In fact, it seems natural at first glance to take 
into consideration the second option, since \eqref{arveson} induces the correspondence 
$u\rho u^{-1}\mapsto u\alpha(u^{-1}\cdot u)u^{-1}$, ($u\in\G_A$); that is, the canonical map 
$\Sg(\rho)\to\Sc(\alpha)$. 
Nevertheless, in such a case we cannot be sure that the algebraic 
isomorphism $\Sc(\alpha)\simeq\G_A/\G(\alpha)$ entails a structure of smooth homogeneous manifold. 
Namely, 
the Lie algebra of $\G(\alpha)$ is $\gg(\alpha)=\{X\in A\mid X-\alpha(X)\in\alpha(A)'\}$, see 
Lemma~\ref{differen}, and it is not clear that $\gg(\alpha)$ is topologically complemented in $A$.

On the other hand, by looking at $\alpha$ as a $*$-representation 
and by using a bit of the structure of von Neumann algebras, 
it is possible to relate the orbits $\Sg(\rho)$ and $\Sg(\alpha)$ to each other,  
as we are going to see.

Since $A$ is Connes-amenable (that is, $A$ is injective, see \cite{Ru02}) an appropriate virtual diagonal 
$M$ for $A$ can be fixed, so that the mapping 
$$
E_\alpha(T):=\int_{A\otimes A}\alpha(a)T\alpha(b)\ dM(a,b), \quad (T\in A),
$$
is a conditional expectation from $A$ onto $B_\alpha:=E_\alpha(A)$. 
This endows $\Sg(\alpha)$ 
with the corresponding smooth homogeneous manifold structure. 

Similarly, $\Oc_N$ is a nuclear (amenable) $C^*$-algebra, see \cite{Cu77}, and therefore we can 
fix a contractive virtual diagonal $M_N$ for $\Oc_N$ so that the mapping 
$$
E_\rho(T):=\int_{\Oc_N\otimes \Oc_N}\rho(s)T\rho(t)\ dM_N(s,t), \quad (T\in A),
$$
is a conditional expectation from $A$ onto $B_\rho:=E_\rho(A)$. 
In this way, 
we regard $\alpha\colon A\to\Bc(\Hc)$ and $\rho\colon\Oc_N\to A$ as 
special cases of Example~\ref{amenable}.

There are two $C^*$-subalgebras of $A$ which are important in the study of the endomorphism $\alpha$. 
These are $\{a\in A\mid\alpha(a)=a\}$ and $\bigcap_{n>0}\alpha^n(A)$.
Recall that 
$B_\alpha=\{a\in A\mid E_\alpha(a)=a\}=\alpha(A)'$ and 
$B_\rho=\{a\in A\mid E_\rho(a)=a\}=\rho(\Oc_N)'$, 
see \cite{CG98}. 
(In this notation, $\gg(\alpha)=(I-\alpha)^{-1}( E_\alpha(A))$.) 
By Proposition~3.1(i) in \cite{Lac93a} 
we have  
$B_\rho=\rho(\Oc_N)'=\{a\in\Bc(\Hc)\mid \alpha(a)=a\}$ (which means in particular that 
$B_\rho={\mathbb C}\1$ if the representation $\rho$ is irreducible 
or, equivalently, if $\alpha$ is an \emph{ergodic} endomorphism of $A=\Bc(\Hc)$). 
Hence, 
for every $b\in B_\rho$ we get $b=\alpha(b)\in E_\alpha(A)'$. 
Shortly,
\begin{equation}\label{bicomm}
B_\rho\subseteq\alpha(A)''=\alpha(A).
\end{equation}
Note that $\alpha(A)'$ is $w^*$-closed in $A$ and so it is a von Neumann subalgebra of $A$. 
Moreover $\alpha(A)'$ is the range of the conditional expectation $E_\alpha$ and then it turns out to be 
injective. Let $\Delta$ be the anti-unitary operator ($\Delta$ is an antilinear isometry with $\Delta^2=\id$) 
which appears in the standard form of $\alpha(A)'$, see \cite{Ta03}, volume II. 
Then the mapping 
$E_{\Delta\alpha\Delta}\colon A\to A$ given by 
$$
E_{\Delta\alpha\Delta}(T):=\Delta E_\alpha(\Delta T\Delta)\Delta\quad (T\in A)
$$
is a conditional expectation such that $\alpha(A)''=E_{\Delta\alpha\Delta}(A)$. 
It follows that 
$\alpha(A)''$ is also injective, see \cite{Ta03}, volume III. 
The above notation is not by chance since 
$E_{\Delta\alpha\Delta}$ corresponds exactly to the expectation defined by the virtual diagonal $M$ 
and representation $\Delta\alpha\Delta$. 
Thus we have
\begin{equation}\label{berhoalpha}  
B_\rho\subseteq B_{\Delta\alpha\Delta}. 
\end{equation}
Passing to quotients, \eqref{berhoalpha} implies that we obtain a canonical surjection 
$\Sg(\rho)\to\Sg(\Delta\alpha\Delta)$. 
In fact, it is a submersion mapping: 
The above map is 
$\G_A$-equivariant and its tangent map at $\rho$ is implemented by the bounded projection
$$
\id-E_{\Delta\alpha\Delta}\colon(\id-E_\rho)(a)\mapsto(\id-E_{\Delta\alpha\Delta})(a).
$$
from $A/B_\rho$ onto $A/B_{\Delta\alpha\Delta}$. 

On the other hand, the involutive transformation $\beta\mapsto\Delta\beta(\cdot)\Delta$ 
carries diffeomorphicaly similarity (unitary) orbits onto similarity (unitary) orbits 
of the space of representations $\Rep^\omega(A)$. 
In particular 
$\Sg(\alpha)$ and $\Sg(\Delta\alpha\Delta)$ are diffeomorphic through the map
$$
a\alpha a^{-1}\mapsto(\Delta a\Delta)(\Delta\alpha\Delta)(\Delta a^{-1}\Delta),\ 
 \Sg(\alpha)\to\Sg(\Delta\alpha\Delta)
$$
(note that $u\alpha=\alpha u$ if and only if 
$(\Delta u\Delta)(\Delta\alpha\Delta)=(\Delta\alpha\Delta)(\Delta u\Delta)$). 
Putting all the above facts together we obtain the analytic submersion 
$\Sg(\rho)\to\Sg(\alpha)$ given by 
$$
a\rho a^{-1}\mapsto(\Delta a\Delta)\alpha(\Delta a^{-1}\Delta),\quad (a\in\G_A).
$$
At this point one can return to \eqref{bicomm}, which also tells us that 
\begin{equation}\label{bealpharho}  
B_\alpha=\alpha(A)'\subseteq B_\rho'.
\end{equation}
Thus we can proceed as formerly, just replacing $\alpha$ with $\rho$. 
In particular, 
the von Neumann algebra $B_\rho=\rho(\Oc_N)'$ is injective and then its commutant subalgebra 
$B_\rho'$ in $A$ is also injective. 
In effect, it is the image $B_\rho'=E_{R\rho R}(A)$ of the 
conditional expectation $E_{R\rho R}$ defined by the virtual diagonal $M_N$ and representation 
$R\rho R$, where $R$ is the anti-unitary operator involved in the standard form of $\rho(\Oc_N)'$. 
Hence, there exists an analytic submersion $\Sg(\alpha)\to\Sg(\rho)$ given by 
$$
u\rho u^{-1}\mapsto(R u R)\rho(R u^{-1} R),\quad (u\in\G_A).
$$
In summary, we have proved the following result.

\begin{theorem}\label{bisubmer} 
Let $\Delta$, $R$ be the anti-unitary operators associated 
with the standard
forms of the injective von Neumann algebras $\alpha(A)'$, $\rho(\Oc_N)'$ respectively. 
Then $B_\rho\subseteq B_{\Delta\alpha\Delta}$ and $B_\alpha\subseteq B_{R\rho R}$, and the mappings
$$\begin{aligned} 
q_\rho\colon 
 & a\rho a^{-1}\mapsto(\Delta a\Delta)\alpha(\Delta a^{-1}\Delta),\ \Sg(\rho)\to\Sg(\alpha), \\
q_\alpha\colon 
 & u\alpha u^{-1}\mapsto(R u R)\rho(R u^{-1} R),\ \Sg(\alpha)\to\Sg(\rho)
 \end{aligned}
$$
are analytic submersions.
\end{theorem}

\begin{remark}\label{nodifeo}
\normalfont
The joint action of suitably chosen mappings in the proposition yields new submersions 
$$ 
a\rho a^{-1}\mapsto(V a V^{-1})\rho(V a^{-1} V^{-1}),\quad \Sg(\rho)\to\Sg(\rho)
$$
and
$$
u\alpha u^{-1}\mapsto (V^{-1} u V)\alpha(V^{-1} u^{-1} V),\quad \Sg(\alpha)\to\Sg(\alpha),
$$
where $V$ is the unitary operator $V=R\Delta$. 
Such submersions need not be diffeomorphisms. 
\qed
\end{remark} 

Because of the inclusion $B_\rho\subseteq B_{\Delta\alpha\Delta}$ we have 
$E_{\Delta\alpha\Delta}\circ E_\rho=E_\rho$. Set $F_\rho:=E_\rho E_{\Delta\alpha\Delta}$. 
Then $F_\rho$ is also a conditional expectation and  
$F_\rho$ and $E_\rho$ are equivalent:
$F_\rho E_\rho=E_\rho$ and $E_\rho F_\rho=F_\rho$, so that $F_\rho(A)=B_\rho$. 
In addition, $F_\rho$ and $E_{\Delta\alpha\Delta}$ commute: 
$E_{\Delta\alpha\Delta} F_\rho=F_\rho=F_\rho E_{\Delta\alpha\Delta}$. 

Let $\Phi\colon A\to\Bc(\Hc_0)$ be a completely positive map, 
for some Hilbert space  $\Hc_0$. Put $\Phi_\rho:=\Phi\circ F_\rho$. 
Then $\Phi_\rho F_\rho=\Phi_\rho$ and $\Phi_\rho E_{\Delta\alpha\Delta}=\Phi_\rho$. 
Applying Stinepring's dilation theorem we find Hilbert spaces  
$\Hc_J(\Phi_\rho)$, for $J=A,B_{\Delta\alpha\Delta},B_\rho$ with 
$\Hc_{B_\rho}(\Phi_\rho)\subseteq\Hc_{B_{\Delta\alpha\Delta}}(\Phi_\rho)\subseteq\Hc_A(\Phi_\rho)$, 
and $*$-representations 
$\pi_J\colon J\to\Bc(\Hc_J(\Phi_\rho))$ 
satisfying $\pi_{B_{\Delta\alpha\Delta}}(u)
=\pi_A(u)|_{\Hc_{B_{\Delta\alpha\Delta}}(\Phi_\rho)}$ for each $u\in B_{\Delta\alpha\Delta}$, 
and 
$\pi_{B_\rho}(u)=\pi_{B_{\Delta\alpha\Delta}}(u)|_{\Hc_{B_{\rho}}(\Phi_{\rho})}$ 
for each $u\in B_\rho$. 

\begin{corollary}\label{univcofin}
In the above setting, the following commutative diagram and $\G_A$-equivariant
diagram: 
$$
\begin{CD}
\G_A\times_{\G_{B_\rho}}\Hc_{B_{\rho}}(\Phi_{\rho}) @>{ }>> 
\G_A\times_{\G_{B_\alpha}}\Hc_{B}(\Phi_{\rho}) \\
@V{ }VV @VV{ }V \\
\Sg(\rho) @>{q_\rho}>> \Sg(\alpha)
\end{CD}
$$
whose arrows are $\G_A$-equivariant and compatible with the 
involutive diffeomorphisms $-*$ on both $\Sg(\rho)$ and $\Sg(\alpha)$.
Moreover, the representation $\pi_A$ of $\G_A$ on $\Hc_A(\Phi_{\rho})$ 
can be extended to $A$ and realized as multiplication
on a reproducing kernel Hilbert space formed by holomorphic sections of 
the left-side vector bundle in the diagram.
\end{corollary}

\begin{proof}
Firstly, a diagram similar to that one of the statement, and concerning the algebras 
$B_\rho\subseteq B_{\Delta\alpha\Delta}$, is immediately obtained by mimicing 
the proof of Corollary~\ref{sandwich}. Then using the diffeomorphism $\Delta(\cdot)\Delta$ one gets the wanted
result  or diagram, where the action of $B_\alpha$ on $\Hc_{B}(\Phi_{\rho})$ is 
given just by transferring the action of $B_{\Delta\alpha\Delta}$ through $\Delta(\cdot)\Delta$ (note that 
$B_{\Delta\alpha\Delta}=\Delta B_\alpha\Delta$).
\end{proof}

A diagram entirely analog to the previous one of the corollary can be obtained 
by interchanging roles of representations $\rho$ and $\alpha$.
 
\begin{remark}\label{hooking}
\normalfont
Using the representation $\rho$, we can make a link between 
Corollary~\ref{fourier} and the preceding setting. 
Let $\tau$ be the gauge automorphism group of $\Oc_N$, and let $E_\tau$ be 
the expectation defined
by \eqref{note2} for $m=0$. 
Corollary~\ref{fourier} applies to completely positive mappings 
$\Phi\colon\Oc_N\to\Bc(\Hc_0)$ such that 
$\Phi\circ E_\tau=\Phi$. Assume that $\rho$ is a $*$-representation of 
$\Oc_N$ on a von Neumann algebra $A$, 
and $E_\rho\colon A\to A$ the conditional expectation associated with some, 
fixed, virtual diagonal $M$ of norm one 
for $A$. Let $\Phi\colon E_\rho(A)\to\Bc(\Hc_0)$ be completely positive, 
and let consider 
$\Phi_{\rho,\tau}:=\Phi\circ E_\rho\circ\rho\circ E_\tau$. 
Then $\Phi_{\rho,\tau}E_\tau=\Phi_{\rho,\tau}$ and 
we obtain
Hilbert spaces and their decompositions like those of 
Corollary~\ref{fourier}, associated with the representation $\rho$ 
and algebra $A$. 
\qed
\end{remark}

Finally, the subalgebra $\rho(\Fc_N)'=\bigcap_{n>0}\alpha^n(A)$ 
(see Proposition~3.1(ii) in \cite{Lac93a})
suggested us to form sequences of vector bundles in the following manner.
Let $\alpha\colon A\to A$ be a normal, $*$-representation where $A=\Bc(\Hc)$ 
as above. 
Then $\alpha^*(A_*)\subseteq A_*$ where $A_*$ denotes the predual of $A$ 
formed by the trace-class operators on $\Hc$, 
and $\alpha^*$ is the transpose mapping of $\alpha$. For $n\in\N$ we are going 
to consider the iterative mappings 
$\beta_n:=\alpha^n\circ\rho$ and corresponding expectations denoted by 
$E_n:=E_{\beta_{n}}$ and put $E_0=E_\rho$. Then, for $\xi\in A_*$ and 
$T\in A$,
$$ 
\begin{aligned}
(E_n\circ\alpha)(T)(\xi)&=\int_{\Oc_N\otimes \Oc_N}(\alpha^n\rho)(s)\ \alpha(T)\ (\alpha^n\rho)(t)(\xi)\ dM_N(s,t)
\\ 
&=\int_{\Oc_N\otimes \Oc_N}(\alpha^{n-1}\rho)(s)\ T\ (\alpha^{n-1}\rho)(t)\ (\alpha^*\xi)\ dM_N(s,t) \\
& =E_{n-1}(T)(\alpha^*\xi)=(\alpha\circ E_{n-1})(T)(\xi),
\end{aligned}
$$
see \cite{CG98}. 
More specifically, $(E_\rho\alpha)(T)
=\alpha(\int_{\Oc_N\otimes \Oc_N} \rho(s)\ T\ \rho(t)\ dM_N(s,t))=\varphi(T)\in\C$, where 
$\varphi$ is the state given by $\varphi(T)
:=\int_{\Oc_N\otimes \Oc_N} \rho(s)\ T\ \rho(t)\ dM_N(s,t)\in\Bc(\Hc)'={\mathbb C}\1$, $T\in A$. 
Hence
\begin{equation}\label{intertwin1}  
E_n\circ\alpha=\alpha\circ E_{n-1},\ \qquad {n\in\N}
\end{equation}
whence, by a reiterative process and since $\alpha E_\rho=E_\rho$, we get
\begin{equation}\label{finalrho}  
E_n\circ\alpha^n= E_\rho,\ \qquad {n\in\N}
\end{equation}
Hence we get $E_nE_\rho=E_\rho$ and therefore $B_\rho\subseteq B_n$, 
where $B_n:=E_n(A)$, for all $n$. 
Further, we have
$\alpha E_n(A)=E_{n+1}\alpha(A)\subseteq E_{n+1}(A)$ by \eqref{intertwin1}, 
that is, $\alpha(B_n)\subseteq B_{n+1}$, $n\in\N$. 

Now consider a countable family $(\Phi_n)_{n\ge0}$
of completely positive mappings $\Phi_n\colon A\to\Bc(\Hc_0)$, for some Hilbert space $\Hc_0$, 
such that
\begin{equation}\label{compatibi}
\Phi_{n+1}\circ \alpha=\Phi_n\ , \qquad \Phi_n\circ E_n=\Phi_n,\ \qquad {n\in\N}
\end{equation}
Such a family exists. 
Take for instance $\phi_n:=E_\rho E_1\cdots E_n$, 
and a completely positive map $\Phi\colon A\to\Bc(\Hc_0)$. 
Then the family 
$\Phi_n:=\Phi\circ\phi_n$, $n\ge0$, satisfies \eqref{compatibi}. 
In these conditions the diagram 
$$
\begin{CD}
A @>{\alpha}>> A @>{\alpha}>> \cdots @>{\alpha}>> A @>{\alpha}>>  A @>{\Phi_{n+1}}>> \Bc(\Hc_0)\\ 
@VV{E_0}V @VV{E_1}V 
@VVV @VV{E_n}V  @VV{E_{n+1}}V @VV{\id}V \\ 
B_0 @>{\alpha}>> B_1 @>{\alpha}>> \cdots @>{\alpha}>> B_n @>{\alpha}>> B_{n+1} @>{\Phi_{n+1}}>> \Bc(\Hc_0)
\end{CD}
$$
is commutative in each of its subdiagrams.

For every $n\ge0$, by applying Theorem~\ref{stinerealiz} to the conditional expectation 
$E_n\colon A\to B_n$ and mapping $\Phi_n$ one finds the corresponding Hilbert space 
$\Hc_{B_n}(\Phi_n)$ 
for the representation which is the Stinespring dilation of $\Phi_n$. Take a finite set of 
elements $b_j$ in $B_n$. As 
$\Phi_n(b_i^*b_j)=\Phi_{n+1}(\alpha(b_i)^*\alpha(b_j))$ it follows that   
$$
\Vert\sum_{j}b_j\otimes f_j\Vert_{\Phi_n}=\Vert\sum_{j}\alpha(b_j)\otimes f_j\Vert_{\Phi_{n+1}}
$$
for all 
$\sum_{j}b_j\otimes f_j\in B_n\otimes\Hc_0$, see Section~\ref{sect5}. 
Hence, 
$\alpha(\Hc_{B_n}(\Phi_n))\subseteq\Hc_{B_{n+1}}(\Phi_{n+1})$. 
This implies that we have found 
the (countable) system of vector bundle homomorphisms 
$$
\begin{CD}
\G_A\times_{\G_{B_\rho}}\Hc_{B_\rho}(\Phi_0) @ >{\id\otimes\alpha} >> 
\G_A\times_{\G_{B_1}}\Hc_{B_1}(\Phi_1) @ >{\id\otimes\alpha} >>
\cdots @>{\id\otimes\alpha}>>  
\G_A\times_{\G_{B_n}}\Hc_{B_n}(\Phi_n) @> {\id\otimes\alpha} >> 
\cdots \\  
@V{}VV @VV{}V @VV{ }V @VV{}V \\  
\Sg(\rho) @ > {\tilde{\alpha}_0} >> 
\Sg(\alpha^1\rho) @ > {\tilde{\alpha}_1} >> 
\cdots @ > {\tilde{\alpha}_{n-1}} >> 
\Sg(\alpha^n\rho) @ > {\tilde{\alpha}_n} >> \cdots
\end{CD}
$$
where $\tilde{\alpha}_n$ is the canonical submersion induced by 
$\alpha\vert_{B_n}\colon B_n\to B_{n+1}$, $n\ge0$.

Of course the above sequence of diagrams gives rise to the corresponding 
statements about complexifications, and realizations of representations on spaces of 
holomorphic sections.

\bigskip
5) {\bf Non-commutative stochastic analysis}

We have just shown a sample of how to find sequences of homogeneous vector bundles 
of the type dealt with in this paper. As a matter of fact, continuous families of 
such bundles are also available, which could hopefully be of interest in other fields. More precisely,     
the geometric models developed in the present paper 
might prove useful in order to get a better understanding 
of the phenomena described by 
the various theories of non-commutative probabilities. 
By way of illustrating this remark, we shall briefly discuss 
from our geometric perspective 
a few basic ideas related to the stochastic calculus 
on full Fock spaces as developed in \cite{BV00} and \cite{BV02}. 
(See also \cite{VDN92} and \cite{Ev80} for a complementary perspective 
that highlights the role of the Cuntz algebras 
in connection with full Fock spaces.) 

In the paper \cite{BV00}, 
a family of conditional expectations $\{E_t\}_{t>0}$ is 
built on the  von Neumann algebra $A$ 
of bounded operators on
the full Fock space, generated by the annihilation, creation, 
and gauge operators. 
Set $A_t:=E_t(A)$ for $t>0$. 
It is readily seen that $A_t\subseteq A_s$ and 
that $E_tE_s=E_t$ whenever 
$0<s\le t$ (check first for the so-called in \cite{BV00} {\it basic} elements). 
Applying the Stinespring dilation procedure to 
the conditional expectation $E_s$ and completely positive mapping $E_t$ 
one gets Hilbert spaces 
$\Hc_{A_s}(E_t)\subseteq \Hc_A(E_t)$ and 
the consequent Stinespring representations
$\pi_{A_j}\colon A_j\to \Bc(\Hc_{A_j}(E_t))$, where $j=0,t$, and $A_0=A$. 
This entails the commutative diagram  
$$
\begin{CD}
\G_A\times_{\G_t}\Hc_{A_t}(E_t) @>{ }>> 
\G_A\times_{\G_s}\Hc_{A_s}(E_t) @>{ }>>
\G_A\times_{\G_r}\Hc_{A_r}(E_t) \\ 
@V{ }VV @VV{ }V @VV{ }V \\
\G_A/\G_t @>{ }>> \G_A/\G_s @>{ }>> \G_A/\G_r,
\end{CD}
$$
for $r<s<t$, where $\G_j=\G_{A_j}$ for $j=r,s,t$. 
Moreover, as usual, the geometrical framework of the
present paper works to produce a Hilbert space $\Hc_A(E_s,E_t)$, 
formed by holomorphic sections on
$\G_A/\G_s$, which is isometric to $\Hc_A(E_t)$ and enables us 
to realize $\pi_A$ as natural 
multiplication.

On the other hand, from the point of view of 
the quantum stochastic analysis 
(see for instance \cite{Par90} and \cite{BP95}), 
it is worth considering unital completely positive mappings  
$\Phi\colon A\to\Bc(\Hc_0)$ with the following filtration property: 
There exists a family $\{\Phi_t\colon A\to\Bc(\Hc_t)\}_{t\ge0}$ of 
completely positive mappings which approximate $\Phi$ in some sense 
and satisfy 
$\Phi_t\circ E_t=\Phi_t$ for all $t>0$. 
Then we get commutative diagrams 
$$
\begin{CD}
\G_A\times_{\G_t}\Hc_{A_t}(\Phi_t) @>{ }>> 
  \G_A\times_{\G_s}\Hc_{A_s}(\Phi_s) \\
@V{ }VV @VV{ }V \\
\G_A/\G_t @>{ }>> \G_A/\G_s
\end{CD}
$$
whenever $s<t$.
By means of the realizations of the full Fock space 
as reproducing kernel Hilbert spaces of sections in 
appropriate holomorhic vector bundles 
we  find geometric interpretations 
for most concepts tradtionally related to the Fock spaces 
(for instance, annihilation, creation, and gauge operators). 
We thus arrive at the challenging perspective of 
a relationship between the non-commutative stochastic analysis  
and the infinite-dimensional complex geometry, 
which certainly deserves to be understood in more detail.  
For one thing, this might provide useful geometric insights in 
areas like the theory of quantum Markov processes. 

\bigbreak

\noindent\textbf{Acknowledgment.} 
We wish to thank Professor Karl-Hermann Neeb for pointing out to us 
some pertinent references. 
The present research was begun during the first-named author's visit 
at the Department of Mathematics of the University of Zaragoza 
with financial support from Proyecto MTM2004-03036, DGI-FEDER, of the MCYT, Spain. 
Our research was further developed during the visit of the second-named author's visit 
at the Institute of Mathematics "Simion Stoilow" of the Romanian Academy 
with support from the aforementioned grant and from the Softwin Group. 
Partial financial support from 
Proyecto E-64 of the DG Arag\'on, Spain, 
and from 
Grant 2-CEx06-11-22/2006 of the Romanian Government  
is also acknowledged.


\begin{thebibliography}{nunitary}

\small

\bibitem[AMR88]{AMR88}
R.~Abraham, J.E.~Marsden, T.~Ratiu, 
\textit{Manifolds, Tensor Analysis, and Applications}
(second edition), Applied Mathematical Sciences 75,
Springer-Verlag, New York, 1988.

\bibitem[AS94]{AS94}
E.~Andruchow, D.~Stojanoff, 
Geometry of conditional expectations and finite index, 
\textit{Int. J. Math.}
\textbf{5} (1994), no.~2, 169--178.

\bibitem[ACS95]{ACS95}
E.~Andruchow, G.~Corach, D.~Stojanoff, 
A geometric characterization of nuclearity and injectivity, 
\textit{J. Funct. Anal.} 
\textbf{133} (1995), 474-494.

\bibitem[ALRS97]{ALRS97}
E.~Andruchow, A.~Larotonda, L.~Recht, D.~Stojanoff, 
Infinite dimensional homogeneous reductive spaces and finite index 
conditional expectations, 
\textit{Illinois J. Math.}
\textbf{41} (1997), 54-76.

\bibitem[AS01]{AS01}
M.~Argerami, D.~Stojanoff, 
Orbits of conditional expectations, 
\textit{Illinois J. Math.}
\textbf{45} (2001), no.~1, 243--263.

\bibitem[Ar89]{Ar89}
W.B.~Arveson, 
Continuous analogues of Fock space, 
\textit{Memoirs Amer. Math. Soc.} 
\textbf{80} (1989), no.~409. 

\bibitem[Ar00]{Ar00}
W.B.~Arveson, 
The curvature invariant of a Hilbert module over 
${\mathbb C}[z_1,\dots,z_d]$, 
\textit{J. reine angew. Math} \textbf{522} (2000), 173--236. 

\bibitem[BV00]{BV00}
C.~Barnett, A.~Voliotis, 
A conditional expectation for the full Fock space,
\textit{J. Oper. Th.} {\bf 44} (2000), 3-23. 

\bibitem[BV02]{BV02}
C.~Barnett, A.~Voliotis, 
Stochastic and belated integrals for the full Fock space, 
\textit{Soochow J. Math.} \textbf{28} (2002), no.~4, 413--441.

\bibitem[Be06]{Be06}
D.~Belti\c t\u a, 
\textit{Smooth Homogeneous Structures in Operator Theory}, 
Monogr. and Surveys in Pure and Applied Math., 137. 
Chapman \& Hall/CRC Press, 
Boca Raton-London-New York-Singapore, 2006.

\bibitem[BR07]{BR07}
D.~Belti\c t\u a, T.S.~Ratiu, 
Geometric representation theory for unitary groups of operator algebras, 
\textit{Adv. Math.} \textbf{208} (2007), no.~1, 299--317.

\bibitem[BN04]{BN04} 
W.~Bertram, K.-H.~Neeb, 
Projective completions of Jordan pairs, 
Part I. The generalized projective geometry of a Lie algebra, 
\textit{J. Algebra} \textbf{277} (2004), 474--519. 

\bibitem[BN05]{BN05} 
W.~Bertram, K.-H.~Neeb, 
Projective completions of Jordan pairs, 
Part II. Manifold structures and symmetric spaces, 
\textit{Geom. Dedicata} \textbf{112} (2005), 73--113. 

\bibitem[BP95]{BP95} 
B.V.R.~Bhat, K.R.~Parthasarathy, 
Markov dilations of nonconservative dynamical semigroups and a quantum
boundary theory, 
\textit{Ann. Inst. H. Poincar\'e Probab. Statist.} 
\textbf{31} (1995), no.~4, 601--651. 

\bibitem[Bi03]{Bi03} 
R.~Bielawski, 
Complexification and hypercomplexification of manifolds with 
a linear connection, 
\textit{Internat. J. Math.} \textbf{14} (2003), no.~8, 813--824.

\bibitem[Bi04]{Bi04}
R.~Bielawski, 
Prescribing Ricci curvature on complexified symmetric spaces, 
\textit{Math. Res. Lett.} \textbf{11} (2004), no.~4, 435--441.

\bibitem[Bo57]{Bo57} 
R.~Bott, 
Homogeneous vector bundles, 
\textit{Ann. Math.} \textbf{66} (1957), no.~2, 203--248. 

\bibitem[Bo80]{Bo80} 
R.P.~Boyer, 
Representation theory of the Hilbert-Lie group $\U(\Hg)_2$, 
\textit{Duke Math. J.} \textbf{47} (1980), no.~2, 325--344. 

\bibitem[BJ99]{BJ99} 
O.~Bratteli, P.E.T.~Jorgensen, 
Iterated function systems and permutation representations of 
the Cuntz algebra, 
\textit{Mem. Amer. Math. Soc.} \textbf{139} (1999), no.~663. 

\bibitem[Bu81]{Bu81} 
J. W.~Bunce, 
The similarity problem for representations of operator algebras, 
\textit{Proc. Amer. Math. Soc.} {\bf 81} (1981), 409--414.

\bibitem[Ch81]{Ch81} 
E.~Christensen, 
On non selfadjoint representations of $C^*$-algebras, 
\textit{Amer. J. Math.} {\bf 103} (1981), 817--834.

\bibitem[CG98]{CG98}
G.~Corach, J.E.~Gal\'e, 
Averaging with virtual diagonals and geometry of representations,  
in: {\it Banach Algebras '97 (Blaubeuren)}, 
de Gruyter, Berlin, 1998, pp.~87--100. 

\bibitem[CG99]{CG99} 
G.~Corach, J.E.~Gal\'e, 
On amenability and geometry of spaces of bounded representations, 
\textit{J. London Math. Soc.} {\bf 59} (1999), 
no.~2, 311--329.

\bibitem[Cu77]{Cu77}
J.~Cuntz, 
Simple $C^*$-algebras generated by isometries,
\textit{Comm. Math. Phys.} {\bf 57} (1977), no.~2, 173-185.

\bibitem[CK80]{CK80}
J.~Cuntz, W.~Krieger, 
A class of C*-algebras and topological Markov chains, 
\textit{Invent. Math.} \textbf{56}(1980), 251-268.

\bibitem[Da96]{Dav96}
K.R.~Davidson, 
\textit{$C\sp *$-algebras by Example}, 
 Fields Institute Monographs, 6. American Mathematical Society,
Providence, RI, 1996.

\bibitem[Di64]{Di64}
J.~Dixmier, 
\textit{Les $C^*$-alg\`ebres et Leurs Representations}, 
Gauthier-Villars, Paris, 1964. 

\bibitem[DPZ98]{DPZ98}
S.~Doplicher, C.~Pinzari, R.~Zuccante, 
The $C^*$-algebra of a Hilbert bimodule,
\textit{Bollettino UMI. Serie VIII} \textbf{1 B} (1998), 263--282.

\bibitem[Do66]{Do66} 
A.~Douady, 
Le probl\`eme des modules pour les sous-espaces analytiques compacts 
d'un espace analytique donn\'e, 
\textit{Ann. Inst. Fourier (Grenoble)} 
\textbf{16} (1966), no.~1, 1--95.

\bibitem[DG01]{DG01}
M.J.~Dupr\'e, J.F.~Glazebrook, 
The Stiefel bundle of a Banach algebra, 
\textit{Integral Equations Operator Theory} 
\textbf{41} (2001), no.~3, 264--287. 

\bibitem[DG02]{DG02}
M.J.~Dupr\'e, J.F.~Glazebrook, 
Holomorphic framings for projections in a Banach algebra, 
\textit{Georgian Math. J.} 
\textbf{9} (2002), no.~3, 481--494.

\bibitem[ER00]{ER00}
E.G.~Effros, Zh.-J.~Ruan, 
{\it Operator Spaces}, 
London Mathematical Society Monographs. New Series, 23. 
The Clarendon Press, Oxford University Press, New York, 2000.

\bibitem[Ev80]{Ev80}
D.E.~Evans, 
On $\Oc_n$,
\textit{Publ. Res. Inst. Math. Sci.} {\bf 16} (1980), no.~3, 915-927.


\bibitem[Fo95]{Fo95}
G. B.~Folland, 
{\it  A Course in Abstract Harmonic Analysis}, 
Studies in Advances Mathematics. CCR Press Inc., Boca Raton, 1995. 


\bibitem[Ga06]{Ga06}
J.E.~Gal\'e, 
Geometr\'\i a de \'orbitas de representaciones de grupos y \'algebras  
promediables, 
\textit{Rev. R. Acad. Cienc. Exactas F\'\i s. Qu\'\i m. Nat. Zaragoza (2)} 
\textbf{61} (2006), 7--46. 


\bibitem[GW54]{GW54}
L.~Garding, A.~Wightman, 
Representations of the anticommutation relations,  
\textit{Proc. Nat. Acad. Sci., USA} \textbf{40} (1954), 617--621. 
 

\bibitem[GN03]{GN03}
H. Gl\"ockner, K.-H. Neeb,
Banach-Lie quotients, enlargibility, and universal complexifications, 
{\it J. Reine Angew. Math.}  
{\bf 560}(2003), 1--28.

\bibitem[Ha83]{Ha83}
U.~Haagerup, 
Solution of the similarity problem for cyclic representations of 
$C^*$- algebras, 
\textit{Ann. Math. (2)} \textbf{118} (1983), 215--240.


\bibitem[Ki04]{Ki04}
A. A.~Kirillov, 
{\it  Lectures on the Orbit Method}, 
Graduate Studies in Mathematics, 64. American Mathematical Society, 
Providence, Rhode Island, 2004.




\bibitem[IS66]{IS66}
N.~Iwahori, M.~Sugiura, 
A duality theorem for homogeneous manifolds of compact Lie groups,  
\textit{Osaka J. Math.} \textbf{3} (1966), 139--153. 

\bibitem[La93a]{Lac93a} 
M.~Laca, 
Endomorphisms of $\Bc(\Hc)$ and Cuntz algebras, 
\textit{J. Operator Theory} \textbf{30} (1993), no.~1, 85--108.

\bibitem[La93b]{Lac93} 
M. Laca, 
Gauge invariant states of ${\mathcal O}_\infty$,  
\textit{J. Operator Theory} \textbf{30} (1993), no.~2, 381--396.

\bibitem[Ln01]{La01}
S.~Lang,
{\it Fundamentals of Differential Geometry} (corrected second printing),
Graduate Texts in Mathematics, 191. Springer-Verlag,
New-York, 2001.

\bibitem[Ls78]{Las78}
M.~Lassalle, 
S\'eries de Laurent des fonctions holomorphes dans la complexification 
d'un espace sym\'etrique compact, 
\textit{Ann. Scient. \'Ec. Norm. Sup.} \textbf{11}(1978), 167--210. 

\bibitem[LS91]{LS91}
L.~Lempert, R.~Sz\H oke, 
Global solutions of the homogeneous complex Monge-Amp\`ere equation 
and complex structures on the tangent bundle of Riemannian manifolds,  
\textit{Math. Ann.} \textbf{290} (1991), no.~4, 689--712.

\bibitem[MR92]{MR92}
L.E.~Mata-Lorenzo, L.~Recht, 
Infinite-dimensional homogeneous reductive spaces,  
{\it Acta Cient. Venezolana}  
{\bf 43}(1992),  no.~2, 76--90. 

\bibitem[Mo55]{Mo55}
G.D.~Mostow, 
On covariant fiberings of Klein spaces, 
\textit{Amer. J. Math.} \textbf{77} (1955), 247--278.

\bibitem[Mo05]{Mo05} 
G.D.~Mostow, 
A structure theorem for homogeneous spaces, 
\textit{Geom. Dedicata} \textbf{114} (2005), 87--102. 

\bibitem[MS03]{MS03}
P.S.~Muhly, B.~Solel,  
The curvature and index of completely positive maps,
\textit{Proc. London Math. Soc.} {\bf 87} (2003), no.~3, 748-778. 

\bibitem[Ne00]{Ne00}
K.-H.~Neeb,
\textit{Holomorphy and Convexity in Lie Theory}, 
de Gruyter Expositions in Mathematics 28, 
Walter de Gruyter \& Co., 
Berlin, 2000.

\bibitem[Ne02]{Ne02}
K.-H.~Neeb,
A Cartan-Hadamard theorem for Banach-Finsler manifolds,
\textit{Geom. Dedicata} \textbf{95} (2002), 115--156.

\bibitem[Ne04]{Ne04}
K.-H.~Neeb, 
Infinite-dimensional groups and their representations,
in: \textit{Lie Theory},
Progr. Math., 228, Birk\-h\"auser, Boston, MA, 2004,
pp.~213--328.

\bibitem[OR03]{OR03}
A.~Odzijewicz, T.S.~Ratiu,
Banach Lie-Poisson spaces and reduction,
\textit{Commun. Math. Phys.} \textbf{243} (2003), 1--54.
{\bf 95}(2002), 115--156.

\bibitem[On60]{On60}
A.L.~Oni\v s\v cik, 
Complex hulls of compact homogeneous spaces, 
\textit{Dokl. Akad. Nauk SSSR} \textbf{130}, 726--729 (Russian); 
translated as
\textit{Soviet Math. Dokl.} \textbf{1} (1960), 88--91. 

\bibitem[Par90]{Par90}
K.R.~Parthasarathy, 
A continuous time version of Stinespring's theorem 
on completely positive maps, 
in: \textit{Quantum probability and applications, V (Heidelberg, 1988)}, 
Lecture Notes in Math. 1442, Springer, Berlin, 1990, pp.~296--300.

\bibitem[Pau02]{Pa02}
V.~Paulsen, 
{\it Completely Bounded Maps and Operator Algebras}, 
Cambridge Studies in Advanced Mathematics, 78. 
Cambridge University Press, Cambridge, 2002.

\bibitem[Pi97]{Pi97}
M.~Pimsner, 
A class of C*-algebras generalizing both Cuntz-Krieger algebras
and cross product by ${\mathbb Z}$, 
in: D.~Voiculescu (ed.) \textit{Free Probability Theory}, 
Fields Institute Communications 12 (1997), 189-212.

\bibitem[Po01]{Po01} 
G.~Popescu, 
Curvature invariant over free semigroup algebras, 
\textit{Adv. Math.} \textbf{158} (2001), 264--309. 

\bibitem[PR87]{PR87}
H.~Porta, L.~Recht, 
Spaces of projections in a Banach algebra,  
\textit{Acta Cient. Venezolana} \textbf{38} (1987), no.~4, 408--426.

\bibitem[PR94]{PR94}
H.~Porta, L.~Recht, 
Conditional expectations and operator decompositions, 
\textit{Ann. Global Anal. Geom.} \textbf{12} (1994),
no.~4, 335--339. 

\bibitem[Ra77]{Rae77}
I.~Raeburn, 
The relationship between a commutative Banach algebra and its
maximal ideal space,  {\it J. Functional Analysis} {\bf 25}
(1977), no.~4, 366--390.

\bibitem[Ri74]{Ri74} 
M.A.~Rieffel, 
Induced representations of $C^*$-algebras, 
\textit{Adv. Math.} \textbf{13} (1974), 176--257.

\bibitem[Ru02]{Ru02}
V.~Runde,
\textit{Lectures on Amenability}, 
Lecture Notes in Mathematics 1774, 
Springer, Berlin, 2002.

\bibitem[Sa71]{Sak71}
S.~Sakai, 
\textit{$C^*$-algebras and $W^*$-algebras}, 
Ergebnisse der Mathematik und ihrer Grenzgebiete,
Band 60, 
Springer-Verlag, New York-Heidelberg, 1971.

\bibitem[Sc64]{Sc64}
L. Schwartz,
Sous espaces Hilbertiens d'espaces vectoriels 
to\-po\-lo\-giques et noyaux associ\'es 
(noyeux reproduisants), 
{\it J. Analyse Math.} 
{\bf 13}(1964),  115--256.



\bibitem[Se57]{Se57}
I. E. Segal,
The structure of a class of representations of the unitary 
group on a Hilbert space,
{\it Proc. Amer. Math. Soc.} 
{\bf 8}(1957),  197--203.


\bibitem[Sh62]{Sh62}
D.~Shale,
Linear symmetries of the free boson field,
\textit{Trans. Amer. Math. Soc.} \textbf{103} (1962), 149--167.


\bibitem[St55]{St55}
W. F.~Stinespring, 
Positive functions on $C^*$-algebras, 
\textit{Proc. Amer. Math. Soc.} \textbf{6} (1955), 211--216.


\bibitem[SV02]{SV02}
S.~Str\u atil\u a, D.~Voiculescu
\textit{Representations of AF-algebras and of the group $\U(\infty)$}, 
Lecture Notes in Mathematics 486, 
Springer, Berlin, 1975.


\bibitem[Sz04]{Sz04} 
R.~Sz\H oke, 
Canonical complex structures associated to connections 
and complexifications of Lie groups, 
\textit{Math. Ann.} \textbf{329} (2004), no.~3, 553--591.


\bibitem[Ta03]{Ta03}
M.~Takesaki, 
\textit{Theory of operator algebras II, III}, 
Encyclopaedia of Mathematical Sciences,
vol. 125 and 127, 
Springer-Verlag, New York-Heidelberg, 2003.


\bibitem[To57]{To57}
J.~Tomiyama, 
On the projection of norm one in $W^*$-algebras, 
\textit{Proc. Japan Acad.} 
\textbf{33} (1957), 608--612.


\bibitem[Up85]{Up85}
H.~Upmeier,
{\it Symmetric Banach Manifolds and Jordan $C^*$-algebras},
North-Holland Mathematics Studies, 104. 
Notas de Matem\`atica, 96. North-Holland Publishing Co.,
Amsterdam,
1985.


\bibitem[VDN92]{VDN92}
D.V.~Voiculescu, K.J.~Dykema, A.~Nica,
\textit{Free Random Variables}  
(A noncommutative probability approach to 
free products with applications to random matrices, 
operator algebras and harmonic analysis on free
groups). 
CRM Monograph Series, 1. American Mathematical Society, 
Providence, RI, 1992.


\bibitem[We80]{We80}
R.O.~Wells, Jr.,  
{\it Differential Analysis on Complex Manifolds} (second edition),  
Graduate Texts in Mathematics, 65. Springer-Verlag, New York-Berlin, 1980.


\end{thebibliography}
\end{document}